\documentclass[11pt,a4paper,twoside,leqno]{article}

\usepackage{times}
\usepackage{amsmath,amssymb,amsthm,amsfonts,epsfig}
\usepackage{a4wide}
\usepackage{exscale}   
\usepackage{enumerate}
\usepackage{latexsym}

\usepackage[all]{xy} 

\usepackage{stmaryrd}

\newcommand{\restr}{\rvert}   
\newcommand{\abs}[1]{\left\lvert#1\right\rvert}   
\newcommand{\norm}[1]{\left\lVert#1\right\rVert}   

\newcommand{\demph}[1]{{\it #1}}

\DeclareMathOperator{\id}{Id}
\DeclareMathOperator{\ev}{ev}

\newcommand{\Cont}{{\mathcal C}} 
\DeclareMathOperator{\supp}{supp}
\DeclareMathOperator{\top2}{top}

\DeclareMathOperator{\Lin}{L} 
\DeclareMathOperator{\Leb}{L}

\DeclareMathOperator{\Komp}{K} 

\DeclareMathOperator{\KTh}{K}
\DeclareMathOperator{\KK}{KK}

\newcommand{\red}{\text{$r$}}
\DeclareMathOperator{\Cred}{C^*_{\red}}

\newcommand{\C}{\ensuremath{{\mathbb C}}}
\newcommand{\E}{\ensuremath{{\mathbb E}}}

\newcommand{\N}{\ensuremath{{\mathbb N}}}

\newcommand{\R}{\ensuremath{{\mathbb R}}}

\newcommand{\mA}{\ensuremath{{\mathcal A}}}
\newcommand{\mB}{\ensuremath{{\mathcal B}}}

\newcommand{\mD}{\ensuremath{{\mathcal D}}}

\newcommand{\mF}{\ensuremath{{\mathcal F}}}
\newcommand{\mG}{\ensuremath{{\mathcal G}}}
\newcommand{\mH}{\ensuremath{{\mathcal H}}}

\newcommand{\mN}{\ensuremath{{\mathcal N}}}

\newcommand{\mR}{\ensuremath{{\mathcal R}}}

\newcommand{\mU}{\ensuremath{{\mathcal U}}}

\newcommand{\fG}{\ensuremath{{\mathfrak G}}}

\theoremstyle{plain}

\newtheorem {theorem} {Theorem}[section]
\newtheorem {lemma}[theorem] {Lemma}
\newtheorem {proposition} [theorem]{Proposition}
\newtheorem {corollary} [theorem]{Corollary}

\newtheorem* {theorem*} {Theorem}
\newtheorem* {proposition*} {Proposition}
\newtheorem* {lemma*}{Lemma}

\theoremstyle{definition}

\newtheorem {definition} [theorem]{Definition}
\newtheorem {defprop} [theorem]{Definition and Proposition}

\newtheorem {examples}[theorem] {Examples}

\newtheorem {Xample}[theorem] {Example}

\DeclareMathOperator{\KKban}{KK^{\ban}}
\DeclareMathOperator{\KKbanG}{KK^{\ban}_G}

\DeclareMathOperator{\Eban}{\E^{\ban}}
\DeclareMathOperator{\EbanG}{\E^{\ban}_G}

\DeclareMathOperator{\ban}{ban}

\newcommand{\LazyAnd}{\quad \text{and} \quad}

\newcommand{\rmd}{{\, \mathrm{d}}}
\newcommand{\unital}[1]{\widetilde{#1}} 

\newcommand{\ketbra}[2]{\big|#1\big\rangle\big\langle #2 \big|} 

\newcommand{\lAngle}{\big\langle\!\big\langle}
\newcommand{\rAngle}{\big\rangle\!\big\rangle}

\DeclareMathOperator{\RKK}{\mR KK}

\newcommand{\RKKban}{\RKK^{\ban}} 
\newcommand{\RKKbanG}{\RKK^{\ban}_{G}} 

\DeclareMathOperator{\uEgG}{\underline{\rm E}\gG}

\newcommand{\gG}{\mG} 

\newcommand{\EbanW}[1]{\E^{\ban}_{#1}} 
\newcommand{\KKbanW}[1]{\KK^{\ban}_{#1}} 

\newcommand{\GreenJulgAbstiegEban}[3]{J_{#1,#2}^{#3}} 
\newcommand{\GreenJulgAbstiegKKban}[2]{J_{#1}^{#2}} 

\newcommand{\GreenJulgAufstiegEban}[3]{M_{#1,#2}^{#3}} 
\newcommand{\GreenJulgAufstiegKKban}[2]{M_{#1}^{#2}} 

\newcommand{\Gelfand}[1]{\fG\!\left(#1\right)} 
\newcommand{\ContSect}{\Gamma}

\newcommand{\Field}[1]{\FieldPur\left(#1\right)} 
\newcommand{\FieldPur}{\mathfrak{F}} 
\newcommand{\cA}{\mA}

\DeclareMathOperator{\MoritabanG}{Mor^{\ban}_G} 

\begin{document}

\title{A Generalised Green-Julg Theorem for\\ Proper Groupoids and Banach Algebras}
\author{Walther Paravicini}
\date{June 1, 2008}
\maketitle

\begin{abstract}
\noindent The Green-Julg theorem states that $\KTh_0^G(B) \cong \KTh_0(\Leb^1(G,B))$ for every compact group $G$ and every $G$-C$^*$-algebra $B$. We formulate a generalisation of this result to proper groupoids and Banach algebras and deduce that the Bost assembly map is surjective for proper Banach algebras.

\medskip


\medskip

\noindent \textsc{AMS 2000} \textit{Mathematics subject classification:} Primary 19K35; 43A20; 22A22; Secondary 22D15

\end{abstract}

\noindent In analogy to the definition of the assembly map of Baum-Connes, one can construct a homomorphism $\mu_{\mA}^B$ from $\KTh^{\top2}_*(\gG, B)$ to $\KTh_*(\mA(\gG,B))$, where $\gG$ is a locally compact Hausdorff groupoid with Haar system, $B$ is a $\gG$-C$^*$-algebra and $\mA(\gG)$ is an unconditional completion of $\Cont_c(\gG)$, that is, a completion with respect to a submultiplicative norm $\norm{\cdot}_{\mA}$ such that $\norm{f}_{\mA}$ only depends on the function $\gamma\mapsto \abs{f(\gamma)}$. This construction was discussed in \cite{Lafforgue:06} and an obvious generalisation to $\gG$-Banach algebras instead of $\gG$-C$^*$-algebras was given in \cite{Paravicini:07:Induction:arxiv}.

It is well-known that the Baum-Connes conjecture is true for proper C$^*$-coefficients, and this result is a key ingredient for the so-called Dirac-dual-Dirac method\footnote{See \cite{Kasparov:81:Conspectus}.} which is applied to prove the conjecture for certain classes of groups and arbitrary C$^*$-coefficients. The Baum-Connes conjecture for proper C$^*$-coefficients is also used to show the Bost conjecture for proper C$^*$-coefficients in \cite{Lafforgue:02}.

The main result of the present article asserts that the Bost assembly map $\mu_{\mA}^B$ is split surjective if the $\gG$-Banach algebra $B$ is proper (and $\mA(\gG)$ satisfies some mild condition). This is a first positive result for coefficients which are not C$^*$-algebras; the proof does not make use of C$^*$-algebraic methods either but rests on a generalised version of the Green-Julg theorem for Banach algebras:

Let $\gG$ be a \emph{proper} locally compact Hausdorff groupoid with unit space $X$ and assume that $\gG$ carries a Haar system. In \cite{Tu:99}, the following C$^*$-algebraic theorem is proved which reduces to the classical Green-Julg theorem if $\gG$ is a compact group.\footnote{Actually, Proposition 6.25 of \cite{Tu:99} is more general than cited here: It allows C$^*$-algebras in the first variable that are of a more general form. We confine ourselves to ``trivial'' coefficients in the first variable. Note that this theorem of Tu also generalises Theorem 5.4 in \cite{KasSkan:03}.}

\begin{theorem*}[Tu] If $\gG$ is $\sigma$-compact and $B$ is a $\gG$-C$^*$-algebra, then there is a canonical isomorphism
\begin{equation}\label{EquationTheoremOfTu}
\KK_{\gG}(\Cont_0(X),\ B) \ \cong \ \KK_{X/\gG}(\Cont_0(X/\gG),\, B \rtimes_r \gG).
\end{equation}
\end{theorem*}

\noindent In order to translate this theorem into the setting of Banach algebras, we choose the language of V.~Lafforgue's bivariant $\KTh$-theory $\KKban$, introduced in \cite{Lafforgue:02} and \cite{Lafforgue:06}. More precisely, we proceed as follows:
\begin{itemize}
\item We replace the $\gG$-C$^*$-algebra $B$ by a $\gG$-Banach algebra, so the left-hand side of (\ref{EquationTheoremOfTu}) should then be replaced by\footnote{Actually, it \emph{should} be replaced by $\KKbanW{\gG}(\C_X, B)$ where $\C_X$ denotes the constant field over $X$ with fibre $\C$. We will sometimes identify $\Cont_0(X)$ and $\C_X$ to obtain statements of theorems which look familiar.} $\KKbanW{\gG}(\Cont_0(X), B)$.
\item The crossed product of $B$ with $\gG$ is replaced by $\mA(\gG, B)$, where $\mA(\gG)$ is some unconditional completion of $\Cont_c(\gG)$ as, for example, $\Leb^1(\gG)$.
\item For technical reasons, we do not use $\KKbanW{X/\gG}$ on the right-hand side but a variant called $\RKKban$ which is defined in the first section of this article.
\end{itemize}

\noindent This way, we obtain the following conjecture:
\begin{equation}\label{Equation:GeneralisedGreenJulg:vorne}
\KKbanW{\gG}(\Cont_0(X),\, B) \ \cong \ \RKKban (\Cont_0(X/\gG);\, \Cont_0(X/\gG),\, \mA(\gG, B)).
\end{equation}
We show this conjecture under some mild regularity conditions: Firstly, $B$ should be a \emph{non-degenerate} $\gG$-Banach algebra, i.e., the span of $BB$ is dense in $B$. Secondly, we want $\gG$ to carry a cut-off function (which is automatic if $X/\gG$ is $\sigma$-compact). Thirdly, we want the unconditional completion $\mA(\gG)$ to be \emph{regular} (this notion will be explained in Paragraph~\ref{Subsection:RegularCompletions}; the completion $\Leb^1(\gG)$ and its symmetrised version $\Leb^1(\gG)\cap \Leb^1(\gG)^*$ are regular). Under these conditions, we define a homomorphism from the left-hand side to the right-hand side of (\ref{Equation:GeneralisedGreenJulg:vorne}) and show that it is split surjective. For split injectivity, we need that $\mA(\gG)$ satisfies some additional regularity condition (which is true if $\mA(\gG)$ equals $\Leb^1(\gG)$ or its symmetrised version). The proof of the injectivity part is only sketched in this article, and the reader is referred to \cite{Paravicini:07} for the details. Note that the above-mentioned surjectivity of the Bost assembly map for proper Banach algebras follows already from the surjectivity part of the generalised Green-Julg theorem.

\medskip

\noindent In the first section of this article, we introduce a variant $\RKKban$ of $\KKban$ for Banach algebras which carry an action of $\Cont_0(X)$, where $X$ is a locally compact Hausdorff space. This theory serves as a recipient for the descent homomorphism and also appears on the right-hand side of the generalised Green-Julg theorem. We prove that the spectral radius of an element in a $\Cont_0(X)$-Banach algebra can be calculated in the fibres of the algebra over points in $X$ if $X$ is second countable.

The second section discusses shortly the concept of a monotone completion needed in the proofs in Section~4. We also prove that $\Leb^1(G, \Cont_0(\R, B))$ and $\Cont_0(\R, \Leb^1(G,B))$ are isomorphic in $\KTh$-theory (and similar results along these lines). 

Section~3 introduces proper groupoids, cut-off functions and cut-off pairs in preparation of Section~4, where we prove the surjectivity part of the generalised Green-Julg theorem and sketch how to prove the injectivity part.

In Section~5, we define what a proper Banach algebra is and conclude from the results of Section~4 that the Bost assembly map is split surjective if the coefficients are proper Banach algebras.

\medskip

\noindent Most of the results of this article are contained in the doctoral thesis \cite{Paravicini:07} which comprises full proofs and all technical details; I would like to thank my Ph.D.\ supervisor Siegfried Echterhoff. I also thank Vincent Lafforgue, who has drawn my attention to the study of the Bost conjecture for proper Banach algebra coefficients, for his helpful advice. This research has been supported by the Deutsche Forschungsgemeinschaft (SFB 478).

\smallskip

\noindent Notation: All Banach spaces and Banach algebras that appear in this article are supposed to be complex. References which explain the necessary notation and the concepts to understand Banach algebras that carry actions of groupoids are \cite{Lafforgue:06} and \cite{Paravicini:07:Induction:arxiv}.

\section{$\Cont_0(X)$-Banach algebras and $\RKKban$-theory}

Let $X$ be a locally compact Hausdorff space. The notion of a $\Cont_0(X)$-C$^*$-algebra is well-known in the literature, and it has already been generalised to the concept of a $\Cont_0(X)$-Banach algebra.\footnote{See \cite{Blanchard:96}.} For $\Cont_0(X)$-C$^*$-algebras, there is a natural variant of $\KK$-theory called $\RKK$. This section is dedicated to the development of an analogous theory for $\Cont_0(X)$-Banach algebras. This can be thought of as an intermediate step between $\KKban$ for ordinary Banach algebras as defined in \cite{Lafforgue:02} and the variant of $\KKban$ for fields of Banach algebras as defined in \cite{Lafforgue:06}.

The starting point for our definition of $\RKK$ is the following observation: If $A$ and $B$ are $\Cont_0(X)$-C$^*$-algebras and $(E,T)$ is a cycle for $\RKK(\Cont_0(X);A,B)$, then $E$ carries a canonical action of $\Cont_0(X)$ defined through the identification $E\cong E\otimes_B B$, just let $\Cont_0(X)$ act on the second factor. This action is the unique action of $\Cont_0(X)$ on $E$ that is compatible with the module action of $B$. The usual condition on an $\RKK$-cycle, namely that $(\chi a)(e b) = (ae) (\chi b)$ for all $a\in A$, $e\in E$, $b\in B$ and $\chi\in \Cont_0(X)$, then just means that the actions of $\Cont_0(X)$ on $A$ and $E$ should be compatible. So $E$ is what could be called a $\Cont_0(X)$-Hilbert $A$-$B$-module. The corner stone for the definition of $\RKKban$ should hence be the notion of a $\Cont_0(X)$-Banach $A$-$B$-pair (if $A$ and $B$ are $\Cont_0(X)$-Banach algebras). The fundamental notion underlying all this is a notion of a $\Cont_0(X)$-Banach space, which turns out to be rather simple:

\subsection{$\Cont_0(X)$-Banach spaces, $\Cont_0(X)$-Banach algebras, etc.}\label{SubsectionCNullXBanachAlgebras}

A \demph{$\Cont_0(X)$-Banach space} is by definition a non-de\-ge\-ne\-rate Banach $\Cont_0(X)$-module. If $E$ and $F$ are $\Cont_0(X)$-Banach spaces, then we take the bounded linear $\Cont_0(X)$-linear maps from $E$ to $F$ as morphisms from $E$ to $F$. We are going to denote the morphisms from $E$ to $F$ by $\Lin^{\Cont_0(X)}(E,F)$.

If $E$ is a Banach space, then $EX=\Cont_0(X,E)$ is a $\Cont_0(X)$-Banach space with the canonical action of $\Cont_0(X)$.

Let $E_1$ and $E_2$ be $\Cont_0(X)$-Banach spaces. Let $E_1 \times E_2$ be the product Banach space (with the sup-norm). Then $E_1 \times E_2$ is a $\Cont_0(X)$-Banach space with the obvious product action. Similarly, there is a notion of the sum $E_1 \oplus E_2$ of $\Cont_0(X)$-Banach spaces $E_1$ and $E_2$ using the sum-norm. It is compatible with the $\Cont_0(X)$-tensor product that we are going to define below. Let $F$ be another $\Cont_0(X)$-Banach space. A $\C$-bilinear map $\mu\colon E_1\times E_2 \to F$ is called \demph{$\Cont_0(X)$-bilinear} if $\mu$ is $\Cont_0(X)$-linear in every component. There is a universal space $E_1 \otimes_{\Cont_0(X)} E_2$ for continuous $\Cont_0(X)$-bilinear maps on $E_1 \times E_2$, called the \demph{$\Cont_0(X)$-tensor product}. It can be constructed as a quotient of the projective tensor product $E_1 \otimes^{\pi} E_2$ and is itself a $\Cont_0(X)$-Banach space in an obvious way.

\begin{definition} A \demph{$\Cont_0(X)$-Banach algebra} $B$ is a Banach algebra $B$ which is at the same time a $\Cont_0(X)$-Banach space such that the multiplication of $B$ is $\Cont_0(X)$-bilinear.
\end{definition}
\noindent A homomorphism of $\Cont_0(X)$-Banach algebras $\varphi \colon A\to B$ is simply a $\Cont_0(X)$-linear homomorphism $\varphi$ of Banach algebras.

\emph{ For the rest of Subsection \ref{SubsectionCNullXBanachAlgebras}, let $A$, $B$ and $C$ be $\Cont_0(X)$-Banach algebras.}

We define the \demph{fibrewise unitalisation} of $B$ to be $B \oplus \Cont_0(X)$. The norm on $B \oplus \Cont_0(X)$ is the sum-norm and multiplication is given by $(b,\varphi) \cdot (c,\psi) := (bc+\psi b + \varphi c,\ \varphi\psi)$ for all $b,c\in B$, $\varphi,\psi\in \Cont_0(X)$. In the theory of $\Cont_0(X)$-Banach algebras, the fibrewise unitalisation is the adequate substitute for the ordinary unitalisation, e.g.\ it should be used in the definition of pushouts along homomorphisms of $\Cont_0(X)$-Banach algebras. We wont stress this technical point in what follows.

A \demph{$\Cont_0(X)$-Banach $B$-module} is a Banach $B$-module $E$ which is at the same time a $\Cont_0(X)$-Banach space such that the module action is $\Cont_0(X)$-bilinear. We define $\Cont_0(X)$-Banach $B$-$C$-bimodules analogously. Let $E$, $F$ be $\Cont_0(X)$-Banach $B$-modules. Then we write $\Lin_B^{\Cont_0(X)}(E,F)$ for the subspace of $\Lin_B(E,F)$ of operators which are also $\Cont_0(X)$-linear. Note that, if $E$ is a non-degenerate Banach $B$-module, then all elements of $\Lin_B(E,F)$ are automatically $\Cont_0(X)$-linear.

There is also an obvious notion of homomorphisms with coefficient maps between $\Cont_0(X)$-Banach modules, compare the definition in \cite{Paravicini:07:Morita:erschienen}.

Let $E$ be a right $\Cont_0(X)$-Banach $B$-module and let $F$ be a left $\Cont_0(X)$-Banach $B$-module. The \demph{balanced $\Cont_0(X)$-tensor product} $E\otimes_B^{\Cont_0(X)}F$ of $E$ and $F$ over $B$ is defined to be the universal object for the $B$-balanced $\Cont_0(X)$-multilinear maps on $E \times F$. It can be obtained by taking $E \otimes_B F$ and dividing out elements of the form $e \varphi \otimes f -e\otimes \varphi f$. Note that, if $E$ or $F$ is $B$-non-degenerate, then it is not hard to show that the usual balanced tensor product and the balanced $\Cont_0(X)$-tensor product agree: $E\otimes_B^{\Cont_0(X)} F= E \otimes_B F.$

The \demph{pushout} along homomorphisms of $\Cont_0(X)$-Banach algebras is defined as in the ordinary case, compare \cite{Lafforgue:02}, page 12, but using the fibrewise unitalisation defined above. It has the expected (functorial) properties.

\begin{definition} Let $B$ be a $\Cont_0(X)$-Banach algebra. A \demph{$\Cont_0(X)$-Banach $B$-pair} $E$ is a $B$-pair $E$ such that $E^<$ and $E^>$ are $\Cont_0(X)$-Banach $B$-modules and such that the inner product is $\Cont_0(X)$-bilinear. If $A$ is another $\Cont_0(X)$-Banach algebra, then a Banach $A$-$B$-pair $E$ is a \demph{$\Cont_0(X)$-Banach $A$-$B$-pair} if it is a $\Cont_0(X)$-Banach $B$-pair and the actions of $A$ on $E^<$ and $E^>$ are $\Cont_0(X)$-bilinear.
\end{definition}

\noindent For example, if $B$ is a $\Cont_0(X)$-Banach algebra, then $(B,B)$ is a $\Cont_0(X)$-Banach $B$-pair.

Let $E$ and $F$ be $\Cont_0(X)$-Banach $B$-pairs. Then an element $T$ of $\Lin_B(E,F)$ is called \demph{$\Cont_0(X)$-linear} if $T^<$ and $T^>$ are $\Cont_0(X)$-linear. The subspace of all $\Cont_0(X)$-linear maps in $\Lin_B(E,F)$ is denoted by $\Lin_B^{\Cont_0(X)}(E,F)$. 

The definitions of concurrent homomorphisms with coefficient maps between $\Cont_0(X)$-Banach pairs, the $\Cont_0(X)$-tensor product of $\Cont_0(X)$-Banach pairs and the pushout of $\Cont_0(X)$-Banach pairs along homomorphisms of $\Cont_0(X)$-Banach algebras are the obvious variation of the corresponding definitions for ordinary Banach pairs, requiring all maps to be $\Cont_0(X)$-linear (compare the discussion for Banach modules above).

\begin{proposition} Let $E$ and $F$ be $\Cont_0(X)$-Banach $B$-pairs. Then $\Komp_B(E,F)$ is always contained
in $\Lin_B^{\Cont_0(X)}(E,F)$, i.e., $\Cont_0(X)$-linearity is automatic for compact operators.
\end{proposition}
\begin{proof}
Let $f^>\in F^>$ and $e^<\in E^<$. Let $T:=\ketbra{f^>}{e^<}$. To show that $T^>$ is $\Cont_0(X)$-linear let $e^>\in
E^>$ and $\varphi\in \Cont_0(X)$. Then
\[
T^>(\varphi e^>) = f^> \langle e^<,\ \varphi e^>\rangle = f^> ( \varphi\langle e^<,\ e^>\rangle) = \varphi (f^> \langle
e^<,\ e^>\rangle ) = \varphi T^>(e^>).
\]
Similarly one shows that $T^<$ is $\Cont_0(X)$-linear. Now the set of all $\Cont_0(X)$-linear elements in $\Lin_B(E,F)$ is a closed subspace, so it contains the whole of $\Komp_B(E,F)$.
\end{proof}

\noindent It is easy to see that $\Komp_B(E,F)$ is a $\Cont_0(X)$-Banach space and that the canonical bilinear map from $F^> \times E^< \to \Komp_B(E,F)$ is $\Cont_0(X)$-bilinear. If $G$ is another $\Cont_0(X)$-Banach $B$-pair, then the composition of elements of $\Komp_B(F,G)$ and $\Komp_B(E,F)$ is $\Cont_0(X)$-bilinear. In particular, $\Komp_B(E)$ is a $\Cont_0(X)$-Banach algebra.

\begin{definition} \label{DefinitionCNullXLocallyCompactOperators} Let $E$ and $F$ be $\Cont_0(X)$-Banach $B$-pairs. Then $T\in \Lin_B(E,F)$ is called \demph{locally compact} if $\chi T$ is compact for all $\chi \in \Cont_0(X)$. 
\end{definition}

\noindent Note that it suffices to check $\chi T\in \Komp_B(E,F)$ for all $\chi\in \Cont_c(X)$. Note also that locally compact operators are automatically $\Cont_0(X)$-linear. The bounded locally compact operators form a closed subset of $\Lin^{\Cont_0(X)}_B(E,F)$.

\subsection{$\RKKbanG \left(\Cont_0(X); A,B\right)$}

\subsubsection{Gradings and group actions}

A \demph{graded $\Cont_0(X)$-Banach space} is a $\Cont_0(X)$-Banach space $E$ endowed with a grading automorphism commuting with the $\Cont_0(X)$-action.

Let $G$ be a locally compact Hausdorff group that acts continuously on $X$. Note that $\Cont_0(X)$ is a $G$-Banach algebra when equipped with the $G$-action $(g\chi )(x):= \chi(g^{-1}x)$, $\chi\in \Cont_0(X)$, $g\in G$, $x\in X$. A \demph{$G$-$\Cont_0(X)$-Banach space} is a $G$-Banach space $E$ which is at the same time a $\Cont_0(X)$-Banach space such that the actions of $G$ and $\Cont_0(X)$ are compatible in the following sense:
\[
g(\chi e) = (g\chi) (ge), \quad \chi\in \Cont_0(X), g\in G, e\in E,
\]
i.e., the product $\Cont_0(X) \times E\to E$ is $G$-equivariant.

From these definitions we also get an obvious definition of a graded $G$-$\Cont_0(X)$-Banach space. Taking this as a starting point one can define graded $G$-$\Cont_0(X)$-Banach algebras and graded equivariant homomorphisms between them, graded $G$-$\Cont_0(X)$-Banach pairs, etc.

\subsubsection{Definition of $\RKKbanG\left(\Cont_0(X);A,B\right)$}

\begin{definition} Let $A$ and $B$ be $G$-$\Cont_0(X)$-Banach algebras. Then the class $\EbanG\left(\Cont_0(X); A,B\right)$ is defined to be the class of pairs $(E,T)$ such that $E$ is a non-degenerate graded $G$-$\Cont_0(X)$-Banach $A$-$B$-pair and, if we forget the $\Cont_0(X)$-structure, the pair $(E,T)$ is an element of $\EbanG\left(A,B\right)$. Note that $T$ in the definition is automatically $\Cont_0(X)$-linear because $E$ is non-degenerate.
\end{definition}

\noindent The constructions one usually performs with $\KKban$-cycles are obviously compatible with the additional $\Cont_0(X)$-structure, so we can form the sum of $\KKban$-cycles and take their pushout along homomorphisms of $G$-$\Cont_0(X)$-Banach algebras. We also have a $\Cont_0(X)$-linear notion of morphisms of $\KKban$-cycles, giving us a $\Cont_0(X)$-linear version of isomorphisms of $\KKban$-cycles. Hence also the notion of homotopy makes sense in the $\Cont_0(X)$-setting so we can formulate the following definition:

\begin{definition} The class of all homotopy classes of elements of $\EbanG\left(\Cont_0(X);A,B\right)$ is denoted by $\RKKbanG\left(\Cont_0(X);A,B\right)$. The sum of cycles induces a law of composition on $\RKKbanG\left(\Cont_0(X);A,B\right)$ making it an abelian group.
\end{definition}

\noindent The fact that the composition on $\RKKbanG\left(\Cont_0(X);A,B\right)$ has inverses can be proved just as in the case without the $\Cont_0(X)$-structure, i.e., Lemme 1.2.5 of \cite{Lafforgue:02} and its proof are compatible with the additional $\Cont_0(X)$-module action. There is an obvious forgetful group homomorphism
\[
\RKKbanG\left(\Cont_0(X);A,B\right)\ \to\ \KKbanG\left(A,B\right).
\]

\subsubsection{A sufficient condition for homotopy}

There is a sufficient condition for the homotopy of $\RKKbanG$-cycles just as there is for $\KKbanG$-cycles, compare \cite{Paravicini:07:Morita:erschienen}, Theorem~3.1. The main idea is that the mapping cylinder of a homomorphism $\Phi$ of $\RKKbanG$-cycles gives a homotopy between the cycles. For this to be true, $\Phi$ has to satisfy a technical condition which says that the operators which are required to be compact in the definition of $\RKKbanG$-cycles can be approximated simultaneously by finite rank operators for both cycles which $\Phi$ connects. This is what is meant by ``$(\Phi, (T,T')) \in \EbanG(\Cont_0(X); \id_A,\id_B)$'' in the following theorem:

\begin{theorem}\label{Theorem:SufficientCondition:RKKbanG} Let $A$ and $B$ be $G$-$\Cont_0(X)$-Banach algebras. Let $(E,T), (E',T')$ be elements of $\EbanG(\Cont_0(X); A,B)$. If there is a $\Cont_0(X)$-linear morphism $\Phi$ from $(E,T)$ to $(E',T')$ (with coefficient maps $\id_A$ and $\id_B$) such that $(\Phi, (T,T')) \in \EbanG(\Cont_0(X); \id_A,\id_B)$, then $(E,T) \sim (E',T')$.
\end{theorem}

\noindent The necessary concepts are explained in \cite{Paravicini:07:Morita:erschienen} for $\KKbanG$; to obtain the result for $\Cont_0(X)$-Banach algebras it suffices to add compatible $\Cont_0(X)$-Banach spaces structures everywhere.

\subsection{Comparison with the $\KKban$-theory for fields of Banach algebras} \label{Subsection:ComparisonKKbanRKKban}

Let $\gG$ be a locally compact Hausdorff groupoid with unit space $X$. In \cite{Lafforgue:06}, V.~Lafforgue has introduced an equivariant $\KKban$-theory for $\gG$-Banach algebras. A $\gG$-Banach algebra is, in particular, an upper semi-continuous field of Banach algebras over $X$. If $A$ is such a field, then one can consider $\ContSect_0(X,A)$, the space of all sections of $A$ which vanish at infinity. The Banach algebra $\ContSect_0(X,A)$ carries a canonical action of $\Cont_0(X)$ making it a $\Cont_0(X)$-Banach algebra. However, it is not clear how to find an elegant way to model a general $\gG$-action on $A$ on the level of elements of $\ContSect_0(X,A)$. Nevertheless, it is rather straightforward in the case that $\gG= G\ltimes X$ where $G$ is a locally compact Hausdorff group acting on $X$. In this case, $\ContSect_0(X,A)$ is a $G$-$\Cont_0(X)$-Banach algebra in a canonical fashion. We have the following result whose proof can be found in \cite{Paravicini:07}, Section~4.7.

\begin{proposition}
Let $A$ and $B$ be $G\ltimes X$-Banach algebras. Then there is a canonical isomorphism
\[
\KKbanW{G\ltimes X}\left(A,B\right) \quad \cong \quad \RKKbanG\left(\Cont_0(X);\ \ContSect_0(X,A),\ \ContSect_0(X,B)\right).
\]
\end{proposition}

\noindent Conversely, start with a $G$-$\Cont_0(X)$-Banach algebra $\cA$. For all $x\in X$, the quotient Banach algebra $\cA_x=\cA / (\Cont_0(X\setminus \{x\})\mA)$ is called the \emph{fibre} of $\cA$ over $x$; it comes with a natural quotient map $\mA \ni a\mapsto a_x \in \mA_x$. One can regard $\Field{\cA}:=(\cA_x)_{x\in X}$ as a $G\ltimes X$-Banach algebra. Let us denote the $G$-$\Cont_0(X)$-Banach algebra $\ContSect_0(X,\Field{\cA})$ by $\Gelfand{\mA}$ and call it the \demph{Gelfand transform} of $\mA$. There is a canonical homomorphism $\iota_{\mA}$ from $\mA$ to $\Gelfand{\mA}$ which sends every $a\in \mA$ to the section $x\mapsto a_x \in A_x$. Sadly enough, $\iota_A$ needs neither be injective nor surjective, we only know that it has dense image; we do not have $\cA \cong \ContSect_0(X,\Field{\cA})$ in general. The homomorphism $\iota_A$ is isometric (and therefore an isomorphism) if and only if the $\Cont_0(X)$-Banach algebra $\cA$ is what is called \emph{locally $\Cont_0(X)$-convex}, i.e.,
\[
\forall \chi_1,\chi_2\in \Cont_0(X), \chi_1,\chi_2\geq 0, \chi_1+\chi_2\leq 1\ \forall a_1,a_2\in \mA:\ \norm{\chi_1 a_1 +  \chi_2a_2} \leq \max\{\norm{a_1},\norm{a_2}\},
\]
see \cite{Gierz:82} and also Appendix A.2 of \cite{Paravicini:07}.

If $A$ is a $G\ltimes X$-Banach algebra, then $\Cont_0(X,A)$ is automatically locally $\Cont_0(X)$-convex. Actually, $A\mapsto \ContSect_0(X,A)$ defines an equivalence of categories between the category of $G\ltimes X$-Banach algebras and the category of locally $\Cont_0(X)$-convex $G$-$\Cont_0(X)$-Banach algebras, the inverse functor being $\Field{\cdot}$. The functor $\mA \mapsto \Gelfand{\mA}$ on the category of $G$-$\Cont_0(X)$-Banach algebras therefore has its values in the subcategory of locally $\Cont_0(X)$-convex $G$-$\Cont_0(X)$-Banach algebras. It is a projector in the sense that $\Gelfand{\Gelfand{\mA}}$ is naturally isomorphic to $\Gelfand{\mA}$.

The functors $\Field{\cdot}$ and $\Gelfand{\cdot}$ can also be applied to $G$-$\Cont_0(X)$-Banach spaces, $G$-$\Cont_0(X)$-Banach pairs etc. It is an interesting fact that $\Field{\cdot}$ is multiplicative in the sense that it intertwines the (fibrewise) tensor product of $G\ltimes X$-Banach spaces and the $\Cont_0(X)$-tensor product of $G$-$\Cont_0(X)$-Banach spaces; this can be proved using the result that the $\Cont_0(X)$-tensor product of locally $\Cont_0(X)$-convex spaces is again locally $\Cont_0(X)$-convex, see \cite{Paravicini:07:CNullX:erschienen}.

If $\mA$ and $\mB$ are arbitrary $G$-$\Cont_0(X)$-Banach algebras, then it is possible to construct a group homomorphism
\[
\RKKbanG(\Cont_0(X); \mA, \mB) \to \RKKbanG(\Cont_0(X); \Gelfand{\mA}, \Gelfand{\mB}).
\]
It is not clear under which conditions this is an isomorphism if $\mA$ and $\mB$ are not locally $\Cont_0(X)$-convex. A first result along these lines is proved in the following section showing that $\mA$ and $\Gelfand{\mA}$ have the same (non-equivariant) $\KTh$-theory.

To conclude, one can say that $\KKbanW{\gG}$ and $\RKKban$ agree on the (equivalent) categories of $G\ltimes X$-Banach algebras / locally $\Cont_0(X)$-convex $G$-$\Cont_0(X)$-Banach algebras, but on the one hand, $\KKbanW{\gG}$ can be extended much further to Banach algebras which carry actions arbitrary groupoids, on the other hand, $\RKKban$ can be extended to $G$-$\Cont_0(X)$-Banach algebras which fail to be locally $\Cont_0(X)$-convex (and that such algebras appear naturally is the raison d'\^{e}tre for this theory).

A much more elaborate discussion of the two concepts can be found in Chapter~4 of \cite{Paravicini:07}.

\subsection{The spectral radius in $\Cont_0(X)$-Banach algebras}

In this section, we analyse to what extend the spectral radius of an element of a $\Cont_0(X)$-Banach algebra is determined by its fibrewise spectral radii. If $A$ is a Banach algebra, then we write $\rho_{A}(a)$ for the spectral radius of $a\in A$ in $A$.

\subsubsection{A formula for the spectral radius}

In this paragraph, let $A$ be a $\Cont_0(X)$-Banach algebra.

\begin{lemma}\label{Lemma:CNullX:Normungleichung}
For all $x_0\in X$, for all $a\in A$  and all $\varepsilon>0$ there is a neighbourhood $V$ of $x_0$ in $X$ such that for all $\chi \in \Cont_c(X)$ such that $0\leq \chi \leq 1$ and $\supp \chi \subseteq V$ we have
\[
\norm{\chi a} \leq \norm{a_{x_0}} + \varepsilon.
\]
\end{lemma}
\begin{proof}
Let $x_0\in X$, $a\in A$ and $\varepsilon>0$. Recall the following formula of Varela (see \cite{Varela:74}; it also follows from Lemme~1.10 of \cite{Blanchard:96} or Lemma~4.2.6 of \cite{Paravicini:07}):
\[
\norm{a_{x_0}} = \inf\{\norm{\chi a}:\ \chi \in \Cont_c(X)\  \exists V\subseteq X \text{ open}: \chi\restr_V =1, 0\leq \chi \leq 1, x_0 \in V\}.
\]
In particular, we can choose a $\chi' \in \Cont_c(X)$ such that $0\leq \chi' \leq 1$ and $\chi'\equiv 1$ on a neighbourhood $V$ of $x_0$ and such that $\norm{\chi' a} \leq \norm{a_{x_0}} + \varepsilon$. Let $\chi \in \Cont_c(X)$ be such that $0\leq \chi \leq 1$ and $\supp \chi \subseteq V$. Then $\chi \chi' = \chi$ and hence
\[
\norm{\chi a} = \norm{\chi \chi' a} \leq \norm{\chi}_{\infty} \norm{\chi'a} \leq \norm{a_{x_0}} + \varepsilon. \qedhere
\]
\end{proof}

\begin{lemma}\label{Lemma:CNullX:SpectralRadiusInequality}
For all $a\in A$, all $x_0\in X$  and all $\varepsilon>0$ there exists a neighbourhood $V$ of $x_0$ such that for all $\chi\in \Cont_c(X)$ such that $0\leq \chi \leq 1$ and $\supp \chi \subseteq V$ we have
\[
\rho_A(\chi a) \leq \rho_{A_{x_0}} (a_{x_0}) + \varepsilon.
\]
\end{lemma}
\begin{proof}
Let $a\in A$, $x_0\in X$  and $\varepsilon>0$. Find a $k\in\N$ such that $\norm{a_{x_0}^k}^{1/k} \leq \rho_{A_{x_0}}(a_{x_0}) + \varepsilon/2$. Apply Lemma~\ref{Lemma:CNullX:Normungleichung} to the element $a^k$ of $A$ to find a neighbourhood $V$ of $x_0$ such that
\[
\norm{ \chi' a^k} ^{1/k} \leq \left(\norm{a_{x_0}^k} + \left(\frac{\varepsilon}{2}\right)^k\right)^{1/k} \leq \norm{a_{x_0}^k}^{1/k} + \frac{\varepsilon}{2}
\]
for all functions $\chi'\in \Cont_c(X)$ such that $0\leq \chi' \leq 1$ and $\supp \chi' \subseteq V$. Let $\chi \in \Cont_c(X)$ such that $0\leq \chi \leq 1$ and $\supp \chi \subseteq V$. Then using the above inequality for $\chi'=\chi^k$ we obtain
\[
\rho_A(\chi a) \leq \norm{\chi^k a^k}^{1/k} \leq \norm{a_{x_0}^k}^{1/k} + \frac{\varepsilon}{2} \leq \rho_{A_{x_0}}(a_{x_0}) + \varepsilon. \qedhere
\]
\end{proof}

\begin{lemma}\label{Lemma:CNullX:Monotonie}
Let $\chi,\chi'\in \Cont_c(X)$ such that $0\leq \chi\leq \chi'$. Let $a\in A$. Then $\norm{\chi a} \leq \norm{\chi'a}$.
\end{lemma}
\begin{proof}
Find a function $\chi''\in \Cont_c(X)$ such that $0\leq \chi'' \leq 1$ and $\chi'' \equiv 1$ on the support of $\chi'$. Let $\varepsilon>0$. Then $\chi' + \varepsilon \chi''$ satisfies $\chi \leq \chi' + \varepsilon \chi''$. Moreover, we have that
\[
x\mapsto \frac{\chi(x)}{\chi'(x) +\varepsilon \chi''(x)}
\]
defines a continuous function on the open set $\{x' \in X: \chi''(x')>0\}$, and we have $\chi(x)=0$ for all $x$ in the open set $X\setminus  \supp \chi$. So defining
\[
\delta(x) := \begin{cases}\frac{\chi(x)}{\chi'(x) +\varepsilon \chi''(x)} &\text{ if } \chi''(x)>0\\
0 &\text{ if } \chi(x)=0\end{cases}
\]
for all $x\in X$ defines a function $\delta\in \Cont_c(X)$ such that $0\leq \delta \leq 1$. We have
\[
\delta (\chi' + \varepsilon \chi'') a = \chi a
\]
and hence
\[
\norm{\chi a} = \norm{\delta (\chi' + \varepsilon \chi'') a} \leq \norm{\delta}_{\infty} \norm{(\chi' + \varepsilon\chi'') a} \leq \norm{\chi'a} +\varepsilon\norm{a}.
\]
Because this is true for all $\varepsilon>0$, we can conclude that $\norm{\chi a} \leq \norm{\chi'a}$.
\end{proof}

\begin{proposition}\label{Proposition:CNullX:Spektralradius}
For every $a\in A$, the function
\[
x\mapsto \rho_{A_x}(a_x)
\]
is upper semi-continuous and vanishes at infinity. If $X$ is second countable, then
\[
\rho_A(a) = \max_{x\in X} \rho_{A_x}(a_x)
\]
for all $a\in A$.
\end{proposition}
\begin{proof}
For the first assertion, let $a\in A$, let $x_0\in X$ and let $\varepsilon>0$.  Find a neighbourhood $V$ of $x_0$ as in Lemma~\ref{Lemma:CNullX:SpectralRadiusInequality}. Let $\chi\in \Cont_c(X)$ be such that $0\leq \chi \leq 1$, $\supp \chi \subseteq V$ and $\chi \equiv 1$ on a neighbourhood $U$ of $x_0$. Let $x\in U$. Then
\[
\rho_{A_x}(a_x) \leq \rho_A(\chi a) \leq \rho_{A_{x_0}} (a_{x_0}) +\varepsilon
\]
follows from $(\chi a)_x = a_x$  and from Lemma~\ref{Lemma:CNullX:SpectralRadiusInequality}. From this we see that $x\mapsto \rho_{A_x}(a_x)$ is upper semi-continuous. From $\rho_{A_x}(a_x)\leq \norm{a_x}$ for all $x\in X$ it follows that $x\mapsto \rho_{A_x}(a_x)$ vanishes at infinity because $x\mapsto \norm{a_x}$ does.

Now for the second assertion. We first reduce to the case that $X$ is compact. Consider the one-point compactification $X^+$. The $\Cont_0(X)$-Banach algebra $A$ is also a $\Cont(X^+)$-Banach algebra with $A_{\infty}=\{0\}$. We therefore have $\rho_{A^+}(a) = \rho_A(a)$ and $\max_{x\in X^+} \rho_{A_x}(a_x) = \max_{x\in X} \rho_{A_x}(a_x)$ for all $a\in A$. If $X$ is second countable, then $X^+$ is compact and metrisable. Hence it suffices to consider the case that $X$ is compact and metrisable.

So let $X$ be a compact metric space. Let $a\in A$. Define
\[
m:= \max_{x\in X} \rho_{A_x}(a_{x}).
\]
We show that $\rho_A(a) \leq m + \varepsilon$ for all $\varepsilon>0$.

Let $\varepsilon>0$. For all $x\in X$ find an (open) neighbourhood $V_x$ as in Lemma~\ref{Lemma:CNullX:SpectralRadiusInequality} for $a$, $\varepsilon$ and $x$. Then $(V_x)_{x\in X}$ is an open covering of $X$. Find a finite subset $S$ of $X$ such that $(V_s)_{s\in S}$ is also a covering of $X$. Find some $\delta>0$ such that every subset of $X$ of diameter less than $\delta$ is contained in one of the sets $V_s$. Find a finite refinement $\mU$ of $\{V_s: s\in S\}$ (i.e. a finite set of open subsets of $X$ which covers $X$ and such that every element of $\mU$ is contained in a $V_s$) such that every element of $\mU$ has diameter less than $\delta/2$ (to produce such a refinement first take any open cover of $X$ by sets of diameter less than $\delta/2$, find a finite subcover; this finite subcover is automatically a refinement of $\{V_s: s\in S\}$).

Define $\mN:= \{\Delta \subseteq \mU:\ \bigcap \Delta\neq \emptyset\}$. This is a finite (combinatorial) simplicial complex. Note that $\bigcup \Delta$ has diameter less than $\delta$ for all $\Delta \in \mN$, so $\bigcup \Delta$ is contained in a $V_s$ with $s\in S$.

Let $(\chi_{U})_{U \in \mU}$ be a continuous partition of unity subordinate to the finite cover $\mU$.

If $\Delta\in \mN$, then $\chi_{\Delta}:=\sum_{U \in \Delta} \chi_U$ satisfies $0 \leq \chi_{\Delta}\leq 1$  and is supported in a set $V_s$ with $s\in S$, hence
\begin{equation}\label{Equation:chiDelta}
\rho_A(\chi_{\Delta}\ a) \leq \rho_{A_s} (a_s) +\varepsilon \leq m +\varepsilon.
\end{equation}

If $x\in X$, then $\Delta_x:= \{U \in \mU:\ x\in U\}$ is in $\mN$. So for all $x\in X$:
\[
\sum_{U \in \mU} \chi_U(x) = 1 = \sum_{U \in \Delta_x} \chi_U(x) = \chi_{\Delta_x}(x)
\]
In particular, we have
\[
1= \left(\sum_{U\in \mU} \chi_U(x)\right)^k \leq \sum_{\Delta \in \mN} \left( \chi_\Delta (x)\right)^k
\]
for all $k\in \N$ and $x\in X$. It hence follows from
\[
a^k = \left(\sum_{U\in \mU} \chi_U \right)^k a^k
\]
and Lemma~\ref{Lemma:CNullX:Monotonie} that
\[
\norm{a^k} = \norm{\left(\sum_{U\in \mU} \chi_U \right)^k a^k} \leq \sum_{\Delta \in \mN} \norm{\left( \chi_\Delta \right)^k a^k}
\]
for all $k\in \N$. It follows that
\[
\norm{a^k}^{1/k} \leq \left(\sum_{\Delta \in \mN} \norm{\left( \chi_\Delta  a\right)^k }\right)^{1/k}
\]
for all $k\in \N$. The left-hand side approaches $\rho_A(a)$ if $k\to \infty$, the right-hand side converges to
\[
\max_{\Delta \in \mN} \ \rho_{A} \left(\chi_\Delta a\right) \stackrel{(\ref{Equation:chiDelta})}{\leq} m +\varepsilon.
\]
So we have shown
\[
\rho_A(a) \leq m +\varepsilon
\]
for all $\varepsilon>0$, so $\rho_A(a) \leq m$.
\end{proof}

\subsubsection{Consequences of the spectral radius formula}

\begin{definition}
Let $A$ and $A'$ be Banach algebras and $\varphi\colon A\to A'$ be a contractive homomorphism. Then $\varphi$ is called \demph{isoradial} if $\rho_{A'}(\varphi(a)) = \rho_A(a)$ for all $a\in A$. If $\varphi$ is isoradial and has dense range, then $\varphi$ is called \demph{full}.
\end{definition}

\noindent It is a well-known fact that full homomorphisms are isomorphisms in $\KTh$-theory, see for example \cite{CMR:07}. We use this fact in the proof of the following proposition.

\begin{proposition}\label{Proposition:CNullX:FullAndDense}
Let $A$ and $A'$ be $\Cont_0(X)$-Banach algebras and let $\varphi\colon A\to A'$ be a $\Cont_0(X)$-linear contractive homomorphism of Banach algebras. Assume that $\varphi_x\colon A_x \to A'_x$ is isoradial for all $x\in X$ and that $X$ is second countable. Then $\varphi$ is an isoradial homomorphism. Moreover, if $\varphi_x$ has dense image for all $x\in X$, then $\varphi_* \colon \KTh_*(A)\to \KTh_*(A')$ is an isomorphism.
\end{proposition}
\begin{proof}
If $a\in A$, then by Lemma~\ref{Proposition:CNullX:Spektralradius}:
\[
\rho_{A'}(\varphi(a)) = \max_{x\in X} \rho_{A'_x}(\underbrace{\varphi(a)_x}_{=\varphi_x(a_x)}) = \max_{x\in X} \rho_{A_x} (a_x) = \rho_A(a).
\]
So $\varphi$ is isoradial.

Now assume that $\varphi(A)$ is fibrewise dense in $A'$.

We first consider the case that $A'$ is locally $\Cont_0(X)$-convex (see the discussion in Section~\ref{Subsection:ComparisonKKbanRKKban}). The space $\varphi(A)$ is not only fibrewise dense in $A'$, but also invariant under $\Cont_0(X)$, so it is dense in $A'$ by the Stone-Weierstrass Theorem for locally $\Cont_0(X)$-convex $\Cont_0(X)$-Banach spaces (an early variant of this is Theorem~7.9 of \cite{Hofmann:72}; see also Proposition~2.3 in \cite{DuGil:83}). Hence $\varphi_*$ is an isomorphism in this case.

Now let $A'$ be arbitrary and let $\iota_{A'} \colon A' \to \Gelfand{A'}$ be the canonical homomorphism of $A'$ into its Gelfand transform $\Gelfand{A'}$ (compare Subsection~\ref{Subsection:ComparisonKKbanRKKban}). By construction, $\iota_{A'}$ is a fibrewise isomorphism and $\Gelfand{A'}$ is locally $\Cont_0(X)$-convex. In particular, $\iota_{A'}$ is an isomorphism in $\KTh$-theory by the preceding part of the proof. Also $\iota_{A'} \circ \varphi$ is a homomorphism with full fibres, so it is also an isomorphism in $\KTh$-theory. So $\varphi$ is an isomorphism in $\KTh$-theory, as well.
\end{proof}

\noindent We have shown and used the following fact in the proof of the preceding proposition, but it is certainly worth to be stated explicitly:

\begin{corollary}
Let $A$ be a $\Cont_0(X)$-Banach algebra and $\Gelfand{A}$ be the Gelfand transform of $A$. Then the canonical map from $A$ to $\Gelfand{A}$ is full and therefore an isomorphism in $\KTh$-theory:
\[
\KTh_*(A) \cong \KTh_*(\Gelfand{A}).
\]
\end{corollary}

\subsection{Special case: $X$ compact}

We conclude this section by discussing how $\RKKbanG$ reduces to ordinary $\KKbanG$-theory if $X$ is a \emph{compact} space:

Let $G$ be a locally compact Hausdorff group and $X$ be a compact Hausdorff space on which $G$ acts. Let $A$ be a non-degenerate $G$-Banach algebra and let $B$ be a non-degenerate $G$-$\Cont(X)$-Banach algebra. Then the projective tensor product $A\otimes \Cont(X)$ is a non-degenerate $G$-$\Cont(X)$-Banach algebra.

There is a canonical forgetful homomorphism
\[
\RKKbanG\left(\Cont(X); A\otimes \Cont(X), B\right)\ \to \ \KKbanG\left(A\otimes \Cont(X),\ B\right).
\]
Secondly, there is a canonical homomorphism $j_A$ of $G$-Banach algebras from $A$ to $A\otimes \Cont(X)$, namely the map $a\mapsto a\otimes 1$. This gives a group homomorphism from $\KKbanG(A\otimes \Cont(X); B)$ to $\KKbanG(A,B)$. Let
\[
\kappa\colon \RKKbanG\left(\Cont(X); A\otimes \Cont(X),\ B\right) \to \KKbanG\left(A,B\right)
\]
be the composition of these two homomorphisms.

\begin{proposition} The homomorphism $\kappa$ is an isomorphism.
\end{proposition}
\begin{proof}
We first prove surjectivity: Let $(E,T) \in \EbanG(A,B)$. Instead of defining a $\Cont(X)$-structure on $E$, which we do not know how to do, we define a structure on the cycle $(E \otimes_B B,\ T\otimes 1) \in \EbanG(A,B)$, where $E\otimes_B B = (B\otimes_B E^<, E^>\otimes_B B)$. Note that $(E \otimes_B B,\ T\otimes 1) = (E,T) \otimes_B \MoritabanG(\id_B)$, so it is homotopic to $(E,T)$, see \cite{Paravicini:07:Morita:erschienen}, assertion (5) of Proposition~5.28. On $E^> \otimes_B B$ we define a canonical $\Cont(X)$-structure: if $e^>\in E^>$ and $b\in B$ and $\varphi \in \Cont(X)$, then $\varphi(e^>\otimes b) := e^>\otimes (\varphi b)$. This makes $E^> \otimes_B B$ a right $G$-$\Cont(X)$-Banach $B$-module. We proceed similarly on the left-hand side. It is easy to see that $E\otimes_B B$ is a $G$-$\Cont(X)$-Banach $B$-pair with this $\Cont(X)$-action. The operator $T \otimes 1$ is clearly $\Cont(X)$-linear (which is automatic anyway, because $E\otimes_B B$ is non-degenerate).

Now we have to define an action of $A\otimes \Cont(X)$ on $E\otimes_B B$: If $a\in A$, $\chi \in \Cont(X)$, $e^>\in E^>$ and $b\in B$, then we define $(a \otimes \chi)(e^> \otimes b):= (ae^>) \otimes (\chi b)$. This gives an action of $A \otimes \Cont(X)$ on $E^> \otimes_B B$ making it a $G$-$\Cont(X)$-Banach $A\otimes \Cont(X)$-$B$-bimodule. A similar definition can be made for the left-hand side. We check that $A\otimes \Cont(X)$ acts on $E\otimes_B B$ by elements of $\Lin_B(E\otimes_B B)$. Let therefore be $a\in A$, $\chi \in \Cont(X)$, $e^<\in E^<$, $e^> \in E^>$ and $b^<,b^>\in B$. Then
\begin{eqnarray*}
\left\langle b^< \otimes e^<,\ (a \otimes \chi) (e^> \otimes b^>) \right\rangle &=& \left\langle b^< \otimes e^<,\ (ae^>) \otimes (\chi b^>)\right\rangle = b^< \left\langle e^<,\ ae^>\right\rangle (\chi b^>)\\
&=& (\chi b^<) \left\langle e^<a, \ e^>\right\rangle b^> = \left\langle (b^<\otimes e^<)(a \otimes\chi),\ e^>\otimes b^>\right\rangle.
\end{eqnarray*}
By trilinearity and continuity of both sides, this equation can be extended from the elementary tensors to all of $A \otimes \Cont(X)$, $B \otimes_B E^<$ and $E^>\otimes_B B$. So $E\otimes_B B$ is in $\EbanG(\Cont(X); A\otimes \Cont(X), B)$. Applying $\kappa$ to it means forgetting the $\Cont(X)$-structure and reducing the $A\otimes \Cont(X)$-action back to the $A$-action on $E\otimes_B B$, so we are back to where we started. Hence $\kappa$ is surjective.

\medskip

\noindent The same argument shows that $\kappa$ is injective: Let $(E_0,T_0)$ and $(E_1,T_1)$ be elements of the class $\EbanG(\Cont(X); A\otimes \Cont(X), B)$ such that $\kappa(E_0,T_0)$ and $\kappa(E_1,T_1)$ are homotopic in $\EbanG(A,B)$. Find $(E,T) \in \EbanG(\Cont(X);A\otimes \Cont(X),\ B[0,1])$ such that $\kappa(E,T) \in \EbanG(A,B[0,1])$ is a homotopy from $\kappa(E_0,T_0)$ to $\kappa(E_1,T_1)$. Now $\ev^B_{i,*}(E,T)$ is contained in $\EbanG(\Cont(X); A\otimes \Cont(X),\ B)$ for all $i\in \{0,1\}$ and $\kappa(\ev^B_{i,*}(E,T))$ is isomorphic (in $\EbanG(A,B)$) to $(E_i,T_i)$. Now $E_i$ is a non-degenerate $B$-pair, so it is easy to see that the $\Cont(X)$-structure on $E$ is unique. Hence the isomorphism between $\kappa(\ev^B_{i,*}(E,T))$ and $(E_i,T_i)$ must be $\Cont(X)$-linear. Also the action of $A\otimes \Cont(X)$ is uniquely determined by the actions of $A$ and $\Cont(X)$, so the isomorphism from $\kappa(\ev^B_{i,*}(E,T))$ and $(E_i,T_i)$ must also respect this structure. In other words, it is an isomorphism of cycles in $\EbanG(\Cont(X); A\otimes \Cont(X),B)$. So $(E_0,T_0)$ and $(E_1,T_1)$ are homotopic. Hence $\kappa$ is injective.
\end{proof}

\noindent If we take $A$ to be $\C$ with the trivial $G$-action, then $A\otimes \Cont(X)$ is isomorphic to $\Cont(X)$. The proposition then reduces to the following statement:

\begin{corollary}\label{Corollary:RKKAndKKforCompactBaseSpace} Let $B$ be a non-degenerate $G$-$\Cont(X)$-Banach algebra. If $X$ is compact, then
\[
\RKKbanG\left(\Cont(X);\ \Cont(X), B\right) \quad \cong \quad \KKbanG\left(\C,B\right).
\]
This, together with Théorème~1.2.8 of \cite{Lafforgue:02}, implies that, if $X$ is compact and $G$ is the trivial group, then
\[
\RKKban\left(\Cont(X);\ \Cont(X), B\right) \quad \cong \quad \KKban\left(\C,B\right) \quad \cong \quad \KTh_0(B).
\]
\end{corollary}

\section{Monotone completions}\label{SectionMonotoneCompletions}

In \cite{Lafforgue:02} and \cite{Lafforgue:06}, the notion of an unconditional completion was introduced which is a special case of what we propose to call a \demph{monotone} completion. Already the article \cite{Lafforgue:02} provides us with some interesting examples of monotone completions which are not unconditional completions.\footnote{For example $H^2(G,A)$ defined after Lemme 1.6.5 or the ``normalised'' completions $\Leb^{p,l}_{\text{\rm norm}}(G,A)$ appearing in 4.5.} The difference simply is that an unconditional completion is required to carry a product making it a Banach algebra whereas an unconditional completion is a Banach space without any product.

\subsection{Definition}

Let $Y$ be a locally compact Hausdorff space.

\begin{definition} \label{DefinitionMonotoneSeminorm}
A semi-norm $\norm{\cdot}_{\mH}$ on $\Cont_c(Y)$ is called \emph{monotone} if the following condition holds:
\begin{equation}\label{PropertyMonotoneNorm}
\forall \varphi_1,\varphi_2 \in \Cont_c(Y):\quad \left(\forall y \in Y:\ \abs{\varphi_1(y)} \leq \abs{\varphi_2(y)}\right) \Rightarrow \norm{\varphi_1}_{\mH} \leq \norm{\varphi_2}_{\mH}.
\end{equation}
Let $\mH(Y)$ denote the (Hausdorff-)completion of $\Cont_c(Y)$ with respect to this semi-norm; this Banach space is called a \demph{monotone completion} of $\Cont_c(Y)$.
\end{definition}

\noindent \emph{For the rest of this section, let $\mH(Y)$ be a monotone completion of $\Cont_c(Y)$.}

For technical reasons and as for unconditional norms, we extend monotone norms to a larger class of functions on $Y$:

\begin{definition} Let $\mF_c\left(Y\right)$ be the set of all (locally) bounded functions $\varphi\colon Y\to \R$ with compact support. Let $\mF_c^+\left(Y\right)$ be the set of elements of $\mF_c\left(Y\right)$ which are non-negative. Define
\[
\norm{\varphi}_\mH := \inf\left\{\norm{\psi}_\mH:\ \psi \in \Cont_c(Y),\ \psi\geq \varphi\right\}
\]
for all $\varphi \in \mF_c^+\left(Y\right)$.
\end{definition}

\noindent Note that, by Property (\ref{PropertyMonotoneNorm}), the new semi-norm agrees on $\Cont_c^+(Y)$ with the semi-norm we started with. For all $\varphi_1,\varphi_2,\varphi\in \mF_c^+(Y)$ and all $c\geq 0$, we have
\begin{enumerate}
    \item $\varphi_1 +\varphi_2 \in \mF_c^+(Y)$ and $\norm{\varphi_1 +\varphi_2}_{\mH} \leq \norm{\varphi_1}_{\mH} +\norm{\varphi_2}_{\mH}$;
    \item $c\varphi \in \mF_c^+(Y)$ and $\norm{c\varphi}_{\mH} = c \norm{\varphi}_{\mH}$;
    \item if $\varphi_1 \leq \varphi_2$, then $\norm{\varphi_1}_{\mH} \leq \norm{\varphi_2}_{\mH}$.
\end{enumerate}
Hence we can use the extended semi-norm to define a semi-norm on sections of u.s.c.~fields of Banach spaces.

\emph{For the rest of this section, let $E$ be a u.s.c.~field of Banach spaces over $Y$.}

\begin{definition}\label{DefinitionMonotoneCompletionOfFieldsOfBanachSpaces} We define the following semi-norm on $\ContSect_c(Y, E)$:
\[
\norm{\xi}_{\mH}:= \norm{y \mapsto \norm{\xi(y)}_{E_{y}}}_{\mH}.
\]
The Hausdorff completion of $\ContSect_c(Y,E)$ with respect to this semi-norm will be denoted by $\mH(Y,E)$.
\end{definition}

\noindent Note that the function $y \mapsto \norm{\xi(y)}$ appearing in the preceding definition is not necessarily continuous. However, it has compact support and is non-negative upper semi-continuous, so we can apply the extended semi-norm on $\mF_c^+(Y)$ to it.

If $E$ is the trivial bundle over $Y$ with fibre $E_0$, then $\ContSect_c(Y,E)$ is $\Cont_c(Y, E_0)$. The completion $\mH(Y,E)$ of $\Cont_c(Y,E_0)$ could hence also be denoted as $\mH(Y,E_0)$ and might be considered as a sort of tensor product of $\mH(Y)$ and $E_0$. If in particular $E_0 =\C$, then $\mH(Y,E) = \mH(Y,\C)= \mH(Y)$.

\begin{definition}\label{DefinitionMonotoneCompletionOfLinearMaps} Let $F$ be another u.s.c.~field of Banach spaces over $Y$ and let $T$ be a bounded continuous field of linear maps from $E$ to $F$. Then $\xi \mapsto T\circ \xi$ is a linear map from $\ContSect_c(Y,E)$ to $\ContSect_c(Y,F)$ such that $\norm{T\circ \xi}_{\mH} \leq \norm{T} \norm{\xi}_{\mH}$. Hence $T$ induces a canonical continuous linear map from $\mH(Y,E)$ to $\mH(Y,F)$ with norm $\leq \norm{T}$.
\end{definition}

\noindent This way, we define a functor from the category of u.s.c.~fields of Banach spaces over $Y$ to the category of Banach spaces, which is linear and contractive on the morphism sets.

Note that the canonical map from $\ContSect_c(Y, E)$ to $\mH(Y, E)$ is continuous if we take the inductive limit topology on $\ContSect_c(Y, E)$ and the norm topology on $\mH(Y,E)$. It follows that, if a subset of $\ContSect_c(Y,E)$ is dense in $\ContSect_c(Y,E)$ for the inductive limit topology, then its canonical image in $\mH(Y,E)$ is dense for the norm topology.

\subsection{An application to monotone and unconditional completions}

In this section, we prove that $\mA(G, \Cont_0(\R, B))$ and $\Cont_0(\R, \mA(G,B))$ are isomorphic in $\KTh$-theory, where $G$ is a locally compact Hausdorff group, $B$ is a $G$-Banach algebra and $\mA(G)$ is an unconditional completion of $\Cont_c(G)$ (this is \ref{Corollary:UnconditionalAndSuspension:Group}). This result already appeared in the work of V.~Lafforgue, see \cite{Lafforgue:02}, Section~1.7, and also a variant for groupoids can be found in \cite{Lafforgue:06}. However, no proof has been published yet.

Here, the original, direct proof that V.~Lafforgue has indicated to me is generalised in several directions. Firstly, we replace $\R$ by a general second countable locally compact Hausdorff space $X$. The main part of the argument is now a statement about $\Cont_0(X)$-Banach algebras, namely the spectral radius formula of Proposition~\ref{Proposition:CNullX:Spektralradius} and its consequence Proposition~\ref{Proposition:CNullX:FullAndDense}. To complete the proof of Lafforgue's result, what is left to show is that the fibres of $\mA(G, \Cont_0(X, B))$ and $\Cont_0(X, \mA(G,B))$ are isometrically isomorphic. This is mainly a statement not about Banach algebras, but about Banach spaces.

So we replace the group $G$ (or, more generally, the groupoid $\gG$) with an arbitrary locally compact Hausdorff space $Y$ and the unconditional completion $\mA(G)$ with a monotone completion $\mH(Y)$. Let $E=(E_y)_{y\in Y}$ be an upper semi-continuous field of Banach spaces over $Y$. Finally, let $X$ be another locally compact Hausdorff space which for this part of the argument does not need to be second countable. The result we are going to show can now be formulated as follows:

\begin{proposition*}
The two $\Cont_0(X)$-Banach spaces $\mH(Y, EX)$ and $\mH(Y,E)X$ have the same fibres over points in $X$.
\end{proposition*}

\subsubsection*{A more precise formulation}

Let $\pi_1\colon Y\times X \to Y$ and $\pi_2\colon Y\times X \to X$ be the canonical projections. Note that $\pi_1^*E=(E_y)_{(y,x)\in Y\times X}$ is an upper semi-continuous field of Banach spaces over $Y\times X$. The pushforward $\pi_{1,*}(\pi_1^* E)$ is an upper semi-continuous continuous field over $Y$, the fibre over $y\in Y$ being isomorphic to $E_yX=\Cont_0(X, E_y)$. We denote this pushforward by $EX=(E_yX)_{y\in Y}$.

We form the Banach space $\mH(Y, EX)$ and compare it to $\Cont_0(X, \mH(Y,E))=\mH(Y,E)X$. Note that there is a canonical contractive linear map $\iota$ from $\mH(Y,EX)$ to $\mH(Y,E)X$.

The second space is actually a $\Cont_0(X)$-Banach space. If $x\in X$, then the fibre of $\mH(Y,E)X$ over $x$ is canonically isomorphic to $\mH(Y,E)$. On the other hand, there is a canonical contractive linear map $e_x$ from $\mH(Y, EX)$ to $\mH(Y,E)$: If $\ev_x^E$ denotes the canonical evaluation map at $x$ from $EX$ to $E$, then $e_x=\mH(Y, \ev_x^E)$. The following diagram commutes
\[
\xymatrix{
\mH(Y,EX) \ar[r]^{\iota} \ar[dr]^{e_x}& \mH(Y,E)X\ar[d]^{\ev_x^{\mH(Y,E)}}\\
& \mH(Y,E).
}
\]
Note that there is a canonical $\Cont_0(X)$-structure on $\mH(Y, EX)$: If $\xi \in \ContSect_c(Y, EX)$ and $\chi\in \Cont_0(X)$, then $(\chi \xi)(y) := \chi \cdot \xi(y)$ for all $y\in Y$; to interpret this formula note that $\xi(y)$ is contained in the fibre $(EX)_y$ of $EX$ over $y\in Y$ which is isomorphic to $\Cont_0(X,E_y)$ (as mentioned above), hence there is a canonical product between $\Cont_0(X)$ and $(EX)_y$. The linear map $\iota$ is clearly $\Cont_0(X)$-linear and hence it induces a canonical homomorphism $\iota_x$ from the fibre $\mH(Y, EX)_x$ to $(\mH(Y,E)X)_x\cong \mH(Y,E)$ for every $x\in X$. We can hence extend the above diagram to the following commutative square
\[
\xymatrix{
\mH(Y,EX) \ar[r]^{\iota} \ar[dr]^{e_x} \ar[d]_{\pi_x}& \mH(Y,E)X\ar[d]^{\ev_x^{\mH(Y,E)}}\\
\mH(Y,EX)_x \ar[r]_{\iota_x}& \mH(Y,E)
}
\]
where we write $\pi_x$ for the canonical projection which maps $f\in \mH(Y,EX)$ to $f_x\in \mH(Y,EX)_x$.

So a more precise formulation of the above proposition is:

\begin{proposition} For all $x_0\in X$, the map $\iota_{x_0}\colon \mH(Y,EX)_{x_0}\to \mH(Y,E)$ is an isometric isomorphism.
\end{proposition}

\subsubsection*{A proof of this proposition}

Let $x_0\in X$. Because $\iota(\ContSect_c(Y, EX))$ is dense in $\mH(Y,E)X$, we know that $\iota_{x_0}$ has dense image. Hence it suffices to show that $\iota_{x_0}$ is isometric. We now show:

For all $f\in \mH(Y, EX)$  and all $\varepsilon>0$ there is a function $\chi \in \Cont_c(X)$ such that $0\leq \chi \leq 1$, $\chi(x_0)= 1$ and such that
\begin{equation}\label{Equation:NormInequality:YetAnotherOne}
\norm{\chi f}_{\mH} \leq \norm{e_{x_0}(f)}_{\mH} + \varepsilon.
\end{equation}
This is sufficient because $\iota_{x_0}(f_{x_0}) = e_{x_0}(f)$ and $f_{x_0} = (\chi f)_{x_0}$, which implies
\[
\norm{\iota_x(f_{x_0})} \leq \norm{f_{x_0}} \leq  \norm{(\chi f)_{x_0}} = \norm{\chi f}_{\mH} \stackrel{!}{\leq} \norm{e_{x_0}(f)}_{\mH} + \varepsilon = \norm{\iota_{x_0}(f_{x_0})} +\varepsilon.
\]
If we prove this for arbitrary $\varepsilon>0$, then we have shown that $\iota_{x_0}$ is isometric (note that the first inequality follows from the fact that $\norm{\iota_{x_0}} \leq \norm{\iota} \leq 1$).

We first treat the case that $f\in \ContSect_c(Y \times X,\pi_1^*E ) \subseteq \ContSect_c(Y, EX)$. Let $\varepsilon>0$. The set $K:=\pi_1(\supp f)\subseteq Y$ is compact. Let $U$ be a compact neighbourhood of $K$ in $Y$. Because $\mH(Y)$ is a monotone completion of $\Cont_c(Y)$, we can find a constant $C\geq 0$ such that $\norm{\xi}_{\mH}\leq C\norm{\xi}_{\infty}$ for all $\xi \in \Cont_c(Y)$ such that $\supp \xi \subseteq U$.

Because $(y,x) \mapsto \norm{f(y,x)-f(y,x_0)}_{E_y}$ is upper semi-continuous and vanishes on the compact set $K\times \{x_0\}$, we can choose a neighbourhood $V$ of $x_0$ in $X$ such that
\[
\sup_{(y,x)\in K \times V} \norm{f(y,x) - f(y,x_0)}_{E_y} \leq \frac{\varepsilon}{C}.
\]
Let $\chi\in \Cont_c(X)$ be a function such that $0\leq \chi \leq 1$, $\chi(x_0)=1$ and $\supp \chi \subseteq V$.

Choose a function $\delta_K \in \Cont_c(\gG)$ such that $0\leq \delta_K \leq 1$, $\delta_K \equiv 1$ on $K$ and $\supp \delta_K \subseteq U$. For all $y\in K$ we have
\begin{eqnarray*}
\sup_{x \in X} \norm{\chi(x) f(y, x)} &= & \sup_{x \in V} \norm{\chi(x) f(y, x)}\leq \sup_{x \in V} \norm{ f(y, x)}\\
&\leq & \sup_{x\in V} \left(\norm{f(y,x)-f(y,x_0)} + \norm{f(y,x_0)} \right) \leq \norm{f(y,x_0)} + \frac{\varepsilon}{C} \delta_K(y),
\end{eqnarray*}
and for $y\in Y\setminus K$ we have
\[
\sup_{x \in X} \norm{\chi(x) f(y, x)}=0 \leq \norm{f(y,x_0)} + \frac{\varepsilon}{C}\delta_K(y).
\]
Because $\norm{\delta_K}_{\mH} \leq C$, we have
\begin{eqnarray*}
\norm{\chi  f}_{\mH} &=& \norm{y\mapsto \sup_{x\in X} \norm{\chi(x) f(y,x)}}_{\mH} \leq  \norm{y\mapsto  \norm{f(y,x_0)} + \frac{\varepsilon}{3C} \delta_K(y) }_{\mH}\\  &\leq& \norm{y\mapsto  \norm{f(y,x_0)}}_{\mH} +\frac{\varepsilon}{C}\norm{\delta_K }_{\mH} \leq \norm{e_{x_0} (f)}_{\mH} + \varepsilon.
\end{eqnarray*}
We now treat the general case, so let $f$ be an arbitrary element of $\mH(Y, EX)$. Let $\varepsilon>0$. Then we can find an $f'\in \ContSect_c(Y\times X, \pi_1^* E)$ such that $\norm{f-f'}_{\mH} \leq \varepsilon/3$. Note that this also implies $\norm{e_{x_0}(f)-e_{x_0}(f')}_{\mH} \leq \varepsilon/3$. By the first part of the proof we can find a function $\chi \in \Cont_c(X)$ such that $0\leq \chi \leq 1$, $\chi(x_0)=1$ and
\[
\norm{\chi f'}_{\mH} \leq \norm{e_{x_0}(f')}_{\mH} + \frac{\varepsilon}{3}.
\]
Now
\[
\norm{\chi f}_{\mH} \leq  \norm{\chi f'}_{\mH} + \frac{\varepsilon}{3} \leq  \norm{e_{x_0}(f')}_{\mH} + \frac{2 \varepsilon}{3} \leq \norm{e_{x_0}(f)}_{\mH} + \varepsilon. \qedhere
\]

\subsubsection*{Some corollaries}

\begin{corollary}
Let $G$ be a locally compact Hausdorff group and let $B$ be a $G$-Banach algebra. Let $\mA(G)$ be an unconditional completion of $\Cont_c(G)$. Let $X$ be a locally compact Hausdorff space. Then the canonical homomorphism of $\Cont_0(X)$-Banach algebras
\[
\iota\colon \mA(G, BX) \to \mA(G,B)X
\]
is an isometric isomorphism on the fibres. If $X$ is second countable, then this means that $\iota$ is full and hence
\[
\iota_* \colon \KTh_*(\mA(G,BX)) \cong \KTh_*(\mA(G,B)X).
\]
\end{corollary}
\begin{proof}
Take $(G,X,\mA,B_G)$ instead of $(Y,X,\mH,E)$ in the above proposition. Here $B_G$ denotes the trivial field of Banach spaces over $G$ with fibre $B$.
\end{proof}

\begin{corollary}\label{Corollary:UnconditionalAndSuspension:Group} In particular, taking $X=\R$ in the preceding corollary, we obtain an isomorphism
\[
\KTh_*(\mA(G,SB)) \cong \KTh_*(S\mA(G,B)) \cong \KTh_{*+1}(\mA(G,B)),
\]
where $S$ denotes the suspension functor for Banach algebras $A\mapsto \Cont_0(\R,A)$.
\end{corollary}

\begin{corollary}
Let $\gG$ be a locally compact Hausdorff groupoid equipped with a Haar system and let $B$ be a $\gG$-Banach algebra. Let $\mA(\gG)$ be an unconditional completion of $\Cont_c(\gG)$. Let $X$ be a locally compact Hausdorff space. Then the canonical homomorphism of $\Cont_0(X)$-Banach algebras
\[
\iota\colon \mA(\gG, BX) \to \mA(\gG,B)X
\]
is an isometric isomorphism on the fibres (note that $BX$ is a $\gG$-Banach algebra in a canonical fashion). If $X$ is second countable, then $\iota$ is full, and hence
\[
\iota_* \colon \KTh_*(\mA(\gG,BX)) \cong \KTh_*(\mA(\gG,B)X).
\]
\end{corollary}
\begin{proof}
Take $(\gG,X,\mA,r^*B)$ instead of $(Y,X,\mH,E)$ in the above proposition.
\end{proof}

\begin{corollary}\label{Corollary:UnconditionalAndSuspension:Groupoid} We have an isomorphism
\[
\KTh_*(\mA(\gG,SB)) \cong \KTh_*(S\mA(\gG,B)) \cong \KTh_{*+1}(\mA(\gG,B)).
\]
\end{corollary}

\section{Proper groupoids}\label{Section:ProperGroupoids}

A locally compact Hausdorff groupoid is called \demph{proper} if the following map is proper, i.e., inverses of compact sets are compact:
\[
\gG \to \gG^{(0)}\times \gG^{(0)},\ \gamma \mapsto (r(\gamma), s(\gamma)).
\]
We collect some examples:
\begin{enumerate}
\item Let $G$ be a locally compact Hausdorff group acting from the left on a locally compact Hausdorff space $X$. Then the transformation groupoid $G \ltimes X$ is proper if and only if the action of $G$ on $X$ is proper.
\item More generally, if $\gG$ is a locally compact Hausdorff groupoid and $X$ is a left $\gG$-space, then $\gG\ltimes X$ is proper if and only if $X$ is a proper $\gG$-space.
\item A locally compact Hausdorff \emph{group} is proper (as a groupoid) if and only if it is compact.
\item If the range and source maps of a locally compact Hausdorff groupoid are equal, the groupoid can be regarded as a bundle of groups. If such a groupoid is proper, then all the fibres are compact groups.
\end{enumerate}

\noindent \emph{For the remainder of Section~\ref{Section:ProperGroupoids}, let $\gG$ be a locally compact proper Hausdorff groupoid with unit space $X$ and carrying a Haar system $\lambda$.}

In particular, this means that $X/\gG$ is a locally compact Hausdorff space.

\subsection{Proper groupoids and the descent}

Let $\mA(\gG)$ be an unconditional completion of $\Cont_c(\gG)$ and let $E$ be a $\gG$-Banach space.\footnote{See \cite{Lafforgue:06} or \cite{Paravicini:07:Induction:arxiv} for the definitions of these concepts.} For all $\xi \in \ContSect_c(\gG, r^*E)$ and $\chi \in \Cont_0(X/\gG)$ define
\[
(\chi \xi) (\gamma):= \chi(\pi(\gamma)) \xi(\gamma)
\]
for all $\gamma\in \gG$, where $\pi$ denotes the (open) projection map $\pi\colon \gG \to X/\gG$. This defines a module action of $\Cont_0(X/\gG)$ on $\ContSect_c(\gG, r^*E)$ which lifts to a module action on $\mA(\gG, E)$. More precisely, $\mA(\gG,E)$ is a non-degenerate Banach $\Cont_0(X/\gG)$-module, i.e., it is a $\Cont_0(X/\gG)$-Banach space. Note that, depending on the choice of $\mA(\gG)$, the $\Cont_0(X/\gG)$-Banach space $\mA(\gG,E)$ does not have to be locally $\Cont_0(X/\gG)$-convex; it is however in important cases, e.g.~if $\mA(\gG)=\Leb^1(\gG)$.

The convolution product and also the descent of continuous linear maps respects the $\Cont_0(X/\gG)$-structure; in particular, if $B$ is a $\gG$-Banach algebra, then $\mA(\gG,B)$ is not only a Banach algebra but a $\Cont_0(X/\gG)$-Banach algebra, and if $E$ is a $\gG$-Banach $B$-pair, then $\mA(\gG,E)$ is a $\Cont_0(X/\gG)$-Banach $\mA(\gG,B)$-pair, etc.
Let $A$ and $B$ be $\gG$-Banach algebras. It is not hard to show that the descent homomorphism from $\KKbanW{\gG}(A,B)$ to $\KKban(\mA(\gG,A), \mA(\gG,B))$ introduced in Section~1.3 of \cite{Lafforgue:06} is indeed a homomorphism
\[
j_{\mA}\colon \KKbanW{\gG}(A,B) \to \RKKban(\Cont_0(X/G);  \mA(\gG,A), \mA(\gG,B)).
\]
Because $\mA(\gG,B)$ does not have to be locally $\Cont_0(X/\gG)$-convex in general, it seems advisable to use $\RKKban$ instead of $\KKbanW{X/\gG}$, the version of $\KKban$ for fields over $X/\gG$.

\subsection{Cut-off functions and cut-off pairs}

\begin{definition}\footnote{Compare \cite{Tu:99}, Définition 6.7.}
A continuous function $c\colon X\to [0,\infty[$ is called \demph{cut-off function} for $\gG$ if
\begin{enumerate}
\item $\forall x\in X:\ \int_{\gG^x} c(s(\gamma)) \rmd \lambda^x(\gamma) =1$;
\item $r\colon \supp (c\circ s) \to X$ is proper.
\end{enumerate}
\end{definition}
\noindent The latter condition means that $\supp c \cap \gG K$ is compact for all compact subsets $K$ of $X$.

Recall from \cite{Tu:04} that there is a cut-off function for $\gG$ if  $X/\gG$ is $\sigma$-compact.

Given a cut-off function $c$, one often uses the function $c^{\frac{1}{2}}$ in the theory of C$^*$-algebras. In the Banach algebra setting, the exponent $\frac 1 2$ is no longer the inevitable choice, also $c^{1/p}$ with $1<p<\infty$ can appear quite naturally. Because we are dealing with Banach pairs rather than Banach modules, it even makes sense to extend the notion of a cut-off function as follows:

\begin{definition}\label{Definition:CutOffPair} A \demph{cut-off pair} for $\gG$ is a pair $(c^<, c^>)$ such that
\begin{enumerate}
    \item $c^< \in \Cont(X)_{\geq 0}$ with $r\colon \supp(c^<\circ s)\to X$ proper;
    \item $c^> \in \Cont(X)_{\geq 0}$ with $r\colon \supp(c^>\circ s)\to X$ proper;
    \item $\forall x\in X: \int_{\gG^x} c^<(s(\gamma))\ c^>(s(\gamma)) \rmd \lambda^{x} (\gamma)=1$.
\end{enumerate}
\end{definition}

\noindent In particular, $x\mapsto c^<(x) c^>(x)$ is a cut-off function. Conversely, if $c$ is a cut-off function for $\gG$ and $p, p'\in ]1,\infty[$ such that $\frac{1}{p} + \frac{1}{p'}=1$, then $(c^{1/p'},\ c^{1/p})$ is a cut-off pair.  We can even cover the case $p=1$:

\begin{proposition}\label{Proposition:CutOffPairLOneCNull}
If $\gG$ is such that $X/\gG$ is $\sigma$-compact and $c$ is a cut-off function for $\gG$, then there exists a function $d\in \Cont(X)$ with $\norm{d}_\infty =1$ such that $(d,c)$ is a cut-off pair.
\end{proposition}
\begin{proof}
Let $(K_n)_{n\in \N}$ be an exhausting sequence of compacts in $X/\gG$ such that $K_n$ is contained in the interior of $K_{n+1}$ for all $n\in \N$. Define $L_n:= \supp c \ \cap \ \pi^{-1}(K_n)$ for all $n\in \N$ (where $\pi$ denotes the canonical surjection from $X$ to $X/\gG$). Then the $L_n$ are all compact. Recursively, find functions $f_1, f_2,f_3 \ldots$ such that $f_n \in \Cont_c(\pi^{-1}(K_n))$, $0\leq f_n\leq 1$ and $f_n \restr_{L_n} \equiv 1$ and $f_n \subseteq f_{n+1}$ for all $n\in \N$. Define $f:= \bigcup_{n\in \N} f_n$. Then this is a well-defined continuous function on $X$ such that $0\leq f\leq 1$. It satisfies $f\restr_{\supp c} \equiv 1$. Moreover, it satisfies the support condition: Let $K\subseteq X/\gG$ be compact. Find an $n\in \N$ such that $K\subseteq K_n$. Then the closed set $\pi^{-1}(K)$ is contained in $\pi^{-1}(K_n)$, so $\pi^{-1}(K) \cap \supp f$ is contained in $\pi^{-1}(K_n) \cap \supp f= \pi^{-1}(K_n) \cap \supp f_n= \supp f_n$. Now $\supp f_n$ is a compact subset of $\pi^{-1}(K_n)$, so $\pi^{-1}(K)\cap \supp f$ is compact as a closed subset of a compact subset.
\end{proof}

\noindent On the level of functions with compact support, we can define a homomorphism from $\Cont_c(X/\gG)$ to $\Cont_c(\gG)$ quite generally; it is a delicate question for which completions of $\Cont_c(\gG)$ this homomorphism can be extended continuously to $\Cont_0(X/\gG)$.

\begin{defprop}
Let $(c^<, c^>)$ be a cut-off pair for $\gG$. For all $\chi \in \Cont_c(X/\gG)$, define
\[
(\varphi(\chi))(\gamma) := c^>(r(\gamma))\ \chi (\pi(\gamma))\ c^<(s(\gamma))
\]
for all $\gamma\in \gG$. Then $\varphi(\chi) \in \Cont_c(\gG)$, and $\varphi$ is a continuous homomorphism of algebras from $\Cont_c(X/\gG)$ to $\Cont_c(\gG)$ (with the convolution product).
\end{defprop}
\begin{proof}
Let $\pi\colon X\to X/\gG$ denote the quotient map and let $K\subseteq X/\gG$ be the support of $\chi$. Then $K_1:=\supp c^<\cap \pi^{-1} (K)$ is compact in $X$ and so is $K_2:=\supp c^>\cap \pi^{-1} (K)$. So $\{\gamma \in \gG:\ s(\gamma) \in K_1,\ r(\gamma)\in K_2\}$ is compact and contains the support of $\varphi(\chi)$. So $\varphi(\chi) \in \Cont_c(\gG)$.

Let $\chi_1,\chi_2\in \Cont_c(\gG)$. Then for all $\gamma\in \gG$:
\begin{eqnarray*}
&& \left(\varphi(\chi_1)*\varphi(\chi_2)\right)(\gamma)\\
&=& \int_{\gG^{r(\gamma)}} c^>(r(\gamma'))\, \chi_1(\pi(\gamma'))\, c^<(s(\gamma')) \, c^>(r(\gamma'^{-1} \gamma))\, \chi_2(\pi(\gamma'^{-1}\gamma))\, c^<(s(\gamma'^{-1} \gamma)) \rmd \lambda^{r(\gamma)}(\gamma')\\
&=& c^>(r(\gamma))\ (\chi_1 \chi_2)(\pi(\gamma))\ c^>(s(\gamma)) \underbrace{\int_{\gG^{r(\gamma)}} c^<(s(\gamma'))\, c^>(s(\gamma')) \rmd\lambda^{r(\gamma)}(\gamma')}_{=1} = (\varphi(\chi_1 \chi_2))(\gamma).\qedhere
\end{eqnarray*}
\end{proof}
\noindent In the C$^*$-algebra case, the interesting cut-off pair is of course $(c^{\frac{1}{2}},\ c^{\frac{1}{2}})$, where $c$ is a cut-off function for $\gG$. In this case,\footnote{See Proposition 6.23 in \cite{Tu:99} for a proof.} the homomorphism $\varphi \colon \Cont_c(X/\gG) \to \Cont_c(\gG)$ preserves the involution and can be extended to a $*$-homomorphism from $\Cont_0(X/\gG)$ to $\Cred(\gG)$. The pullback along this $*$-ho\-mo\-morphism gives a homomorphism of groups from $\KK_{X/\gG}(\Cred(\gG),\, B \rtimes_r \gG)$ to $\KK_{X/\gG}(\Cont_0(X/\gG),\, B \rtimes_r \gG)$.

Can the same homomorphism $\varphi\colon \Cont_c(X/\gG) \to \Cont_c(\gG)$ be extended to a homomorphism from $\Cont_0(X/\gG)$ to $\mA(\gG)$ if $\mA(\gG)$ is an unconditional completion of $\Cont_c(\gG)$? This would come in handy in the construction of a homomorphism from $\KKbanW{\gG}(\Cont_0(X),\, B)$ to $\RKKban(\Cont_0(X/\gG);\, \Cont_0(X/\gG),\, \mA(\gG, B))$ where $B$ is a $\gG$-Banach algebra, see Section~\ref{Section:GreenJulg}. One could simply take the descent homomorphism and compose it with the pullback along $\varphi$.

Apparently, $\varphi$ is not bounded even for rather elementary unconditional completions like $\Leb^1(\gG)$ and rather simple cut-off pairs. The construction works for C$^*$-algebras because the choice of the cut-off pair is compatible with the norm on $\Cred(\gG)$ which is defined through the action of $\Cont_c(\gG)$ on $\Leb^2(\gG)$. We have to find another way to define the homomorphism for our generalised Green-Julg theorem in Section~\ref{Section:GreenJulg} if we do not want to deal with the technical problems that come with unbounded homomorphisms or with the compression of a Banach algebra by an unbounded projection.

\subsection{Automatic equivariance}

There is another feature of proper groupoids which will prove very convenient in the upcoming sections:

\begin{proposition}\label{Proposition:ProperWLOGOperatorEquivariant} Let $A$ and $B$ be $\gG$-Banach algebras (with $\gG$ being proper and allowing a cut-off function). Then the operators and homotopies in the definition of $\KKbanW{\gG}(A,\ B)$ can be assumed to be $\gG$-equivariant.
\end{proposition}

\begin{proof}
The basic idea here, as in the proof of the corresponding result for C$^*$-algebras, is to use the cut-off function and the integration with respect to the Haar system to make given operators equivariant, compare the discussion before Proposition~6.24 in \cite{Tu:99}. On a technical level, we do this by integrating fields of operators with compact support; note that we define this integration pointwise:

Let $E$ and $F$ be $\gG$-Banach $B$-pairs. Let $T=(T^<,T^>) \in \Lin_{r^*B}\left(r^*E,\ r^*F\right)$ have compact support. Then
\[
\int_{\gG^x} T_{\gamma} \rmd \lambda^{x}(\gamma):= \left(\int_{\gG^x} T_{\gamma}^< \rmd \lambda^{x}(\gamma),\ \int_{\gG^x} T^>_{\gamma} \rmd \lambda^{x}(\gamma)\right)
\]
is a continuous field of linear operators from $E$ to $F$. The same definition makes sense if $T$ has proper support, i.e., if the support of $(\chi \circ r) \cdot T$ is compact for all $\chi \in \Cont_c(X)$. The operator $\int_{\gG^x} T_{\gamma} \rmd \lambda^{x}(\gamma)$ is compact if $T\in \Komp_{r^*B}\left(r^*E,\ r^*F\right)$ has compact support.

We can use this procedure to produce equivariant operators. Fix a cut-off function $c$ for $\gG$. For all $T\in \Lin_B(E,F)$, we define
\[
T^{\gG}_x=\int_{\gG^x} c(s(\gamma)) \ \gamma T_{s(\gamma)} \rmd \lambda^{x}(\gamma),\qquad x\in X,
\]
Then $T^{\gG}$ is an equivariant element of $\Lin_B(E,F)$. The construction commutes with the pushout: If $B'$ is another $\gG$-Banach algebra and $\varphi\colon B\to B'$ is a $\gG$-equivariant homomorphism, then $\varphi_* \left(T^{\gG}\right) = \left(\varphi_*(T)\right)^{\gG}$ as elements of $\Lin^{\gG}_{B'} \left(\varphi_*(E),\varphi_*(F)\right)$.

Now let $(E,T)\in \EbanW{\gG}\left(A,\ B\right)$. Then $\left(E, T^{\gG}\right)$ is in $\EbanW{\gG}\left(A ,\ B\right)$ and homotopic to $(E,T)$.  To see this, let $a\in \ContSect_c(X,A)$. For all $x\in X$, we have
\begin{eqnarray*}
a_x \left(T_x - T^{\gG}_x \right) &=& \int_{\gG^x}  c(s(\gamma)) a_{r(\gamma)} \left(T_{r(\gamma)} - \gamma T_ {s(\gamma)}\right) \rmd \lambda^{x}(\gamma).
\end{eqnarray*}
The family $\gamma \mapsto  c(s(\gamma)) a_{r(\gamma)} (T_{r(\gamma)} - \gamma T_{s(\gamma)})$ is compact and of compact support, so the integral is compact. So $T$ and $T^{\gG}$ ``differ by a compact operator''. By Lemma~3.19 of \cite{Paravicini:07:Induction:arxiv}, $(E,T^{\gG})$ is a $\KKban$-cycle and homotopic to $(E,T)$. 

We have a similar result for homotopies: If $(E_0,T_0)$ and $(E_1,T_1)$ are homotopic in $\EbanW{\gG}\left(A,\ B\right)$ and if $T_0$ and $T_1$ are equivariant, then there is an equivariant homotopy between them.

This shows that the map $(E,T) \mapsto (E,T^{\gG})$ is bijective on the level of $\KKbanW{\gG}$-classes.
\end{proof}

\section{A generalised Green-Julg theorem}\label{Section:GreenJulg}

\emph{In Section~\ref{Section:GreenJulg}, let $\gG$ be a locally compact Hausdorff groupoid with Haar system $\lambda$. Let $\mA(\gG)$ be an unconditional completion of $\Cont_c(\gG)$ and let $B$ be a non-degenerate $\gG$-Banach algebra.}

\noindent The generalised Green-Julg theorem that we prove in this section asserts that we have an isomorphism
\begin{equation}\label{Equation:GeneralisedGreenJulg:hinten}
\KKbanW{\gG}(\Cont_0(X),\, B) \cong \RKKban (\Cont_0(X/\gG);\, \Cont_0(X/\gG),\, \mA(\gG, B))
\end{equation}
if $\gG$ is a proper groupoid. We construct this isomorphism only under certain conditions, more precisely, we proceed as follows:
\begin{enumerate}
\item We define a natural homomorphism $\GreenJulgAbstiegKKban{\mA}{B}$ from the left-hand side to the right-hand side of (\ref{Equation:GeneralisedGreenJulg:hinten}) in case that $\gG$ admits a cut-off function.
\item We define a natural homomorphism $\GreenJulgAufstiegKKban{\mA}{B}$ in the other direction in case that $\mA(\gG)$ is what we call \emph{regular}.
\item We show $\GreenJulgAbstiegKKban{\mA}{B} \circ \GreenJulgAufstiegKKban{\mA}{B} =\id$ if both conditions are satisfied.
\item We sketch how to show $\GreenJulgAufstiegKKban{\mA}{B}\circ \GreenJulgAbstiegKKban{\mA}{B}  =\id$ if $\mA(\gG)$ satisfies some additional regularity condition.
\end{enumerate}
Note that already the split surjectivity of $\GreenJulgAbstiegKKban{\mA}{B}$ is an interesting result as it implies the split surjectivity of the Bost-map with proper coefficients for many unconditional completions, as shown in Section~\ref{Section:BostConjectureAndProperAlgebras}. We state the surjectivity part of the generalised Green-Julg theorem for further reference:

\begin{theorem}\label{Theorem:NaturalSplitPartOfGreenJulg} Let $\mA(\gG)$ be a regular unconditional completion of $\Cont_c(\gG)$. Let there exists a cut-off function for $\gG$. Then the natural homomorphism
\[
\GreenJulgAbstiegKKban{\mA}{B}\colon \KKbanW{\gG}(\Cont_0(X),\, B) \ \to \ \RKKban (\Cont_0(X/\gG);\, \Cont_0(X/\gG),\, \mA(\gG, B))
\]
is split surjective (with natural split $\GreenJulgAufstiegKKban{\mA}{B}$) for all non-degenerate Banach algebras $B$.
\end{theorem}

\noindent The complete generalised Green-Julg theorem will be stated in Paragraph~\ref{Subsection:SketchOfProofOfGreenJulg} after some additional technical concept is introduced which is needed for the formulation of the conditions under which we can show the injectivity result; its proof is rather lengthy and will only be sketched.

\subsection{The homomorphism $\GreenJulgAbstiegKKban{\mA}{B}$}\label{Subsection:GreenJulgAbstieg}

\emph{In \ref{Subsection:GreenJulgAbstieg}, assume that $\gG$ is proper and admits a cut-off function.}

\subsubsection{The algebraic construction of $\GreenJulgAbstiegKKban{\mA}{B}$ on the level of sections with compact support}

Let $E$ be a $\gG$-Banach $B$-pair. Define the operations
\[
(e^> \beta)(x) := \int_{\gG^x} \gamma e>(s(\gamma)) \gamma \beta(\gamma^{-1}) \rmd \lambda^{x}(\gamma)
\]
and
\[
(\beta e^<)(x):= \int_{\gG^x} \beta(\gamma) \gamma e^<(s(\gamma)) \rmd \lambda^{x}(\gamma),
\]
where $x\in X$, and the $\ContSect_c(\gG,\ r^*B)$-valued bracket
\[
\lAngle e^<,\ e^>\rAngle (\gamma):= \left\langle e^<(r(\gamma)),\ \gamma e^>(s(\gamma))\right\rangle_{E_{r(\gamma)}},
\]
where $\gamma\in \gG$, for all $e^<\in \ContSect_c(X,E^<)$, $e^>\in \ContSect_c(X,E^>)$ and $\beta \in \ContSect_c\left(\gG,\ r^*B\right)$.

This turns $\ContSect_c(X, E^>)$ into a right $\ContSect_c\left(\gG, r^*B\right)$-module and $\ContSect_c(X,E^<)$ into a left $\ContSect_c\left(\gG,\ r^*B\right)$-module. These module actions are separately continuous, and they are non-degenerate for the inductive limit topologies if $E$ is non-degenerate. The bracket is $\C$-bilinear and $\ContSect_c\left(\gG,\ r^*B\right)$-linear on the left and on the right. Moreover, it is separately continuous for the inductive limit topologies.

Moreover, there are canonical actions of $\Cont(X/\gG)$ on the modules $\ContSect_c(X,E^<)$ and $\ContSect_c(X,E^>)$ given by
\[
(\chi e^>)(x):= \chi(\pi(x)) e^>(x)
\]
for all $\chi \in \Cont\left(X/\gG\right)$, $e^>\in \ContSect_c(X,E^>)$ and $x\in X$ (and analogously for the left-hand side). The module actions and the bracket are compatible with these actions.

\medskip

Let $E$ and $F$ be $\gG$-Banach $B$-pairs and let $T$ be a $\gG$-equivariant continuous field of operators from $E$ to $F$. For all $e^>\in \ContSect_c(X,E^>)$, define
\[
(\ContSect_c(X, T^>)e^>)(x) := T^>_x (e^>(x))
\]
for all $x\in X$. Then $e^>\mapsto \ContSect_c(X,T^>)e^>$ is $\C$-linear, $\Cont\left(X/\gG\right)$-linear, $\ContSect_c\left(\gG,\ r^*B\right)$-linear on the right and continuous for the inductive limit topology. The same formula defines a similar map $f^<\mapsto \ContSect_c(X,T^<)f^<$ on the left-hand side. The pair $\left(f^<\mapsto \ContSect_c(X,T^<)f^<,\ e^>\mapsto \ContSect_c(X,T^>)e^>\right)$ of linear operators is formally adjoint with respect to the brackets on $\left(\ContSect_c(X,E^<),\ \ContSect_c(X,E^>)\right)$ and $\left(\ContSect_c(X,F^<),\ \ContSect_c(X,F^>)\right)$:
\[
\lAngle f^<\ContSect_c(X,T^<),\ e^>\rAngle =\lAngle f^<,\ \ContSect_c(X,T^>)e^>\rAngle.
\]

\subsubsection{The analytic part of the construction of $\GreenJulgAbstiegKKban{\mA}{B}$} \label{SubsectionAnalyticPartOfTheConstructionOfTheHomomorphism}

In the C$^*$-world, the right module $\ContSect_c(\gG,\ r^*B)$-action and the inner product on $\ContSect_c(X,E)$ is sufficient to define the structure  $B\rtimes_r \gG$-Hilbert module if $E$ is a Hilbert $B$-module. There can only be one norm on $\ContSect_c(X,E)$ which completes to a Hilbert module and the bracket actually gives such a norm.

In the Banach-world, the situation is more complicated. If $B$ is a $\gG$-Banach algebra and $E$ is a $\gG$-Banach $B$-pair, then we will see that there are several ways to complete $\ContSect_c(X,E^<)$ and $\ContSect_c(X,E^>)$ to give a $\Cont_0(X/\gG)$-Banach $\mA(\gG,B)$-pair. However, it turns out that every (monotone) pair of such completions will give rise to the same homomorphism $\GreenJulgAbstiegKKban{\mA}{B}$.

Let $\mD^<(X)$ and $\mD^>(X)$ be monotone completions of $\Cont_c(X)$. Assume that the pair $\mD(X):=\left(\mD^<(X),\ \mD^>(X)\right)$ satisfies the following compatibility conditions with $\mA(\gG)$:
\begin{enumerate}
\item[{\bf (D1)}]
\begin{eqnarray*}
\forall \chi^< \in \Cont_c(X),\ \beta\in \Cont_c(\gG):\ && \norm{\beta \chi^<}_{\mD^<} \leq \norm{\beta}_{\mA} \norm{\chi^<}_{\mD^<} \quad \text{and} \\
\forall \chi^> \in \Cont_c(X),\ \beta\in \Cont_c(\gG):\ && \norm{\chi^>\beta}_{\mD^>} \leq \norm{\chi^>}_{\mD^>}\norm{\beta}_{\mA}.
\end{eqnarray*}
\item[{\bf (D2)}]  \qquad \qquad \quad $\forall \chi^< \in \Cont_c(X),\ \chi^>\in \Cont_c(X):\ \quad  \norm{\lAngle \chi^<,\chi^>\rAngle}_{\mA} \leq \norm{\chi^<}_{\mD^<}\norm{\chi^>}_{\mD^>}$.
\end{enumerate}
We can extend the actions of $\Cont_c(\gG)$ on $\Cont_c(X)$ from the left and from the right and also the inner product to continuous bilinear maps which turn $\mD(X)$ into a Banach $\mA(\gG)$-pair. Note that the action of $\Cont_0(X/\gG)$ on $\Cont_c(X)$ also gives a continuous non-degenerate action of $\Cont_0(X/\gG)$ on $\mD^<(X)$ and $\mD^>(X)$ making $\mD(X)$ a $\Cont_0(X/\gG)$-Banach $\mA(\gG)$-pair.

Let $E=(E^<,E^>)$ be a $\gG$-Banach $B$-pair. On $\ContSect_c(X,E^<)$ define $\norm{\xi^<}_{\mD^<}:= \norm{\ x\mapsto \norm{\xi^<(x)}\ }_{\mD^<}$ as in Definition~\ref{DefinitionMonotoneCompletionOfFieldsOfBanachSpaces} and define a semi-norm $\norm{\cdot}_{\mD^>}$ on $\ContSect_c(X,E^>)$ similarly. Then the actions of $\ContSect_c(\gG, \ r^*B)$ on $\ContSect_c(X,E^<)$ and on $\ContSect_c(X,E^>)$ and the bracket satisfy
\[
\norm{\beta \xi^<}_{\mD^<} \leq \norm{\beta}_{\mA} \norm{\xi^<}_{\mD^<},\ \ \norm{\xi^> \beta}_{\mD^>}\leq \norm{\xi^>}_{\mD^>} \norm{\beta}_{\mA},\ \ \norm{\lAngle \xi^<,\ \xi^>\rAngle}_{\mA} \leq \norm{\xi^<}_{\mD^<} \norm{\xi^>}_{\mD^>}
\]
for all $\beta\in \ContSect_c(\gG,\ r^*B)$, $\xi^<\in \ContSect_c(X,E^<)$ and $\xi^>\in \ContSect_c(X,E^>)$. As in Definition~\ref{DefinitionMonotoneCompletionOfFieldsOfBanachSpaces}, write $\mD^<(X,E^<)$ for the completion of $\ContSect_c(X,E^<)$ for the semi-norm $\norm{\cdot}_{\mD^<}$; define $\mD^>(X,E^>)$ analogously. With the extensions of the actions of $\ContSect_c(\gG,\ r^*B)$ and the extension of the bracket,
\[
\mD(X,E):= \left(\mD^<(X,E^<),\ \mD^>(X,E^>)\right)
\]
is a $\Cont_0(X/\gG)$-Banach $\mA(\gG,B)$-pair.

If $F$ is another $\gG$-Banach $B$-pair and $T \in \Lin_B(E,F)$ is $\gG$-equivariant, then $\ContSect_c(X, T^>)$ is a bounded linear map from $\ContSect_c(X,E^>)$ to $\ContSect_c(X,F^>)$ with norm less than or equal to $\norm{T^>}$, so it extends to a bounded $\C$-linear, $\Cont_0\left(X/\gG\right)$-linear and $\mA(\gG, B)$-linear map $\mD(X,T^>)$ from $\mD(X, E^>)$ to $\mD(X,F^>)$ of the same norm. Similarly, one gets a bounded linear map $\mD(X,T^<)$ from $\mD(X,F^<)$ to $\mD(X,E^<)$. Together, this defines a linear operator
\[
\mD(X, T):= \left(\mD(X,T^<),\ \mD(X,T^>)\right)\ \in \ \Lin_{\mA(\gG,B)}^{\Cont_0(X/\gG)}\left(\mD(X,E),\ \mD(X,F)\right)
\]
of norm less than or equal to $\norm{T}$. The assignment $E\mapsto \mD(X,E)$ and $T\mapsto \mD(X,T)$ is a contractive functor from the category $\gG$-Banach $B$-pairs and bounded $\gG$-equivariant operators to the category of $\Cont_0\left(X/\gG\right)$-Banach $\mA(\gG,B)$-pairs. Similarly, one can define $\mD(X,\Phi)$ for $\gG$-equivariant concurrent homomorphisms.

We omit the longsome proof of the following result which can be found in \cite{Paravicini:07}, Section~7.2.3:

\begin{proposition}\label{PropositionProperDescentRespectsLocallyCompact}
Let $S\in \Lin_B(E,F)$ be bounded, $\gG$-equivariant and locally compact. Then $\mD(X,S)$ is locally compact in the sense of Definition~\ref{DefinitionCNullXLocallyCompactOperators}, i.e., $\chi \mD(X,S)$ is compact for all $\chi \in \Cont_c\left(X/\gG\right)$.
\end{proposition}

\subsubsection{The construction for $\KKban$-cycles}

\begin{theorem}
Let $(E,T)$ be a cycle in $\EbanW{\gG}\left(\Cont_0(X),\ B\right)$. We assume that  $T$ is $\gG$-equivariant, compare Proposition~\ref{Proposition:ProperWLOGOperatorEquivariant}. Equip $\mD\left(X, E\right)$ with the obvious grading operator. Then $\mD\left(X, T\right)$ is odd and
\[
\GreenJulgAbstiegEban{\mA}{\mD}{B} (E,T):=\left(\mD(X,E),\ \mD(X,T)\right) \in \Eban\left(\Cont_0\left(X/\gG\right);\ \Cont_0\left(X/\gG\right),\ \mA\left(\gG, B\right)\right).
\]
\end{theorem}
\begin{proof}
The important property that we have to check is that $\mD(X,T)^2-1$ is locally compact. But
\[
\mD(X,T)^2-1 = \mD(X, \ T^2-1),
\]
and $T^2-1$ is locally compact. Since $T^2-1$ is also $\gG$-equivariant, we can apply Proposition~\ref{PropositionProperDescentRespectsLocallyCompact} which implies that $\mD(X, \ T^2-1)$ is locally compact.
\end{proof}

\noindent It can be shown that the map $\GreenJulgAbstiegEban{\mA}{\mD}{B}$ is compatible with the pushforward along equivariant homomorphisms of $\gG$-Banach algebras, with homotopies and with the sum of cycles; see \cite{Paravicini:07}, Section 7.2.3, for the proofs.

As a consequence of these results, we have:

\begin{proposition} The map $(E,T) \mapsto \left(\mD(X, E),\ \mD(X,T)\right)$ gives rise to a group-homomorphism from
\[
\GreenJulgAbstiegEban{\mA}{\mD}{B} \colon \KKbanW{\gG}\left(\Cont_0(X),\ B\right)\ \to \ \RKKban\left(\Cont_0\left(X/\gG\right);\ \Cont_0\left(X/\gG\right),\ \mA\left(\gG, B\right)\right)
\]
which is natural in the non-degenerate $\gG$-Banach algebra $B$.
\end{proposition}

\subsubsection{Uniqueness and existence}\label{Subsubsection:UniquenessAndExistenceOfDX}

The following uniqueness result was shown in \cite{Paravicini:07}:

\begin{defprop}\label{Defprop:GreenJulgAbstieg} Let $\mD(X)$ and $\mD'(X)$ be pairs of monotone completions of $\Cont_c(X)$ which both satisfy (D1) and (D2). Then $\GreenJulgAbstiegEban{\mA}{\mD}{B} = \GreenJulgAbstiegEban{\mA}{\mD'}{B}$ as homomorphisms from $\KKbanW{\gG}\left(\Cont_0(X),\ B\right)$ to $\RKKban\left(\Cont_0\left(X/\gG\right);\ \Cont_0\left(X/\gG\right),\ \mA\left(\gG, B\right)\right)$. We hence write $\GreenJulgAbstiegKKban{\mA}{B}$ for this homomorphism.
\end{defprop}

\noindent Another question is whether such pairs $\mD(X)$ of monotone completions exist. We have a positive answer because we have assumed that $\gG$ admits a cut-off function; there are even quite a few such completions: for every cut-off pair $c$, we construct a compatible pair of monotone completions that we call $\mA^c(X)$.

So let $c=(c^<,c^>)$ be a cut-off pair for $\gG$. Let $E$ be a $\gG$-Banach $B$-pair. Define
\[
j_{E,c}^<\colon \ContSect_c\left(X,E^<\right) \to \ContSect_c\left(\gG,\ E^<\right),\ e^<\mapsto \left(\gamma \mapsto c^<(s(\gamma)) e^<(r(\gamma))\right)
\]
and
\[
j_{E,c}^>\colon \ContSect_c\left(X,E^>\right) \to \ContSect_c\left(\gG,\ E^>\right),\ e^>\mapsto \left(\gamma \mapsto c^>(r(\gamma)) \gamma e^>(s(\gamma))\right).
\]
Then $j_{E,c}=(j_{E,c}^<,\ j_{E,c}^>)$ is a pair of injective maps such that
\begin{enumerate}
    \item $j_{E,c}^<$ is $\C$-linear, $\ContSect_c(X/\gG)$-linear and $\ContSect_c\left(\gG,\ r^*B\right)$-linear on the left,
    \item $j_{E,c}^>$ is $\C$-linear, $\ContSect_c(X/\gG)$-linear and $\ContSect_c\left(\gG,\ r^*B\right)$-linear on the right,
    \item for all $e^<\in \ContSect_c(X,E^<)$ and $e^>\in \ContSect_c(X,E^>)$, we have
    \[
    \left\langle j_{E,c}^<(e^<),\ j_{E,c}^>(e^>)\right\rangle_{\ContSect_c(\gG,r^*B)} = \lAngle e^<,\ e^>\rAngle.
    \]
\end{enumerate}
Define a $\Cont_0(X/\gG)$-Banach $\mA(\gG,B)$-pair $\mA^c(X,E)=\left(\mA^c(X,E^<),\ \mA^c(X,E^>)\right)$ by pulling back the norms of $\mA(\gG, E)$ along $j_{E,c}$ and completing $\ContSect_c(X,E)$ for this norms. Alternatively, one could take the closure of the image of $j_{E,c}$. The norms on the left and the right part are given by
\[
\norm{e^<}_{\mA^c(X,E^<)}:= \norm{j_{E,c}^<(e^<)}_{\mA(\gG,E^<)} = \Big\|\ \gamma \mapsto c^<(s(\gamma)) \norm{e^<(r(\gamma))}\ \Big\|_{\mA}
\]
and
\[
\norm{e^>}_{\mA^c(X,E^>)}:= \norm{j_{E,c}^>(e^>)}_{\mA(\gG,E^>)} = \Big\|\ \gamma \mapsto c^>(r(\gamma)) \norm{e^>(s(\gamma))}\ \Big\|_{\mA}
\]
for all $e^<\in \ContSect_c(X,E^<)$ and $e^>\in \ContSect_c(X,E^>)$.

Note that the norms depend on $\mA(\gG)$ as well as on $c$. The pair $\mA^c(X)=\left((\mA^c)^<(X),\ (\mA^c)^>(X)\right)$ is a pair of monotone completions of $\Cont_c(X)$ satisfying (D1) and (D2). Note that $\GreenJulgAbstiegEban{\mA}{\mA^c}{B}$ as a homomorphism from $\KKbanW{\gG}\left(\Cont_0(X), B\right)$ to $\RKKban\left(\Cont_0\left(X/\gG\right);\ \Cont_0\left(X/\gG\right),\ \mA\left(\gG, B\right)\right)$ does not depend on $c$ by \ref{Defprop:GreenJulgAbstieg}; without the detour via more general compatible pairs $\mD(X)$ of monotone completions this latter fact seems to be hard to prove.

\subsection{Monotone completions as analogues of $\Leb^2(\gG,B)$} \label{SectionMonotoneCompletionsAsGeneralisationsOfLTwo}

If $B$ is a $\gG$-C$^*$-algebra, then there is a canonical $\gG$-Hilbert $B$-module $\Leb^2(\gG,B)$ with a left action of $B\rtimes_r \gG$. We want to find an analogue of this module for the case that $B$ is a general $\gG$-Banach algebra. Apparently, it is not sufficient (or not systematic, at least) to just consider pairs of the type $(\Leb^2(\gG,B),\ \Leb^2(\gG,B))$; we want to treat rather general unconditional completions, so it seems appropriate to consider rather general completions of the space $\ContSect_c(\gG,\ r^*B)$ as substitutes of $\Leb^2(\gG,B)$; our treatment should at least cover pairs of the form $(\Leb^1(\gG,B), \ \ContSect_0(\gG,\ B))$ or $(\Leb^{p'}(\gG,B),\ \Leb^p(\gG,B))$ for $p,p'\in ]1,\infty[$ with $1/p+1/p'=1$ (compare the precise definitions below).

Our substitute for $\Leb^2(\gG)$ is a general pair of monotone completions of $\Cont_c(\gG)$ which satisfies some compatibility conditions with $\mA(\gG)$ and the action of $\gG$; we will  usually denote such a pair by $\mH(\gG)$, and write $\mH(\gG,B)$ for its version with coefficients in $B$. It seems advisable to even consider pairs of the form $\mH(\gG, E)$ where $E$ is a $\gG$-Banach $B$-pair because this makes the constructions a bit clearer. The important result is that (under certain conditions) the unconditional completion $\mA(\gG,B)$ acts on $\mH(\gG,B)$ by locally compact operators. This allows us to use the tensor product $\otimes_{\mA(\gG,B)} \mH(\gG,B)$ to turn $\mA(\gG,B)$-pairs into $\gG$-Banach $B$-pairs preserving locally compact operators between them.

Recall that, in this section, $\gG$ denotes a locally compact Hausdorff groupoid with left Haar system $\lambda$ and $X$ denotes the unit space of $\gG$. Recall also that $\mA(\gG)$ is an unconditional completion of $\Cont_c(\gG)$ and that $B$ is a non-degenerate $\gG$-Banach algebra.

Let $\mH(\gG) = \left(\mH^<(\gG),\ \mH^>(\gG)\right)$ be a pair of monotone completions of $\Cont_c(\gG)$ such that the bilinear map
\[
\langle \cdot,\cdot\rangle_{\Cont_c(X)}\colon \Cont_c(\gG) \times \Cont_c(\gG) \to \Cont_c(X), \ (\varphi^<,\varphi^>) \mapsto \left(x\mapsto \int_{\gG^x} \varphi^<(\gamma)\ \varphi^>(\gamma^{-1}) \rmd \lambda^x(\gamma) \right)
\]
satisfies
\begin{itemize}
\item [{\bf (H1)}] \qquad  \qquad $\forall \varphi^<,\varphi^>\in \Cont_c(\gG): \quad \norm{\langle \varphi^<,\ \varphi^>\rangle_{\Cont_c(X)}}_\infty \leq \norm{\varphi^<}_{\mH^<} \norm{\varphi^>}_{\mH^>}$.
\end{itemize}
In this case, $\langle\cdot,\cdot\rangle_{\Cont_c(X)}$ can be extended to a continuous bilinear map $\langle\cdot,\cdot\rangle_{\Cont_0(X)}\colon \mH^<(\gG) \times \mH^>(\gG) \to \Cont_0(X)$ which is $\Cont_0(X)$-bilinear if we consider the following actions of $\Cont_0(X)$:
\[
(\chi \xi^<)(\gamma) := \chi(r(\gamma)) \xi^<(\gamma) \LazyAnd (\xi^>\chi)(\gamma):= \xi^>(\gamma) \chi(s(\gamma))
\]
for all $\chi \in \Cont_0(X)$, $\xi^<\in \Cont_c(\gG)\subseteq \mH^<(\gG)$, $\xi^>\in \Cont_c(\gG) \subseteq \mH^>(\gG)$ and $\gamma\in \gG$.

\begin{examples}\label{Examples:PairsOfMonotoneCompletions}
Let $p\in [1,\infty[$. Define the norm
\[
\norm{\chi^<}_{p,r} := \sup_{x\in X} \left(\int_{\gG^x} \abs{\chi^<(\gamma)}^p \rmd \lambda^{x}(\gamma) \right)^{\frac{1}{p}}
\]
for all $\chi^<\in \Cont_c(\gG)$. The corresponding monotone completion is called $\Leb^p_r(\gG)$. Note that $\Leb^1(\gG)=\Leb^1_r(\gG)$. Secondly, define
\[
\norm{\chi^>}_{p,s} := \sup_{x\in X} \left(\int_{\gG^x} \abs{\chi^>(\gamma^{-1})}^p \rmd \lambda^{x}(\gamma) \right)^{\frac{1}{p}}
\]
for all $\chi^>\in \Cont_c(\gG)$. The corresponding monotone completion is called $\Leb^p_s\left(\gG\right)$
\begin{enumerate}
\item The pairs $(\Leb^1(\gG), \ \Cont_0(\gG))$ and $(\Cont_0(\gG),\ \Leb^1_s(\gG))$ are pairs of monotone completions of $\Cont_c(\gG)$ satisfying (H1).
\item If $p,p'\in ]1,\infty[$ such that $\frac{1}{p} +\frac{1}{p'}=1$, then $(\Leb^{p'}_r(\gG),\ \Leb^p_s(\gG))$ also satisfies (H1).
\item In particular, this applies to $(\Leb^{2}_r(\gG),\ \Leb^2_s(\gG))$.
\end{enumerate}
\end{examples}

\noindent Now let $E$ be a $\gG$-Banach $B$-pair. Define a right action of $\ContSect(X,B)$ on $\ContSect_c(\gG, \ r^*E^>)$ by
\[
(\xi^> \beta)(\gamma) := \xi^>(\gamma)\, \gamma \beta(s(\gamma)),\quad \xi^>\in \ContSect_c(\gG,r^*E^>),\,  \beta\in \ContSect(X,B),\, \gamma\in \gG,
\]
and a left action of $\ContSect(X,B)$ on $\ContSect_c(\gG,r^*E^<)$ by
\[
(\beta \xi^<)(\gamma) := \beta(r(\gamma)) \, \xi^<(\gamma),\quad \beta \in \ContSect(X,B),\,  \xi^<\in \ContSect_c(\gG,\ r^* E^<),\, \gamma\in \gG.
\]
These actions define continuous actions of $\ContSect_0(X,B)$ on $\mH^>(\gG,E^>)$ (from the right) and $\mH^<(\gG,E^<)$ (from the left). Define a bilinear map
\begin{eqnarray*}
\langle \cdot,\cdot\rangle_{\ContSect_c(X,B)}\colon \ContSect_c(\gG,\ r^*E^<) \times \ContSect_c(\gG,\ r^*E^>) &\to& \ContSect_c(X,B),\\ (\xi^<,\xi^>)&\mapsto & \left(x\mapsto \int_{\gG^x} \left\langle \xi^<(\gamma),\  \gamma \xi^>(\gamma^{-1})\right\rangle_{E_{r(\gamma)}} \rmd \lambda^x(\gamma) \right).
\end{eqnarray*}
This map extends to a contractive bracket from $\mH^<(\gG,E^<) \times \mH^>(\gG,E^>)$ to $\ContSect_0(X,B)$ which makes $\mH(\gG,E):= \left(\mH^<(\gG,E^<), \mH^>(\gG,E^>)\right)$ a $\Cont_0(X)$-Banach $\ContSect_0(X,B)$-pair. If $E$ is non-degenerate, then so is $\mH(\gG,E)$.

Note that the $\Cont_0(X)$-structures on $\mH^<(\gG,E^<)$ and $\mH^>(\gG,E^>)$ are not the same in general: on the left-hand side it is induced by the range map $r$, on the right-hand side by the source map $s$. This implies that the fibre of $\mH^<(\gG, E^<)$ over some $x\in X$ should be regarded as a completion of $\ContSect_c(\gG^x, E^<_x)$, whereas the fibre of $\mH^>(\gG, E^>)$ over $x$ should be regarded as a completion of $\ContSect_c(\gG_x,\ (r^*E)\restr_{\gG_x})$.

Assume now that $\mH(\gG)$ has also the following properties
\begin{itemize}
\item [{\bf (H2)}]
\begin{eqnarray*}
\forall \chi, \xi^<\in \Cont_c(\gG):\ && \norm{\xi^< * \chi}_{\mH^<(\gG)} \leq  \norm{\xi^<}_{\mH^<(\gG)}\norm{\chi}_{\mA(\gG)} \quad \text{and}\\
\forall \chi, \xi^>\in \Cont_c(\gG):\ && \norm{ \chi * \xi^>}_{\mH^>(\gG)} \leq \norm{\chi}_{\mA(\gG)} \norm{\xi^>}_{\mH^>(\gG)}.
\end{eqnarray*}
\end{itemize}

\noindent Let $A$ be another $\gG$-Banach algebra and let $E$ be a $\gG$-Banach $A$-$B$-pair. For all $a\in \ContSect_c(\gG,\ r^*A)$, all $\xi^<\in \ContSect_c(\gG,\ r^*E^<)$ and all $\xi^> \in \ContSect_c(\gG,\ r^*E^>)$, define
\[
(a \ \xi^>)(\gamma) = (a * \xi^>)(\gamma) = \int_{\gG^{r(\gamma)}} a (\gamma')\ \gamma'\xi^>(\gamma'^{-1}\gamma) \rmd \lambda^{r(\gamma)} (\gamma')
\]
and
\[
(\xi^< \ a)(\gamma) = (\xi^< * a )(\gamma) = \int_{\gG^{r(\gamma)}} \xi^<(\gamma')\ \gamma'a (\gamma'^{-1}\gamma) \rmd \lambda^{r(\gamma)} (\gamma')
\]
for all $\gamma\in \gG$. These actions lift to actions of $\mA(\gG,A)$ on $\mH^>(\gG,E^>)$ and $\mH^<(\gG,E^<)$, respectively. Equipped with them, $\mH(\gG,E)$ becomes a $\Cont_0(X)$-Banach $\ContSect_0(X,B)$-pair on which $\mA(\gG,A)$ by elements of $\Lin^{\Cont_0(X)}_{\ContSect_0(X,B)} \left(\mH(\gG,E)\right)$.

\begin{proposition} \label{Proposition:ActionOfAgGonHgGLocallyCompact}
If $\ContSect(X,A)$ acts on $E$ by locally compact operators and $\gG$ is proper, then $\mA(\gG,A)$ acts on $\mH(\gG,E)$ by locally compact operators.
\end{proposition}
\begin{proof}
Let $a \in \ContSect_c(\gG,\ r^*A)$ have compact support. If we can show that the action of $a$ on $\mH(\gG,E)$, denoted by $\pi(a)\in \Lin_{\ContSect_0(X,B)} \left(\mH(\gG,E)\right)$, is locally compact, then we are done. Let $\chi \in \Cont_c(X)$. We have to show that $\chi \pi(a)$ is compact. It was shown in \cite{Paravicini:07}, Appendix~E.8.3, that an operator which is given by a compact kernel with compact support is compact. We thus prove that $\chi \pi(a)$ is such an operator. Define
\[
k_{(\gamma_1,\gamma_2)} := \chi(s(\gamma_1))\pi_A(a(\gamma_2)) \in \Lin_{B_{r(\gamma_1)}} \left(E_{r(\gamma_1)}\right)
\]
for all $(\gamma_1,\gamma_2) \in \gG *_{r,r} \gG$. Then the action of $\chi \pi(a)$  on $\ContSect_c\left(\gG,\ r^*E^>\right)$ is given by
\begin{eqnarray*}
(\chi \pi(a))^>(\xi^>)(\gamma) &=&\chi(s(\gamma)) \ \int_{\gG^{r(\gamma)}}  a(\gamma') \gamma'\xi^>(\gamma'^{-1}\gamma) \rmd \lambda^{r(\gamma)}(\gamma')\\
&=& \int_{\gG^{r(\gamma)}}  \chi(s(\gamma)) a(\gamma') \gamma'\xi^>(\gamma'^{-1}\gamma) \rmd \lambda^{r(\gamma)}(\gamma')\\
&=& \int_{\gG^{r(\gamma)}}  k_{(\gamma,\gamma')} \gamma'\xi^>(\gamma'^{-1}\gamma) \rmd \lambda^{r(\gamma)}(\gamma')
\end{eqnarray*}
for all $\xi^>\in \ContSect_c\left(\gG,\ r^*E^>\right)$ and $\gamma\in \gG$. A similar calculation for the left-hand side shows that $\chi \pi(a)$ is indeed given by the kernel $k$.

The field of operators $\left(\pi_A(a(\gamma_2))\right)_{(\gamma_1,\gamma_2)\in \gG*_{r,r} \gG}$ is locally compact, so the same is true for $k$. Moreover, the support of $k$ is compact: Since $\gG$ is proper, the set $K:=\{\gamma\in \gG:\ r(\gamma)\in \supp \chi,\ s(\gamma) \in r(\supp a)\}$ is compact. Let $(\gamma_1,\gamma_2)\in \gG *_{r,r} \gG$. Then $k_{(\gamma_1,\gamma_2)} \neq 0$ implies $\gamma_1 \in K$ and $\gamma_2 \in \supp a$. So $(\gamma_1 ,\gamma_2)$ is contained in $K\times \supp a$. Hence $k$ has compact support.
\end{proof}

\noindent As a corollary of Proposition~\ref{Proposition:ActionOfAgGonHgGLocallyCompact} and because $B$ is non-degenerate, we get:

\begin{corollary}\label{Corollary:ActionOnmHLocallyCompact}
If $\gG$ is proper, then $\mA(\gG,B)$ acts on $\mH(\gG,B)$ by locally compact operators.
\end{corollary}

\noindent We now want to put an action of $\gG$ on $\mH(\gG,E)$. Technically, we have to replace $\mH(\gG,E)$ with the u.s.c.~field $\FieldPur(\mH(\gG,E))$ of pairs over $X$, compare Paragraph~\ref{Subsection:ComparisonKKbanRKKban}. So it is a natural to assume:
\begin{itemize}
\item [{\bf (H3)}] The $\Cont_0(X)$-Banach space $\mH^<(\gG)$ is a locally $\Cont_0(X)$-convex with respect to the $\Cont_0(X)$-action induced by $r$ and $\mH^>(\gG)$ is locally $\Cont_0(X)$-convex with respect to the action induced by $s$.
\end{itemize}

\noindent For all $\gamma\in \gG$, define a map $\alpha^<_{\gamma}$ from $\Cont_c\left(\gG^{s(\gamma)}\right)$ to $\Cont_c\left(\gG^{r(\gamma)}\right)$ by
\[
\chi^< \mapsto \alpha^<_{\gamma} (\chi^<)=\gamma \chi^< = \left(\gamma' \mapsto \chi^<(\gamma^{-1} \gamma') \right)
\]
and a map $\alpha^>_{\gamma}$ from $\Cont_c\left(\gG_{s(\gamma)}\right)$ to $\Cont_c\left(\gG_{r(\gamma)}\right)$ by
\[
\chi^> \mapsto \alpha^>_{\gamma}(\chi^>)= \gamma \chi^> = \left(\gamma' \mapsto \chi^>(\gamma'\gamma) \right).
\]
To get an action of $\gG$ on $\FieldPur(\mH(\gG,E))$ we have to assume that $\alpha^<$ and $\alpha^>$ are families of isometric maps, i.e., if we have that
\begin{itemize}
\item[{\bf (H4)}] \begin{eqnarray*}
\forall \chi^< \in \Cont_c(\gG^{s(\gamma)})\ \forall \gamma\in \gG:\ && \norm{\gamma \chi^<}_{\mH^<\left(\gG^{r(\gamma)}\right)} = \norm{\chi^<}_{\mH^<\left(\gG^{s(\gamma)}\right)} \qquad \text{and}\\
\forall \chi^>\in \Cont_c(\gG_{s(\gamma)})\ \forall \gamma\in \gG: && \norm{\gamma \chi^>}_{\mH^>\left(\gG_{r(\gamma)}\right)} = \norm{\chi^>}_{\mH^>\left(\gG_{s(\gamma)}\right)}.
\end{eqnarray*}
\end{itemize}

\noindent Note that all the examples of \ref{Examples:PairsOfMonotoneCompletions} satisfy (H3) and (H4).

The following result is proved in \cite{Paravicini:07}, 7.3.12, the proof involves some additional technical constructions which we prefer to omit here.

\begin{defprop}
Let $E$ be a $\gG$-Banach $B$-pair. Define
\[
\alpha^<_{\gamma} \colon \ContSect_c\big(\gG^{s(\gamma)},\ r^*E^<\big) \to \ContSect_c\big(\gG^{r(\gamma)},\ r^*E^<\big),\ \xi^< \mapsto \gamma \xi^< := \left(\gamma' \mapsto \gamma \xi^<(\gamma^{-1} \gamma')\right),
\]
and
\[
\alpha^>_{\gamma} \colon \ContSect_c\big(\gG_{s(\gamma)},\ r^*E^>\big) \to \ContSect_c\big(\gG_{r(\gamma)},\ r^*E^>\big),\ \xi^> \mapsto \gamma \xi^> := \left(\gamma' \mapsto \xi^>(\gamma'\gamma)\right),
\]
for all $\gamma\in \gG$. Then $\alpha^<_{\gamma}$ and $\alpha^>_{\gamma}$ are isometric for all $\gamma\in \gG$ and extend to isometric isomorphisms $\mH^<(\gG^{s(\gamma)},\ r^*E^<)\to \mH^<(\gG^{r(\gamma)},\ r^*E^<)$ and $\mH^>(\gG_{s(\gamma)},\ r^*E^>)\to \mH^>(\gG_{r(\gamma)},\ r^*E^>)$, respectively. The field $(\alpha^<_{\gamma},\ \alpha^>_{\gamma})_{\gamma\in \gG}$ is a continuous field of isomorphisms making $\FieldPur(\mH(\gG,E))$ a $\gG$-Banach $B$-pair.
\end{defprop}

\noindent Now that we assume that $\mH(\gG)$ does not only satisfy (H1) and (H2) but also (H3) and (H4), we can refine Proposition~\ref{Proposition:ActionOfAgGonHgGLocallyCompact} as follows:

\begin{proposition} Let $E$ be a $\gG$-Banach $A$-$B$-pair. Then $\FieldPur(\mH(\gG,E))$ is a $\gG$-Banach $B$-pair on which $\mA(\gG,A)$ acts by bounded \emph{$\gG$-equivariant} fields of linear operators. If $\gG$ is proper and $\ContSect(X,A)$ acts on $E$ by locally compact operators, then the action of $\mA(\gG,A)$ on $\FieldPur(\mH(\gG,E))$ is by \emph{$\gG$-equivariant} bounded locally compact fields of operators.
\end{proposition}

\noindent Because $B$ is non-degenerate, $\ContSect(X,B)$ acts on $B$ by locally compact operators. Hence we have:
\begin{corollary}
If $\gG$ is proper, then $\mA(\gG,B)$ acts on $\FieldPur{\mH(\gG,B)}$ by locally compact $\gG$-equivariant operators.
\end{corollary}

\noindent To finish this section, we state and prove an extension result which shows that the bracket on $\mH(\gG, E)$ is the restriction of the convolution product. This fact can be used to show that certain algebras are hereditary, see Lemma~\ref{Lemma:CompactSectionHereditaryInCompletion}.

\begin{proposition}\label{Proposition:ExtensionOfBracketToConvolution} Let $E$ be a $\gG$-Banach $B$-module. Then the convolution
\begin{eqnarray*}
\ContSect_c(\gG,\ r^*E^<) \times \ContSect_c(\gG,\ r^*E^>) &\to& \ContSect_c(\gG,r^*B),\\ (\xi^<,\xi^>)&\mapsto & \xi^<*\xi^>= \left(\gamma \mapsto \int_{\gG^{r(\gamma)}} \left\langle \xi^<(\gamma'),\  \gamma' \xi^>(\gamma'^{-1}\gamma)\right\rangle_{E_{r(\gamma)}} \rmd \lambda^{r(\gamma)}(\gamma') \right)
\end{eqnarray*}
extends to a contractive bilinear map
\[
\mH^<(\gG,E^<) \times \mH^>(\gG,E^>) \to \ContSect_0(\gG,r^*B),
\]
(also) written as a convolution product, such that the bracket on $\mH(\gG,E)$ is the composition of this map and the restriction map from $\ContSect_0(\gG, r^*B)$ to $\ContSect_0(X,B)$.
\end{proposition}
\begin{proof}
Let $\xi^<\in \ContSect_c(\gG,r^*E^<)$ and $\xi^>\in \ContSect_c(\gG,r^*E^>)$. For all $\gamma\in \gG$, we have
\[
(\xi^<*\xi^>)(\gamma) = \left\langle\xi^<_{r(\gamma)},\ \gamma \xi^>_{s(\gamma)}\right\rangle_{r(\gamma)}
\]
and hence
\begin{eqnarray*}
\norm{(\xi^<*\xi^>)(\gamma)} &=& \norm{\left\langle\xi^<_{r(\gamma)},\ \gamma \xi^>_{s(\gamma)}\right\rangle_{r(\gamma)}} \leq \norm{\xi^<_{r(\gamma)}}_{\mH^<(\gG,E^<)_{r(\gamma)}} \norm{\gamma \xi^>_{s(\gamma)}}_{\mH^>(\gG,E^>)_{r(\gamma)}}\\
&=&\norm{\xi^<_{r(\gamma)}}_{\mH^<(\gG,E^<)_{r(\gamma)}}\norm{\xi^>_{s(\gamma)}}_{\mH^>(\gG,E^>)_{s(\gamma)}}\leq \norm{\xi^<}_{\mH^<(\gG,E^<)}\norm{\xi^>}_{\mH^>(\gG,E^>)},
\end{eqnarray*}
because $\mH^>(\gG)$ satisfies (H4). Hence the convolution is continuous with norm $\leq 1$ and extends to a map $\mH^<(\gG,E^<) \times \mH^>(\gG,E^>) \to \ContSect_0(\gG,r^*B)$ with the desired properties.
\end{proof}

\subsection{Regular unconditional completions}\label{Subsection:RegularCompletions}

\noindent For simplicity, we introduce the following abbreviation:

\begin{definition} An unconditional completion $\mA(\gG)$ of $\Cont_c\left(\gG\right)$ is said to be \demph{regular} if there exists a pair $\mH(\gG)$ of monotone completions of $\Cont_c(\gG)$ which satisfies (H1)-(H4).
\end{definition}

\noindent Note that there might exist many different such pairs of monotone completions on which a regular unconditional completion acts, the important part of the definition really is the existence of such a pair, not its particular shape.

\begin{examples} Most examples of unconditional completions that we have come across so far are regular for rather obvious reasons:
\begin{enumerate}
\item The unconditional completion $\Leb^1\left(\gG\right)$ acts on the pair $\left(\Leb^1\left(\gG\right),\ \Cont_0\left(\gG\right)\right)$.

\item The symmetrised version $\Leb^1\left(\gG\right) \cap \Leb^1\left(\gG\right)^*$ is also {regular} because the norm defining it dominates the norm $\norm{\cdot}_1$. Moreover, it acts on the pair $\left(\Leb^2_r\left(\gG\right),\ \Leb^2_s\left(\gG\right)\right)$ (see \cite{Renault:80}). It should not be too hard to check that it also acts on $\left(\Leb^{p'}_r\left(\gG\right),\ \Leb^{p}_s\left(\gG\right)\right)$ for all $p,p'\in ]1,\infty[$ such that $\frac{1}{p} + \frac{1}{p'} =1$.

\item The completion $\mA_{\max}\left(\gG\right)$ acts on $\left(\Leb^2_r\left(\gG\right),\ \Leb^2_s\left(\gG\right)\right)$ by definition; see Section~3 of \cite{Lafforgue:06}.

\item If $G$ is a locally compact Hausdorff group acting on some locally compact Hausdorff space $X$, then $\Leb^1\left(G,\ \Cont_0(X)\right)$ is a {regular} completion of $\Cont_c\left(G\ltimes X\right)$ because its norm dominates the norm of the {regular} completion $\Leb^1\left(G\ltimes X\right)$.
\end{enumerate}
\end{examples}

\noindent Regularity is essential in our construction of the homomorphism $\GreenJulgAufstiegKKban{\mA}{B}$ down below. It also makes some arguments in the next chapter simpler (but might perhaps be avoided in some instances).

\subsection{The (inverse) homomorphism $\GreenJulgAufstiegKKban{\mA}{B}$}\label{Subsection:GreenJulgAufstieg}

\emph{In \ref{Subsection:GreenJulgAufstieg}, let $\gG$ be proper and $\mA(\gG)$ be regular.}

Recall that we used the name $\pi$ for the canonical projection from $X$ to $X/\gG$. Let $\pi$ also denote the map from $\gG$ to $X/\gG$ that maps $\gamma$ to $\pi(r(\gamma)) = \pi(s(\gamma))$ (which extends $\pi\colon X \to X/\gG$). If we regard $X/\gG$ as a locally compact Hausdorff groupoid, then the map $\pi \colon \gG \to X/\gG$ is actually a strict morphism of groupoids. If $E$ is a u.s.c.~field of Banach spaces over $X/\gG$, then $\pi^*E$ is a $\gG$-Banach space (with a rather trivial action).\footnote{See \cite{Lafforgue:06} or \cite{Paravicini:07:Induction:arxiv}, Section~4.1, for a definition of the pullback along strict morphisms.}  If $T$ is a continuous field of linear maps between u.s.c.~fields of Banach spaces over $X/\gG$, then $\pi^*T$ is an $\gG$-equivariant continuous field of linear maps between $\gG$-Banach spaces. We use these facts to define our ``inverse homomorphism'':

\begin{enumerate}
\item The first step is the map\footnote{See Subsection~\ref{Subsection:ComparisonKKbanRKKban} or \cite{Paravicini:07}, Chapter~4, for a definition of the functor $\Field{\cdot}$.} $\Field{\cdot}$ which yields a homomorphism
\[
\Field{\cdot}\colon \RKKban\left(\Cont_0(X/\gG);\ \Cont_0(X/\gG),\ \mA(\gG,B)\right) \to \KKbanW{X/\gG} \left(\C_{X/\gG},\ \Field{\mA(\gG,B)}\right).
\]

\item The second step is the pullback homomorphism along $\pi$:
\[
\pi^*\colon \KKbanW{X/\gG} \left(\C_{X/\gG},\ \Field{\mA(\gG,B)}\right)\to \KKbanW{\gG} \left(\C_{X},\ \pi^*\Field{\mA(\gG,B)}\right).
\]
Note that this homomorphism, on the level of cycles, produces cycles with $\gG$-equivariant operator.

\item Pick a pair $\mH(\gG)$ of monotone completions of $\Cont_c(\gG)$ satisfying (H1) - (H4). Note that there is a canonical action of $\pi^*\Field{\mA(\gG,B)}$ on $\Field{\mH(\gG,B)}$.

By Corollary~\ref{Corollary:ActionOnmHLocallyCompact}, the algebra $\mA(\gG,B)$ acts on $\Field{\mH(\gG,B)}$ by locally compact operators. If $\chi \in \Cont_c(X)$ and $a \in \mA(\gG,B)$, then $x \mapsto \chi(x) a_{\pi(x)}$ is a section of $\pi^*(\Field{\mA(\gG,B)})$ with compact support. Such sections act on $\Field{\mH(\gG,B)}$ by locally compact operators with compact support, so they act by compact operators. By a density argument, all sections of $\pi^*(\Field{\mA(\gG,B)})$ that vanish at infinity act on $\Field{\mH(\gG,B)}$ by compact operators. Hence we can regard $\Field{\mH(\gG,B)}$ as a Morita cycle\footnote{ See \cite{Paravicini:07:Induction:arxiv}, Paragraph~3.3.6, for a definition; compare \cite{Paravicini:07:Morita:erschienen}, Definition~5.7.} from $\pi^*\Field{\mA(\gG,B)}$ to $B\cong \Field{\ContSect_0(X,B)}$. The important point is that this Morita cycle carries an action of $\gG$ which makes it a $\gG$-equivariant Morita cycle. Morita cycles act on $\KKban$ from the right, so we get a homomorphism
\[
\otimes_{\pi^*\Field{\mA(\gG,B)}} \Field{\mH(\gG,B)}\colon \KKbanW{\gG} \left(\C_{X},\ \pi^*\Field{\mA(\gG,B)}\right) \to \KKbanW{\gG} \left(\C_{X},\ B\right).
\]
If a cycle has a $\gG$-equivariant operator, then it stays equivariant under this homomorphism.
\end{enumerate}

\noindent The composition of these three homomorphisms gives the desired natural homomorphism
\[
\GreenJulgAufstiegEban{\mA}{\mH}{B}\colon \RKKban\left(\Cont_0\left(X/\gG\right);\ \Cont_0\left(X/\gG\right),\ \mA(\gG,B)\right) \to \KKbanW{\gG}\left(\C_X,\ B\right)
\]
which produces cycles with $\gG$-equivariant operators.

\begin{proposition}
Let $\mH'(\gG)=\left(\mH'^<(\gG),\ \mH'^>(\gG)\right)$ be another pair of monotone completions of $\Cont_c(\gG)$ satisfying (H1) - (H4). Then the natural homomorphisms $\GreenJulgAufstiegEban{\mA}{\mH}{B}$ and $\GreenJulgAufstiegEban{\mA}{\mH'}{B}$ are equal. We call this natural homomorphism $\GreenJulgAufstiegKKban{\mA}{B}$.
\end{proposition}
\begin{proof}
We first consider the case that $\norm{\cdot}_{\mH^<} \leq \norm{\cdot}_{\mH'^<}$ and $\norm{\cdot}_{\mH^>} \leq \norm{\cdot}_{\mH'^>}$. In this case, we have a canonical homomorphism $\Phi$ from $\mH'(\gG,B)$ to $\mH(\gG,B)$ which gives us an equivariant homomorphism $\Field{\Phi}$ from $\Field{\mH'(\gG,B)}$ to $\Field{\mH(\gG,B)}$.  The homomorphism $\Field{\Phi}$ is actually a morphism of equivariant Morita cycles from $\pi^*\Field{\mA(\gG,B)}$ to $B$. A careful revision of the proof that $\pi^*\Field{\mA(\gG,B)}$ acts by compact operators on $\Field{\mH'(\gG,B)}$ and on $\Field{\mH(\gG,B)}$ shows that $\Field{\Phi}$ satisfies the conditions of Theorem~3.20 of \cite{Paravicini:07:Induction:arxiv} and hence induces a homotopy from $\Field{\mH'(\gG,B)}$ to $\Field{\mH(\gG,B)}$. So $\GreenJulgAufstiegEban{\mA}{\mH'}{B} = \GreenJulgAufstiegEban{\mA}{\mH}{B}$ because the tensor product with Morita cycles lifts to homotopy classes.\footnote{See \cite{Paravicini:07:Induction:arxiv}, Paragraph 3.3.6, or \cite{Paravicini:07}, Section 3.8.} 

Now consider the general case. By taking the maximum of the norms on $\mH^<(\gG)$ and $\mH'^<(\gG)$ we define a monotone completion $\mH''^<(\gG)$ of $\Cont_c(\gG)$; similarly, we define $\mH''^>(\gG)$. The pair $\mH''(\gG) := \left(\mH''^<(\gG),\ \mH''^>(\gG)\right)$ also satisfies (H1) - (H4). By the first part of the proof we can conclude $\GreenJulgAufstiegEban{\mA}{\mH}{B} = \GreenJulgAufstiegEban{\mA}{\mH''}{B} = \GreenJulgAufstiegEban{\mA}{\mH'}{B}$.
\end{proof}

\subsection{$\GreenJulgAbstiegKKban{\mA}{B}\circ \GreenJulgAufstiegKKban{\mA}{B}  = \id$ on the level of $\KKban$} \label{SubsectionMainTheoremSurjectivity}

We now prove Theorem~\ref{Theorem:NaturalSplitPartOfGreenJulg}. \emph{Let $\gG$ be proper and let $\mA(\gG)$ be regular. Assume moreover that $\gG$ admits a cut-off function.}

\subsubsection{Idea of the proof}

First, we choose a pair $\mD(X)$ of monotone completions of $\Cont_c(X)$ satisfying (D1) and (D2). Secondly, because $\mA(\gG)$ is {regular}, we can also choose a pair $\mH(\gG) = \left(\mH^<(\gG),\ \mH^>(\gG)\right)$ of monotone completions of $\Cont_c(\gG)$ satisfying (H1) - (H4). Let $(E,T) \in \Eban\left(\Cont_0(X/\gG);\ \Cont_0(X/\gG),\ \mA\left(\gG,B\right)\right)$. We have to show that $(E,T)$ is homotopic to $\GreenJulgAbstiegEban{\mA}{\mD}{B} (\GreenJulgAufstiegEban{\mA}{\mH}{B}(E,T))$. The obvious strategy is to define a morphism from $\GreenJulgAbstiegEban{\mA}{\mD}{B} (\GreenJulgAufstiegEban{\mA}{\mH}{B}(E))$ to $E$ which induces a homotopy; there is a canonical candidate for such a morphism defined on a dense subspace, but this candidate does not extend to a continuous morphism on the entire space: The norms on $\GreenJulgAbstiegEban{\mA}{\mD}{B} (\GreenJulgAufstiegEban{\mA}{\mH}{B}(E))$ and $E$ seem to be difficult to compare in general.

We overcome this problem by constructing a pair $\tilde{E}:=(\tilde{E}^<,\tilde{E}^>)$ of $\C$-vector spaces which are equipped with compatible $\Cont_c(X/\gG)$-module structures and left/right $\ContSect_c\left(\gG,\ r^*B\right)$-module structures and a bilinear map from $\tilde{E}^< \times \tilde{E}^>$ to $\ContSect_c\left(\gG,\ r^*B\right)$. On this pair, which could be called a ``pre-$\mA(\gG,B)$-pair'', we construct a pair of formally adjoint operators $\tilde{T}$. Moreover, we define canonical ``homomorphisms'' $\Phi_E$ from $\tilde{E}$ to $E$ and $\Psi_E$ from $\tilde{E}$ to $\GreenJulgAbstiegEban{\mA}{\mD}{B} (\GreenJulgAufstiegEban{\mA}{\mH}{B}(E))$ which intertwine $\tilde{T}$ and $T$ and $\GreenJulgAbstiegEban{\mA}{\mD}{B} (\GreenJulgAufstiegEban{\mA}{\mH}{B} (T))$, respectively:
\[
\xymatrix{
&(\tilde{E},\tilde{T})\ar[ld]_{\Phi_E}\ar[rd]^{\Psi_E}&\\
\quad (E,T)\quad && \GreenJulgAbstiegEban{\mA}{\mD}{B} (\GreenJulgAufstiegEban{\mA}{\mH}{B} (E,T))
}
\]
One can think of $\tilde{E}$ as a dense subspace of both, $E$ and $\GreenJulgAbstiegEban{\mA}{\mD}{B} (\GreenJulgAufstiegEban{\mA}{\mH}{B} (E))$. Now we put on $\tilde{E}$ the supremum of the semi-norms which are induced by the two homomorphisms, making the homomorphisms continuous. We then show that the completion of $\tilde{E}$ together with the continuous extension of $\tilde{T}$ is in $\Eban\left(\Cont_0(X/\gG);\ \Cont_0(X/\gG),\ \mA\left(\gG,B\right)\right)$ and that the two homomorphisms induce homotopies. Hence also $(E,T)$ and $\GreenJulgAbstiegEban{\mA}{\mD}{B} (\GreenJulgAufstiegEban{\mA}{\mH}{B} (E,T))$ are homotopic.

\subsubsection{The construction of $\tilde{E}$, $\Phi_E$ and $\Psi_E$}

We are going to cut the proof into a series of statements and definitions. In this subsection, let $E$ and $F$ be $\Cont_0(X/\gG)$-Banach $\mA(\gG,B)$-pairs.

\noindent {\bf The pair $\tilde{E}$:} Define
\[
\tilde{E}^>:= E^>\otimes_{\ContSect_c\left(\gG,\; r^*B\right)} \ContSect_c\left(\gG,\; r^*B\right)
\]
and
\[
\tilde{E}^<:= \ContSect_c\left(\gG,\; r^*B\right)\otimes_{\ContSect_c\left(\gG,\; r^*B\right)} E^<.
\]
These vector spaces carry canonical and compatible actions of $\ContSect_c\left(\gG,\ r^*B\right)$ and $\Cont_c\left(X/\gG\right)$. A bracket on $\tilde{E}$ is defined by
\begin{eqnarray*}
\langle \cdot,\cdot\rangle\colon \tilde{E}^<\times \tilde{E}^> &\to& \ContSect_c\left(\gG,\; r^*B\right),\\
\left\langle \beta^< \otimes e^<, \ e^> \otimes \beta^>\right\rangle (\gamma) &:=& \beta^< * \langle e^<,e^>\rangle *\beta^> = \left\langle \beta^<e^<,\ e^>\beta^>\right\rangle.
\end{eqnarray*}

\noindent We check that the bracket has indeed its values in $\ContSect_c(\gG,r^*B)$: The element $\langle e^<,e^>\rangle$ is in $\mA(\gG,B)$ by definition, and we now show that the product $\beta^< * \beta * \beta^>$ is in $\ContSect_c(\gG,B)$ for all $\beta^<,\beta^> \in \ContSect_c(\gG,B)$ and $\beta\in \mA(\gG,B)$. If we regard $\beta^<$ as an element of $\mH^<(\gG,B)$ and $\beta^>$ as an element of $\mH^>(\gG,B)$, then we can conclude from Proposition~\ref{Proposition:ExtensionOfBracketToConvolution} that the map $\beta\mapsto \beta^< * \beta * \beta^>$ is continuous from $\mA(\gG, B)$ to $\ContSect_0(\gG,B)$ because $\mA(\gG)$ acts on $\mH(\gG)$. Moreover, the support of the product $\beta^< * \beta *\beta^>$ is always contained in the set $\{\gamma\in \gG:\ r(\gamma) \in r(\supp \beta^<),\ s(\gamma)\in s(\supp \beta^>)\}$, which is compact because $\gG$ is proper. \footnote{Compare the proof of Lemma~\ref{Lemma:CompactSectionHereditaryInCompletion}.}

\noindent {\bf The map $\Phi_E$:} Define
\[
\Phi_E^>\colon \tilde{E}^> \to E^>,\ e^> \otimes \beta^> \mapsto e^>\beta^>
\]
and
\[
\Phi_E^<\colon \tilde{E}^< \to E^<,\ \beta^< \otimes e^< \mapsto \beta^< e^<.
\]

\noindent Both maps are clearly $\ContSect_c\left(\gG,\; r^*B\right)$- and $\Cont_c\left(X/\gG\right)$-linear. The pair $\Phi_E=\left(\Phi_E^<,\Phi_E^>\right)$ is compatible with the brackets on $\tilde{E}$ and $E$.

\noindent {\bf  The map $\Psi_E$:} Let $e^> \in E^>$ and $\beta^>\in \ContSect_c\left(\gG,r^*B\right)$. Since $\beta^>$ has compact support, the function $x\mapsto (e^>\otimes \beta^>)_x = e^>_{\pi(x)}\otimes \beta^>_x$ is in $\ContSect_c\left(X,\ \pi^*\Field{E^>} \otimes_{\pi^*\Field{\mA(\gG,B)}} \Field{\mH^>(\gG,B)} \right)$; we can regard this function as an element $\Psi_E^>\left(e^> \otimes \beta^>\right)$ of $\mD^>\left(X,\ \pi^*\Field{E^>} \otimes_{\pi^*\Field{\mA(\gG,B)}} \Field{\mH^>(\gG,B)}\right)$; here $\pi\colon X\to X/\gG$ denotes the canonical projection.  This gives rise to a map $\Psi^>_E$ from $\tilde{E}^>$ to $\GreenJulgAbstiegEban{\mA}{\mD}{B} (\GreenJulgAufstiegEban{\mA}{\mH}{B}(E))^>$.

\noindent Similarly we define
\[
\Psi_E^<\left(\beta^< \otimes e^<\right)_x:= \beta^<_x\otimes e^<_{\pi(x)}  \in \mH^<(\gG,B)_x \otimes_{\mA(\gG,B)_{\pi(x)}} E_{\pi(x)}^<
\]
for all $e^<\in E^<$, $\beta^<\in \ContSect_c\left(\gG,r^*B\right)$ and $x\in X$, giving us a $\ContSect_c(\gG, r^*B)$-linear and $\Cont_c(X/\gG)$-linear map $\Psi^<_E$ from $\tilde{E}^<$ to $\GreenJulgAbstiegEban{\mA}{\mD}{B} (\GreenJulgAufstiegEban{\mA}{\mH}{B}(E))^<$.

The pair $\Psi_E=\left(\Psi_E^<, \Psi_E^>\right)$ is compatible with the brackets on $\tilde{E}$ and $\GreenJulgAbstiegEban{\mA}{\mD}{B} (\GreenJulgAufstiegEban{\mA}{\mH}{B}(E))$.

\noindent {\bf The constructions for linear operators:} Let $S\in \Lin_{\mA(\gG,B)}(E,F)$ be an operator between the $\Cont_0(X/\gG)$-Banach $\mA(\gG,B)$-pairs $E$ and $F$. Define
\[
\tilde{S}^>\colon \tilde{E}^>\to \tilde{F}^>,\ \xi^> \otimes \beta^> \mapsto S^>(\xi^>) \otimes \beta^>
\]
and
\[
\tilde{S}^<\colon \tilde{F}^<\to \tilde{E}^<,\ \beta^< \otimes \xi^< \mapsto \beta^< \otimes S^<(\xi^<).
\]

\noindent Note that $\tilde{S}:=\left(\tilde{S}^<,\tilde{S}^>\right)$ is formally adjoint in the following sense:
\begin{eqnarray*}
&& \left\langle \tilde{S}^<\left(\beta^< \otimes \xi^<\right),\ \xi^>\otimes \beta^>\right\rangle = \beta^< * \left\langle S^<(\xi^<), \ \xi^>\right\rangle * \beta^>\\ &=&\beta^< * \left\langle \xi^<, \ S^>(\xi^>)\right\rangle * \beta^> =\left\langle \beta^< \otimes \xi^<,\ \tilde{S}^>\left(\xi^>\otimes \beta^>\right)\right\rangle
\end{eqnarray*}
for all $\beta^<,\beta^>\in \ContSect_c(\gG,r^*B)$, $\xi^<\in \ContSect_c(X,F^<)$ and $\xi^>\in \ContSect_c(X,E^>)$.

By direct calculation one checks:
\begin{enumerate}
\item The maps $\Phi_E$ and $\Phi_F$ intertwine $\tilde{S}$ and $S$ in the obvious sense.

\item The maps $\Psi_E$ and $\Psi_F$ intertwine $\tilde{S}$ and $\GreenJulgAbstiegEban{\mA}{\mD}{B} \left(\GreenJulgAufstiegEban{\mA}{\mH}{B}(S)\right)$.
\end{enumerate}

\subsubsection{Putting a norm on $\tilde{E}$}

If $\tilde{e}^>\in \tilde{E}^>$, then define
\[
\norm{\tilde{e}^>}:= \max\left\{\norm{\Phi_E^>(\tilde{e}^>)},\ \norm{\Psi_E^>(\tilde{e}^>)} \right\}.
\]
This is a semi-norm on $\tilde{E}^>$. Let $\overline{E}^>$ be the (Hausdorff-) completion of $\tilde{E}^>$ with respect to this semi-norm. In an analogous fashion, define a semi-norm on $\tilde{E}^<$ and call the completion $\overline{E}^<$. The actions of $\ContSect_c(\gG,r^*B)$ and $\Cont_c(X/\gG)$ on $\tilde{E}$ extend to non-degenerate actions of $\mA(\gG,B)$ and $\Cont_0(X/\gG)$ on $\overline{E}$. The bracket on $\tilde{E}$ extends to a continuous bracket on $\overline{E}$, making $\overline{E}$ a $\Cont_0(X/\gG)$-Banach $\mA(\gG, B)$-pair.

Now the map $\Phi_E^>$ extends by continuity to a continuous linear map from $\overline{E}^>$ to $E$ which is $\mA(\gG,B)$- and $\Cont_0(X/\gG)$-linear. Similar things can be said about $\Phi^<_E$, $\Psi^>_E$ and $\Psi_E^<$. We get homomorphisms $\Phi_E$ from $\overline{E}$ to $E$ and $\Psi_E$ from $\overline{E}$ to $\GreenJulgAbstiegEban{\mA}{\mD}{B} (\GreenJulgAufstiegEban{\mA}{\mH}{B}(E))$.

Let $S\in \Lin_{\mA(\gG,B)}(E,F)$ as above. Then the map $\tilde{S}^>$ satisfies
\[
\norm{\tilde{S}^>(\tilde{e}^>)} \leq \norm{S^>} \norm{\tilde{e}^>}
\]
for all $\tilde{e}^>\in \tilde{E}^>$ and extends therefore to an operator $\overline{S}^>$ from $\overline{E}^>$ to $\overline{F}^>$. Analogously for $\tilde{S}^<$. We thus get an element $\overline{S} \in \Lin_{\mA(\gG,B)}\left(\overline{E},\overline{F}\right)$ of norm $\leq \norm{S}$. The map $S\mapsto \overline{S}$ is $\C$-linear and functorial. The homomorphisms $\Phi_E$ and $\Phi_F$ intertwine $\overline{S}$ and $S$ in the obvious sense and the homomorphisms $\Psi_E$ and $\Psi_F$ intertwine $\overline{S}$ and $\GreenJulgAbstiegEban{\mA}{\mD}{B} (\GreenJulgAufstiegEban{\mA}{\mH}{B}(S))$.

\noindent By direct comparison of the operators one can show:

\begin{lemma}\label{Lemma:JointCompactnessOfSandSbar} Let $e^<\in \ContSect_0(X,E^<)$, $f^>\in \ContSect_0(X,F^>)$, $\beta^<, \beta^>\in \ContSect_c(\gG,r^*B)$. If
\[
S= \ketbra{f^>\beta^>}{\beta^<e^<}\in\Komp_{\mA(\gG,B)}\left(E,F\right),
\]
then
\[
\overline{S}= \ketbra{f^>\otimes \beta^>}{\beta^<\otimes e^<}\in\Komp_{\mA(\gG,B)}\left(\overline{E},\overline{F}\right)
\]
and
\begin{eqnarray*}
\GreenJulgAbstiegEban{\mA}{\mD}{B} \left(\GreenJulgAufstiegEban{\mA}{\mH}{B}(S)\right) &=& \ketbra{\Psi_F^>(f^>\otimes \beta^>)}{\Psi_E^<(\beta^<\otimes e^<)}\\
&\in&\Komp_{\mA(\gG,B)}\left(\GreenJulgAbstiegEban{\mA}{\mD}{B} \left(\GreenJulgAufstiegEban{\mA}{\mH}{B}(E)\right), \GreenJulgAbstiegEban{\mA}{\mD}{B} \left(\GreenJulgAufstiegEban{\mA}{\mH}{B}(F)\right)\right).
\end{eqnarray*}
It follows for all $S\in \Komp_{\mA(\gG,B)} \left(E,F\right)$ that $\overline{S}$ and $\GreenJulgAbstiegEban{\mA}{\mD}{B} (\GreenJulgAufstiegEban{\mA}{\mH}{B}(S))$ are compact and that $(\overline{S},S) \in \Komp(\Phi_E,\Phi_F)$ as well as  $(\overline{S}, \ \GreenJulgAbstiegEban{\mA}{\mD}{B} (\GreenJulgAufstiegEban{\mA}{\mH}{B}(S)))\in \Komp\left(\Psi_E,\Psi_F\right)$. The precise definition of $\Komp\left(\Phi_E,\Phi_F\right)$ can be found in \cite{Paravicini:07:Morita:erschienen}; compare the discussion around Theorem~\ref{Theorem:SufficientCondition:RKKbanG}.
\end{lemma}

\subsubsection{The proof of $\GreenJulgAbstiegKKban{\mA}{B}\circ \GreenJulgAufstiegKKban{\mA}{B}  = \id$}

Let $(E,T) \in \Eban\left(\Cont_0(X/\gG);\ \Cont_0(X/\gG),\ \mA(\gG,B)\right)$. We show that $\left(\overline{E},\overline{T}\right)$ is homotopic to $(E,T)$ as well as to $\GreenJulgAbstiegEban{\mA}{\mD}{B} \left(\GreenJulgAufstiegEban{\mA}{\mH}{B}(E,T)\right)$.

If $\chi\in \Cont_c(X/\gG)$ and $S:=\chi (T^2-1)$, then $(\overline{S},S)$ is in $\Komp(\Phi_E,\Phi_E)$ and $\left(\overline{S}, \GreenJulgAbstiegEban{\mA}{\mD}{B} \left(\GreenJulgAufstiegEban{\mA}{\mH}{B}(S)\right)\right)\in \Komp\left(\Psi_E, \Psi_E\right)$ by Lemma~\ref{Lemma:JointCompactnessOfSandSbar}. If follows that $(\overline{E},\overline{T})$ is in $\Eban\left(\Cont_0(X/\gG);\ \Cont_0(X/\gG),\ \mA(\gG,B)\right)$ and, using Theorem~\ref{Theorem:SufficientCondition:RKKbanG}, that it is homotopic to $(E,T)$ as well as to $\GreenJulgAbstiegEban{\mA}{\mD}{B} \left(\GreenJulgAufstiegEban{\mA}{\mH}{B}(E,T)\right)$.

\subsection{Sketch of the proof of $\GreenJulgAufstiegKKban{\mA}{B}\circ \GreenJulgAbstiegKKban{\mA}{B}  = \id$} \label{Subsection:SketchOfProofOfGreenJulg}

We first have to introduce an additional technical concept to be able to formulate the precise conditions under which we can show the injectivity part of the generalised Green-Julg theorem:

Let $\mH(\gG)=\left(\mH^<(\gG),\ \mH^>(\gG)\right)$ be a pair of monotone completions of $\Cont_c(\gG)$ satisfying (H1) - (H4). A cut-off pair  $c=(c^<,c^>)$ for $\gG$ is called an \demph{$\mH(\gG)$-cut-off pair} if
\[
\forall x\in X: \ \Big\|\gG_x\ni \gamma \mapsto c^>(r(\gamma))\Big\|_{\mH^>(\gG_x)} = 1 \quad \wedge \quad  \Big\|\gG^x\ni \gamma \mapsto c^<(s(\gamma))\Big\|_{\mH^<(\gG^x)}  = 1.
\]

\begin{examples} Assume that $X/\gG$ is $\sigma$-compact. Let $c$ be a cut-off-function for $\gG$.
\begin{enumerate}
\item Proposition~\ref{Proposition:CutOffPairLOneCNull} gives a $\mH(\gG)$-cut-off pair $(c,d)$ for the pair $\mH(\gG)= \left(\Leb^1(\gG),\ \Cont_0\left(\gG\right)\right)$.
\item If $p,p'\in ]1,\infty[$ such that $\frac{1}{p} + \frac{1}{p'} =1$, then $\left(c^{\frac{1}{p'}},\ c^{\frac{1}{p}}\right)$ is a $\mH(\gG)$-cut-off pair for the pair $\mH(\gG)=\left(\Leb^{p'}_r\left(\gG\right),\ \Leb^p_s\left(\gG\right)\right)$.
\end{enumerate}
\end{examples}

\noindent The technical reason to consider $\mH(\gG)$-cut-off pairs is that they allow us to embed $\ContSect_0(X,E)$ into $\mH(\gG,E)$ as a direct summand, where $E$ is a $\gG$-Banach $B$-pair. This is a Banach algebraic analogue of the fact that every $\gG$-C$^*$-algebra $B$ can be embedded into $\Leb^2(\gG,B)$ as a direct summand, see Proposition~6.21 of \cite{Tu:99}.

We can now formulate our result:

\begin{theorem}[Generalised Green-Julg Theorem]\label{Theorem:GeneralisedGreenJulg}  Let $\gG$ be proper and let $\mA(\gG)$ be an unconditional completion of $\Cont_c(\gG)$ such that there exists a pair $\mH(\gG)$ of monotone completions of $\Cont_c(\gG)$ satisfying (H1) - (H4) and such that there exists an $\mH(\gG)$-cut-off pair for $\gG$. Then there is an isomorphism
\[
\GreenJulgAbstiegKKban{\mA}{B}\colon \KKbanW{\gG}(\Cont_0(X),\ B) \ \cong \ \RKKban (\Cont_0(X/\gG);\, \Cont_0(X/\gG),\, \mA(\gG, B)),
\]
natural in the non-degenerate $\gG$-Banach algebra $B$.
\end{theorem}

\noindent Note that, trivially, the hypotheses of the theorem imply that $\gG$ admits a cut-off function and that $\mA(\gG)$ is regular. Hence the surjectivity part of the theorem has already been settled.

\noindent {\bf Idea of the proof of $\GreenJulgAufstiegKKban{\mA}{B}\circ \GreenJulgAbstiegKKban{\mA}{B}  = \id$:}
Let $B$, $\gG$, $\mA(\gG)$ and $\mH(\gG)$ be as in the theorem and let $c$ be an $\mH(\gG)$-cut-off pair. We want to show that $\GreenJulgAufstiegKKban{\mA}{B}\circ \GreenJulgAbstiegKKban{\mA}{B}   = \id$ as an endomorphism of the group $\RKKban(\Cont_0(X/\gG);\ \Cont_0(X/\gG),\ \mA(\gG,B))$.

Let $(E,T) \in \EbanW{\gG}(\Cont_0(X),\ B)$ with $\gG$-equivariant $T$. The idea is to define a homomorphism $\Phi_E$ from $\GreenJulgAufstiegEban{\mA}{\mH}{B}(\GreenJulgAbstiegEban{\mA}{\mA^c}{B} (E))$ to $E$ that commutes with the operator $\GreenJulgAufstiegEban{\mA}{\mH}{B}(\GreenJulgAbstiegEban{\mA}{\mA^c}{B}(T))$ and $T$. We then show that $\Phi_E$ induces a homotopy by checking the technical conditions of Theorem~3.20 of \cite{Paravicini:07:Induction:arxiv}.  Note that we use the particular pair $\mA^c(X)$ of monotone completions of $\Cont_c(X)$ here, see Paragraph~\ref{Subsubsection:UniquenessAndExistenceOfDX}.

The central ingredient in the construction of $\Phi_E$ is a homomorphism
\[
\mA^c(X,E) \otimes_{\mA(\gG,B)} \mH(\gG,B) \to \ContSect_0(X,E).
\]
To define it, observe that the convolution gives a homomorphism
\[
\mA(\gG,E) \otimes_{\mA(\gG,B)} \mH(\gG,B) \to \mH(\gG, E).
\]
By definition, $\mA^c(X,E)$ embeds into $\mA(\gG,E)$, so we can embed $\mA^c(X,E) \otimes_{\mA(\gG,B)} \mH(\gG,B)$ into $\mA(\gG,E) \otimes_{\mA(\gG,B)} \mH(\gG,B)$. On the other hand, $\ContSect_0(X,E)$ is contained as a direct summand in $\mH(\gG,E)$ because $c$ is a $\mH(\gG)$-cut-off pair, so we can compose with the projection onto this summand to obtain the desired homomorphism. The homomorphism $\Phi_E$ is constructed from it by some standard operations.

The main difficulty of the proof is to check that the homomorphism $\Phi_E$ really gives a homotopy between $\GreenJulgAufstiegEban{\mA}{\mH}{B}(\GreenJulgAbstiegEban{\mA}{\mA^c}{B} (E,T))$ and $(E,T)$. This boils down to some approximation arguments for compact operators which are carried out in detail in \cite{Paravicini:07}, Section~7.8.

\section{The Bost conjecture and proper Banach algebras}\label{Section:BostConjectureAndProperAlgebras}

In this section, let $\gG$ be a locally compact Hausdorff groupoid equipped with a Haar system. Assume moreover that there is a locally compact classifying space $\uEgG$ for proper actions of $\gG$, which is then unique up to homotopy. Let $\mA(\gG)$ be an unconditional completion of $\Cont_c(\gG)$.

\subsection{The Bost conjecture}

If $B$ is a $\gG$-Banach algebra, then there is an obvious definition of a  $\gG$-Banach algebra $SB:= B]0,1[$ with fibres $B_g]0,1[$ for all $g\in \gG^{(0)}$.

\begin{definition} For every $\gG$-Banach algebra $B$, define the \demph{topological $\KTh$-theory} for $\gG$ and $B$ as
\[
\KTh^{\top2,\ban}_0\left(\gG,\ B\right):= \lim_{\to}\, \KKbanW{\gG}\left(\Cont_0(X),\ B\right),
\]
where $X$ runs through the closed proper $\gG$-compact subspaces of $\uEgG$. Define $\KTh^{\top2,\ban}_n\left(\gG,\ B\right):= \KTh^{\top2,\ban}_0\left(\gG,\ S^nB\right)$ for $n\in \N$.
\end{definition}

\noindent Note that, if $X$ is a locally compact Hausdorff left $\gG$-space (with anchor map $\rho$), then we would like to think of $\Cont_0(X)$ as a $\gG$-Banach space. We thus have to take the pushforward $\rho_* \C_X$ instead, which is a field over $\gG^{(0)}$.\footnote{See Paragraph~\ref{Subsubsection:ThePushforwardConstruction} for a sketch of this construction.} To make the notation more familiar, we nevertheless write $\Cont_0(X)$. For more details, consult \cite{Paravicini:07:Induction:arxiv} or \cite{Paravicini:07}.

If $B$ is a $\gG$-C$^*$-algebra, then there is a canonical homomorphism from the C$^*$-algebraic version of topological $\KTh$-theory to the Banach algebraic version:
\[
\KTh^{\top2}_*\left(\gG,\ B\right)\to \KTh^{\top2,\ban}_*\left(\gG,\ B\right).
\]

The Bost assembly map is defined in analogy to the Baum-Connes assembly map:

\begin{definition}\label{Definition:TheBostMap} Let $B$ a $\gG$-Banach algebra. Define the \demph{Bost assembly map} as the homomorphism of abelian groups
\[
\mu_{\mA}^B\colon \KTh^{\top2,\ban}_0\left(\gG,\ B\right) \to \KTh_0\left(\mA\left(\gG,B\right)\right)
\]
which is the direct limit of the group homomorphisms $\mu_{\mA,X}^B$ given by
\[
\KKbanW{\gG}\left(\Cont_0(X),\ B\right) \stackrel{j_{\mA}}{\to}\KKban\left(\mA\left(\gG, \Cont_0(X)\right), \
\mA\left(\gG,B\right)\right) \stackrel{\Sigma(\cdot)\left(\lambda_{X,\gG,\mA}\right)}{\to} \KTh_0\left(\mA\left(\gG, B\right)\right)
\]
where $X$ runs through all closed, $\gG$-compact, proper subspaces of $\uEgG$.
\end{definition}

\noindent Here, $\lambda_{X,\gG, \mA}$ denotes a canonical element of $\KTh_0(\mA(\gG, \Cont_0(X)))$ and $\Sigma(\cdot)$ denotes the action of $\KKban$ on $\KTh$-theory. See \cite{Paravicini:07:Induction:arxiv} or \cite{Paravicini:07} and to some extend \cite{Lafforgue:06} for more details.

To define the Bost-assembly map also for higher $\KTh$-groups note that there is a the canonical homomorphism $\iota_B\colon \mA(\gG, SB) \to S\mA(\gG,B)$ for every $\gG$-Banach algebra $B$. We can define $\mu_{\mA}^B$ also for $\KTh^{\top2,\ban}_1\left(\gG,\ B\right)$ as the composition
\[
\xymatrix{
\KTh^{\top2,\ban}_1\left(\gG,\ B\right) = \KTh^{\top2,\ban}_0\left(\gG,\ SB\right) \ar[r]^-{\mu_{\mA}^B} & \KTh_0(\mA(\gG,SB)) \ar[r]^-{\iota_{B,*}} & \KTh_0(S\mA(\gG,B)) = \KTh_1(\mA(\gG,B)).
}
\]
Proceed inductively to define the assembly map for all $n\in \N_0$. Note that $\iota_B$ is an isomorphism in $\KTh$-theory by Corollary~\ref{Corollary:UnconditionalAndSuspension:Groupoid}.

The \emph{Banach algebraic version of the Bost conjecture} for $\gG$ and $\mA(\gG)$ with coefficients in a $\gG$-Banach algebra $B$ asserts that $\mu^B_{\mA}$ is an isomorphism for all $n\in \N_0$. If $B$ is a $\gG$-C$^*$-algebra, then the assembly map introduced in \cite{Lafforgue:06} factors through the Banach algebraic topological $\KTh$-theory, i.e., it is given by the composition
\[
\xymatrix{
\KTh^{\top2}_*(\gG,B) \ar[r] &  \KTh^{\top2, \ban}_*(\gG,B) \ar[r]^{\mu^B_{\mA}} & \KTh_*(\mA(\gG, B)).
}
\]
The \emph{C$^*$-algebraic version of the Bost conjecture} for $\gG$ and $\mA(\gG)$ with coefficients in a $\gG$-C$^*$-algebra $B$ asserts that the composed assembly map is an isomorphism. The C$^*$-algebraic version of the Bost conjecture for groups is an instance of an ``isomorphism conjecture'' as elaborated in \cite{BarEchLueck:07}.

\subsection{The Bost conjecture and proper groupoids}

\begin{definition}[Hereditary subalgebra]\label{Definition:HereditarySubalgebra} Let $B_0$ be a subalgebra of a complex algebra $B$. Then $B$ is called \demph{hereditary} if $B_0 \ B\ B_0 \subseteq B_0$.
\end{definition}

\noindent The following lemma is a variant of Lemme~1.7.9 of \cite{Lafforgue:02}, a proof can be found in \cite{Paravicini:07}, Lemma~8.2.2.

\begin{lemma}\label{LemmaHereditaryImageAndNilpotentKernel} Let $B$ be a Banach algebra and let $A$ be a topological algebra (with separately continuous multiplication) and let $\varphi \colon A\to B$ be a continuous homomorphism such that $\varphi(A)$ is a dense hereditary subalgebra of $B$ and such that the kernel of $\varphi$ is nilpotent. Then $\varphi\colon \pi_0\big(\unital{A}^{-1}\big) \to \pi_0\big(\unital{B}^{-1}\big)$ is a bijection.
\end{lemma}

\noindent The following lemma is an elaborate version of Lemme~1.7.10 of \cite{Lafforgue:02}; there are two minor differences: The first is that we allow $\norm{\cdot}_1$ and $\norm{\cdot}_2$ to be semi-norms rather than norms (with the restriction that the kernel of the homomorphisms into the completions are nilpotent), and secondly, we do not ask the homomorphism $\psi$ to be injective. The first generalisation is necessary because we want to apply the result to unconditional completions in the groupoid setting where semi-norms appear naturally, the second generalisation might already be necessary in the setting of \cite{Lafforgue:02}, because in the proof of Lemme~1.7.8 there is no explicit argument given why the homomorphism from $\mB(G,B)$ to $\mA(G,B)$ is injective.

The lemma is proved analogously to Lemme~1.7.10 of \cite{Lafforgue:02}, based on our Lemma~\ref{LemmaHereditaryImageAndNilpotentKernel}.

\begin{lemma}\label{Lemma:HereditraySubalgebraNilpotentKernelTheSameKTheory}
Let $A$ be a topological algebra (with separately continuous multiplication). Let $\norm{\cdot}_1$ and $\norm{\cdot}_2$ be continuous semi-norms on $A$ such that the completion of $A$ with respect to both norms is a Banach algebra. Let $\iota_1$ be the canonical continuous homomorphism from $A$ into its completion $B_1$ with respect to $\norm{\cdot}_1$ and define $\iota_2$ and $B_2$ analogously. Assume that $\norm{a}_1 \geq \norm{a}_2$ for all $a\in A$, and let $\psi\colon B_1 \to B_2$ the homomorphism of Banach algebras that we get from this inequality. Assume also that $\iota_i(A)$ is hereditary in $B_i$ and that the kernel of $\iota_i$ is nilpotent for all $i\in \{1,2\}$. Then the map
\[
\psi_* \colon \KTh_*(B_1) \to \KTh_*(B_2)
\]
is an isomorphism.
\end{lemma}

\noindent \emph{For the rest of this section, let $\gG$ be proper and let $B$ be a non-degenerate $\gG$-Banach algebra.}

\begin{lemma}\label{Lemma:CompactSectionHereditaryInCompletion} Let $\mA(\gG)$ be \emph{regular}. Let $\iota$ be the canonical map from $\ContSect_c(\gG,\ r^*B)$ to $\mA(\gG,B)$. Since $\gG$ is proper, $\iota\left(\ContSect_c(\gG,\ r^*B)\right)$ is a hereditary subalgebra of $\mA(\gG,B)$ and the kernel $N$ of $\iota$ satisfies $\ContSect_c(\gG,\ r^*B) \ N \ \ContSect_c(\gG,\ r^*B) =0$; in particular, it is nilpotent with $N^3=0$.
\end{lemma}
\begin{proof}
Let $\mA(\gG)$ act on the equivariant pair $\mH(\gG)$ of locally convex monotone completions of $\Cont_c(\gG)$. Let $\beta^<, \beta^> \in \ContSect_c\left(\gG, \ r^*B\right)$. Let $K_r:= r\left(\supp \beta^<\right)$ and $K_s:= s\left(\supp \beta^>\right)$. The two sets $K_r$ and $K_s$ are compact subsets of $\gG^{(0)}$. Because $\gG$ is proper, the set $K:=\left\{\gamma \in \gG:\ r(\gamma)\in K_r,\ s(\gamma) \in K_s\right\}$ is compact. For all $\beta\in \ContSect_c\left(\gG,\ r^*B\right)$, we have $\supp \left(\beta^< * \beta* \beta^>\right) \subseteq K$. Because $\mA(\gG)$ acts on $\mH(\gG)$, we also have, by \ref{Proposition:ExtensionOfBracketToConvolution} and Property~(H2):
\[
\norm{\beta^< * \beta *\beta^>}_{\infty} \leq \norm{\beta^<}_{\mH^<} \norm{\beta}_{\mA} \norm{\beta^>}_{\mH^>}.
\]
It follows that $\left(\beta^< * \beta_n * \beta^>\right)_{n\in \N}$ is a Cauchy-sequence in $\ContSect_K\left(\gG, \ r^*B\right)$ whenever $\left(\beta_n\right)_{n\in \N}$ is a Cauchy-sequence in $\ContSect_c\left(\gG,\ r^*B\right)$ for the semi-norm $\norm{\cdot}_{\mA}$; in this case, $\left(\beta^< * \beta_n * \beta^>\right)_{n\in \N}$ converges to some element of $\ContSect_K\left(\gG,\ r^*B\right)$, and hence $\iota\left(\beta^< * \beta_n * \beta^>\right)= \iota(\beta^<) \iota(\beta_n) \iota(\beta^>)$ converges to some element in the image of $\iota$ if $n\to \infty$. Thus the image of $\iota$ is hereditary in $\mA(\gG,B)$.

Now let $\beta \in \ContSect_c\left(\gG,\ r^*B\right)$ satisfy $\iota(\beta)=0\in\mA(\gG,B)$. Let $\beta^<,\beta^>$ be elements of $\ContSect_c\left(\gG,\ r^*B\right)$. From (H2) we have $\norm{\beta^< * \beta * \beta^>}_{\infty} \leq \norm{\beta^<}_{\mH^<} \norm{\beta}_{\mA} \norm{\beta^>}_{\mH^>} = 0$, so $\beta^< * \beta * \beta^> =0$. This shows that the kernel $N$ of $\iota$ satisfies $\ContSect_c(\gG,\ r^*B) \ N \ \ContSect_c(\gG,\ r^*B) =0$.
\end{proof}

\noindent As a consequence of the preceding lemmas, the $\KTh$-theory of $\mA(\gG,B)$  does not depend on the particular (regular) completion $\mA(\gG)$ if $\gG$ is proper:

\begin{proposition}\label{Proposition:ProperGroupoidAlgebraKTheory}
Let $\mA'(\gG)$ be another regular unconditional completion of $\Cont_c(\gG)$. Then $\KTh_*\left(\mA(\gG,B)\right)$ and $\KTh_*\left(\mA'(\gG,B)\right)$ are canonically isomorphic.
\end{proposition}
\begin{proof}
Define  $\norm{\chi}_{\mB}:= \max\left\{\norm{\chi}_{\mA},\ \norm{\chi}_{\mA'}\right\}$ for all $\chi\in \Cont_c(\gG)$. Then $\mB(\gG)$ is a regular unconditional completion of $\Cont_c(\gG)$. By the preceding lemmas it follows that $\KTh_*\left(\mB(\gG,B)\right) \cong \KTh_*\left(\mA(\gG,B)\right)$ and $\KTh_*\left(\mB(\gG,B)\right) \cong \KTh_*\left(\mA'(\gG,B)\right)$. The resulting isomorphism $\KTh_*\left(\mA(\gG,B)\right) \cong \KTh_*\left(\mA'(\gG,B)\right)$ does not depend on the particular norm $\norm{\cdot}_{\mB}$, we could have taken any unconditional norm dominating $\norm{\cdot}_{\mA}$ and $\norm{\cdot}_{\mA'}$.
\end{proof}

\begin{Xample}\label{Example:TwoCompletionsOfGroupActionGroupoidSameKtheory}
Let $G$ be a locally compact Hausdorff group acting properly on some locally compact Hausdorff space $X$. Then
$\Leb^1\left(G\ltimes X\right)$ and $\Leb^1 \left(G,\ \Cont_0(X)\right)$ are two regular unconditional completions of $\Cont_c\left(G\ltimes X\right)$. Because $G\ltimes X$ is a proper groupoid, we have a canonical isomorphism
\[
\KTh_*\left(\Leb^1\left(G,\ \Cont_0(X)\right)\right) \cong \KTh_*\left(\Leb^1\left(G\ltimes X\right)\right).
\]
Because the unconditional norm given by $\Leb^1\left(G,\ \Cont_0(X)\right)$ dominates $\norm{\cdot}_1$, the isomorphism in $\KTh$-theory is given by the canonical homomorphism from $\Leb^1\left(G,\ \Cont_0(X)\right)$ to $\Leb^1\left(G\ltimes X\right)$.
\end{Xample}
\noindent Because $\gG$ is proper, the proper $\gG$-space $X=\gG^{(0)}$ is a model for $\uEgG$. If, in addition, $X/\gG$ is compact, then the canonical homomorphism
\[
\KKbanW{\gG}\left(\Cont_0(X),\ B\right) \to \KTh^{\top2,\ban}_0\left(\gG,\ B\right)
\]
is an isomorphism; moreover, the following diagram commutes:
\begin{equation}\label{Diagram:BostForProperCocompactGroupoids}
\xymatrix{
\KKbanW{\gG}\left(\Cont_0(X),\ B\right)\ar[rr]^-{\GreenJulgAbstiegKKban{\mA}{B}}\ar[d]_{\cong}&& \RKKban\left(\Cont_0\left(X/\gG\right); \ \Cont_0\left(X/\gG\right),\ \mA(\gG,B)\right) \ar[d]^{\cong}\\
 \KTh^{\top2,\ban}_0\left(\gG,\ B\right) \ar[rr]^{\mu_{\mA}^B}&& \KTh_0\left(\mA(\gG,B)\right)
}
\end{equation}
\noindent The isomorphism on the right-hand side is the given by the embedding $\C\mapsto\Cont_0\left(X/\gG\right)$ as constant functions (compare Corollary~\ref{Corollary:RKKAndKKforCompactBaseSpace}). Actually, Diagram~(\ref{Diagram:BostForProperCocompactGroupoids}) commutes already on the level of ($\gG$-equivariant) $\KKban$-cycles up to isomorphism.

Applying Theorem~\ref{Theorem:NaturalSplitPartOfGreenJulg} to Diagram~(\ref{Diagram:BostForProperCocompactGroupoids}) yields:

\begin{lemma}\label{Lemma:gGProperAndQuotientCompactBostSurjective}
Let $\gG$ be proper, let $X/\gG$ be compact, let $\mA(\gG)$ be regular and let $B$ be non-degenerate. Then the Bost map $\mu_{\mA}^B$ has a natural split.
\end{lemma}

Note that Diagram~(\ref{Diagram:BostForProperCocompactGroupoids}) gives the result directly only for $*=0$, but we can apply it to $S^nB$ instead of $B$ and use Corollary~\ref{Corollary:UnconditionalAndSuspension:Groupoid} to obtain it for arbitrary degrees.

As the top arrow in the above diagram is not only surjective but bijective if we impose some technical extra condition on $\mA(\gG)$, see Theorem~\ref{Theorem:GeneralisedGreenJulg}, there is also a version of the preceding lemma that asserts that the Bost map is an isomorphism in this case. However, that extended version does not seem to be of great value in the discussion of the Bost conjecture for proper Banach algebras.

\subsection{The Bost conjecture and proper Banach algebras}

\subsubsection{The forgetful map}\label{Subsubsection:ThePushforwardConstruction}

We give a short summary of a pushforward construction introduced in \cite{Paravicini:07}, Section~8.3; see also \cite{Paravicini:07:Trees:arxiv}, Section~1. Let $Y$ be a locally compact Hausdorff left $\gG$-space with anchor map $\rho$  and let $B$ be a $\gG \ltimes Y$-Banach algebra. We now turn $B$ into a $\gG$-Banach algebra; because $B$ is a field of Banach algebras over $Y$, we have to merge all those fibres of $B$ over points of $Y$ that have the same image under $\rho$ to get a field of Banach algebras over $X$:

For all $x\in X$, define\footnote{This definition makes sense if $x\in \rho(Y)$, and can and should be interpreted as $\rho_*(B)_x=0$ if $x\notin \rho(Y)$.}
\[
(\rho_*B)_x:=\ContSect_0\left(Y_x,\ B\restr_{Y_x}\right).
\]
On this family $\rho_*B=((\rho_*B)_x)_{x\in X}$ of Banach algebras over $X$, one can define a structure of a u.s.c.~field of Banach spaces over $X$ such that $\left\{\rho_*(\xi):\ \xi \in \ContSect_0(Y,B)\right\} =\ContSect_0\left(X,\ \rho_*B\right)$, where $\rho_*(\xi)\colon x\mapsto \xi\restr_{Y_x}$ for all $\xi\in \ContSect_0(Y,B)$. Moreover, there is a canonical $\gG$-action on $\rho_*B$ making it a $\gG$-Banach algebra. If $B$ is non-degenerate, then so is $\rho_*B$.

If $B$ is a $\gG\ltimes Y$-Banach algebra, similar definitions can be made for $\gG\ltimes Y$-Banach $B$-pairs to obtain $\gG$-Banach $\rho_*B$-pairs. We call this construction the pushforward construction or forgetful map; the idea is that $\rho_*$ forgets the fine fibration over $Y$ and only remembers the coarser fibration over $X$.

The forgetful map lifts to $\KKban$-cycles, i.e., if $A$ and $B$ be $\gG \ltimes Y$-Banach algebras, then $\rho_*$ gives a homomorphism
\[
\rho_* \colon \KKbanW{\gG \ltimes Y} \left(A,B\right) \to \KKbanW{\gG} \left(\rho_*A,\ \rho_*B\right).
\]
This construction also induces a homomorphism $\KTh^{\top2,\ban}(\gG \ltimes Y, B)$ to $\KTh^{\top2,\ban}(\gG , \rho_* B)$.

The forgetful map is also compatible with the descent. To be more precise, let $\mA(\gG)$ be an unconditional completion of $\Cont_c(\gG)$. For all $\xi \in \Cont_c(\gG\ltimes Y)$, define
\[
\norm{\xi}_{\mA_Y}:= \norm{\gamma\mapsto \sup_{y\in Y_{r_{\gG}(\gamma)}} \abs{\xi(\gamma,y)} }_{\mA}.
\]
This is an unconditional norm on $\Cont_c(\gG\ltimes Y)$. We have $\mA_Y(\gG\ltimes Y, B) \cong \mA(\gG,\rho_*B)$ for all $\gG\ltimes Y$-Banach algebras $B$. In \cite{Paravicini:07:Trees:arxiv}, it is shown that the following diagram is commutative
\begin{equation}\label{Diagram:DescentAndForgetfulMap}
\xymatrix{
\KTh^{\top2,\ban}_*(\gG \ltimes Y, B) \ar[d] \ar[r] & \KTh_*(\mA_Y(\gG \ltimes Y, B))\ar[d] ^{\cong}\\
\KTh^{\top2,\ban}_*(\gG , \rho_* B) \ar[r] & \KTh_*(\mA(\gG, \rho_* B)).
}
\end{equation}

\noindent Finally, the following result was shown in \cite{Paravicini:07}, the proof being somewhat technical (Proposition~8.3.26):

\begin{proposition}\label{Proposition:GoodCompletionVererbung} If $\mA(\gG)$ is regular, then also $\mA_Y(\gG\ltimes Y)$ is regular.
\end{proposition}

\subsubsection{Proper $\gG$-Banach algebras}

\begin{definition} A $\gG$-Banach algebra $B$ is called \demph{proper} if there is a proper locally compact Hausdorff $\gG$-space $Z$ (with anchor map $\rho$) and a $\gG \ltimes Z$-Banach algebra $\hat{B}$ such that the $\gG$-Banach algebra $\rho_* \hat{B}$ is isomorphic to $B$.
\end{definition}

\noindent As for proper C$^*$-algebras one can prove that we can assume without loss of generality that the space $Z$ is equal to $\uEgG$. Note that if $\gG$ itself is proper, then every $\gG$-Banach algebra is proper.

Let $\mA(\gG)$ be a regular unconditional completion of $\Cont_c(\gG)$ and let $B$ be a \emph{proper} non-degenerate $\gG$-Banach algebra. The following proposition generalises Proposition~\ref{Proposition:ProperGroupoidAlgebraKTheory}, which discusses the case that $\gG$ itself is proper. We are going to prove it by reducing it to this special case.

\begin{proposition}\label{Proposition:ProperAlgebraKTheory}
Let $\mA'(\gG)$ be another regular unconditional completion of $\Cont_c(\gG)$. Then $\KTh_*\left(\mA(\gG,B)\right)$ and $\KTh_*\left(\mA'(\gG,B)\right)$ are canonically isomorphic.
\end{proposition}
\begin{proof}
As in the proof of Proposition~\ref{Proposition:ProperGroupoidAlgebraKTheory}, let $\mB(\gG)$ be a regular uncondition completion of $\Cont_c(\gG)$ the norm of which dominates the norms of $\mA(\gG)$ and $\mA'(\gG)$. It suffices to compare $\mA(\gG)$ and $\mB(\gG)$. Let $\psi$ be the canonical homomorphism of Banach algebras from $\mB(\gG,B)$ to $\mA(\gG,B)$. We show that $\psi_*\colon \KTh_*\left(\mB(\gG,B)\right) \to \KTh_*\left(\mA(\gG,B)\right)$ is an isomorphism.

Find a proper locally compact Hausdorff $\gG$-space $Z$ with anchor map $\rho$ and a $\gG \ltimes Z$-Banach algebra $\hat{B}$ such that $\rho_*\hat{B}$ is isomorphic to $B$. Then $\hat{B}$ is non-degenerate. Because $\mA(\gG)$ and $\mB(\gG)$ are regular unconditional completions of $\Cont_c\left(\gG\right)$, also $\mA_Z\left(\gG\ltimes Z\right)$ and $\mB_Z\left(\gG\ltimes Z\right)$ are {regular} unconditional completions of $\Cont_c\left(\gG\ltimes Z\right)$ by Proposition~\ref{Proposition:GoodCompletionVererbung}. Moreover, $\norm{\chi}_{\mB_Z} \geq \norm{\chi}_{\mA_Z}$ for all $\chi\in \Cont_c\left(\gG\ltimes Z\right)$, hence there is a canonical homomorphism $\psi^Z\colon \mB_Z(\gG\ltimes Z,\ \hat{B}) \to \mA_Z(\gG\ltimes Z,\ \hat{B})$. The following diagram commutes
\[
\xymatrix{
\mB_Z\left(\gG \ltimes Z,\ \hat{B}\right)  \ar[rr]^{\psi^Z} \ar[d]_{\cong} && \mA_Z\left(\gG \ltimes Z,\ \hat{B}\right) \ar[d]^{\cong}\\
\mB\left(\gG, B\right)\ar[rr]^{\psi} && \mA\left(\gG, B\right)
}
\]
Hence also the following diagram commutes
\[
\xymatrix{
\KTh_*\left(\mB_Z\left(\gG \ltimes Z,\ \hat{B}\right)\right)  \ar[rr]^{\psi^Z_*} \ar[d]_{\cong} && \KTh_*\left(\mA_Z\left(\gG \ltimes Z,\ \hat{B}\right)\right) \ar[d]^{\cong}\\
\KTh_*\left(\mB\left(\gG, B\right)\right) \ar[rr]^{\psi_*} && \KTh_*\left(\mA\left(\gG, B\right)\right)
}
\]
By Proposition~\ref{Proposition:ProperGroupoidAlgebraKTheory}, $\psi^Z_*$ is an isomorphism, so $\psi_*$ is an isomorphism as well.
\end{proof}

\noindent Recall that $B$ is a non-degenerate proper $\gG$-Banach algebra and that $\mA(\gG)$ is a regular unconditional completion of $\Cont_c(\gG)$.

\begin{theorem}\label{Theorem:MainBostResult} The homomorphism
\[
\mu_{\mA}^B\colon \KTh^{\top2,\ban}_*\left(\gG,\ B\right) \to \KTh_*\left(\mA\left(\gG,B\right)\right)
\]
is split surjective. The split is natural in $B$.
\end{theorem}

\noindent This applies in particular to the regular unconditional completion $\Leb^1(\gG)$ and its symmetrised version $\Leb^1(\gG)\cap\Leb^1(\gG)^*$.

Before we prove Theorem~\ref{Theorem:MainBostResult}, we consider yet another special case:

\begin{lemma}\label{LemmaMainBostResultForGCompactProper} Let $B$ be a non-degenerate proper $\gG$-Banach algebra such that there exists a proper \emph{$\gG$-compact} $\gG$-space $Z$ with anchor map $\rho$ and a $\gG \ltimes Z$-Banach algebra $\hat{B}$ such that $\rho_* \hat{B} \cong B$. Then $\mu_{\mA}^B$ is split surjective, the split being natural in $B$.
\end{lemma}
\begin{proof}
Let $Z$, $\rho$ and $\hat{B}$ be as in the statement of the lemma. By Proposition~\ref{Proposition:GoodCompletionVererbung}, $\mA_Z\left(\gG\ltimes Z\right)$ is a regular unconditional completion of $\Cont_c\left(\gG \ltimes Z\right)$ because $\mA(\gG)$ is regular. So by Lemma~\ref{Lemma:gGProperAndQuotientCompactBostSurjective}, the homomorphism
\[
\mu_{\mA_Z}^{\hat{B}} \colon \KTh^{\top2,\ban}_*\left(\gG\ltimes Z,\ \hat{B}\right) \to \KTh_*\left(\mA_Z\left(\gG\ltimes Z,\ \hat{B}\right)\right)
\]
has a natural split. The diagram
\[
\xymatrix{
\KTh^{\top2,\ban}_*\left(\gG\ltimes Z,\ \hat{B}\right)  \ar[r]_{\mu_{\mA_Z}^{\hat{B}}} \ar[d]_{\rho_*} & \KTh_* \left(\mA_Z\left(\gG \ltimes Z,\ \hat{B}\right)\right)\ar[d]^{\cong} \ar @/_1.5pc/ @{.>} [l] \\
\KTh^{\top2,\ban}_*\left(\gG,\ B\right) \ar[r]_{\mu_{\mA}^{B}} & \ar @/_1.5pc/ @{.>} [l]\KTh_*\left(\mA\left(\gG, B\right)\right)
}
\]
commutes, see Diagram~(\ref{Diagram:DescentAndForgetfulMap}) above. Because the top-arrow has a natural split (dashed arrow), also the bottom-arrow has a natural split. \qedhere
\end{proof}

\begin{proof}[Proof of Theorem~\ref{Theorem:MainBostResult}]
Let $\hat{B}$ be a $\gG \ltimes \uEgG$-Banach algebra and let $\rho\colon \uEgG \to X=\gG^{(0)}$ be the anchor map of the proper action of $\gG$ on $\uEgG$; assume that $\rho_* \hat{B} \cong B$ as $\gG$-Banach algebras. Then $\hat{B}$ is non-degenerate. For every open $\gG$-invariant subspace $U$ of $\uEgG$, define $\hat{B}_U$ to be the $\gG \ltimes \uEgG$-Banach algebra with the following fibres: If $u\in U$, then the fibre over $u$ is $\hat{B}_u$, if $y\in \uEgG \setminus U$, then the fibre over $y$ is zero; the space $\ContSect(\uEgG,\ \hat{B}_U)$ is defined to be the set of all elements of $\ContSect(\uEgG,\ \hat{B})$ that vanish outside $U$. By definition, there is a $\gG \ltimes \uEgG$-equivariant ``injection'' $\hat{\jmath}_U$ from $\hat{B}_U$ to $\hat{B}$. It descends to a $\gG$-equivariant homomorphism $j_U:= \rho_* \hat{\jmath}_U$ from $B_U:= \rho_* \hat{B}_U$ to $B=\rho_* \hat{B}$. We can regard $B_U$ as a subalgebra of $B$.

The $\hat{B}_U$, where $U$ runs through the open $\gG$-invariant subsets of $\uEgG$ such that $\gG \backslash U$ is relatively compact, form a directed system: If $U$ and $V$ are open $\gG$-invariant and $\gG$-relatively compact subsets of $\uEgG$ with $U\subseteq V$, then there is an obvious homomorphism $\hat{\jmath}_{U,V}\colon \hat{B}_U \to \hat{B}_V$ such that $\hat{\jmath}_U = \hat{\jmath}_V \circ \hat{\jmath}_{U,V}$. Also the $B_U$ form a directed system, just take the $j_{U,V}:=\rho_*\hat{\jmath}_{U,V}$ as connecting maps. We can regard $B$ as the direct limit of the $B_U$. More importantly, the $\mA\left(\gG,B_U\right)$ form a directed system with connecting maps $\alpha_{U,V}:=\mA(\gG,\ j_{U,V})\colon \mA(\gG,B_U) \to \mA(\gG,B_V)$. The Banach algebra $\mA(\gG,B)$ is the direct limit of this system with embeddings $\alpha_U:=\mA(\gG,\ j_{U})\colon \mA(\gG,B_U) \to \mA(\gG,B)$. Because the $\KTh$-theory of Banach algebras is continuous, we get:
\vspace{-0.1cm}
\[
\KTh_*\left(\mA(\gG,B)\right) = \lim_{\to} \KTh_*\left(\mA(\gG, B_U)\right)
\vspace{-0.2cm}
\]
where $U$ runs through the $\gG$-invariant open subsets of $\uEgG$ such that $\gG\backslash U$ is relatively compact.

Now let $U$ be such a set. Find a closed set $Z\subseteq \uEgG$ such that $U\subseteq Z$ and $\gG \backslash Z$ is compact. Define $\rho^Z:= \rho\restr_Z$. Then $\hat{B}_U\restr_Z$ is a $\gG \ltimes Z$-Banach algebra and $(\rho^Z)_*\hat{B}_U \restr_Z$ is isomorphic to $B_U$. So $B_U$ satisfies the hypotheses of Lemma~\ref{LemmaMainBostResultForGCompactProper}, so $\mu_{\mA}^{B_U} \colon \KTh^{\top2,\ban}_*\left(\gG,\ B_U\right) \to \KTh_*\left(\mA\left(\gG,B_U\right)\right)$ is split surjective. Let $\sigma_U$ denote the natural split constructed above. It is easy to see, using the naturality of the split, that $\sigma_V \circ (\alpha_{U,V})_* = (j_{U,V})_*\circ \sigma_U$. Define $\tau_U := (j_U)_* \circ \sigma_U \colon \KTh_*\left(\mA\left(\gG, B_U\right)\right) \to \KTh^{\top2,\ban}_*\left(\gG,\ B\right)$. Then $\tau_V = \tau_U \circ (\alpha_{U,V})_*$. The universal property of the direct limit shows that there exists a natural homomorphism $\tau\colon \KTh_0\left(\mA\left(\gG, B\right)\right) \to \KTh^{\top2,\ban}\left(\gG,\ B\right)$ such that $\tau \circ (\alpha_U)_* = \tau_U$ for all $U$.

Note that
\vspace{-0.1cm}
\[
\mu_{\mA}^B \circ \tau_U = \mu_{\mA}^B \circ (j_U)_* \circ \sigma_U =  (\alpha_U)_* \circ \mu_{\mA}^{B_U} \circ \sigma_U = (\alpha_U)_*
\]
because $\sigma_U$ is a split. Passing to the limit shows that $\mu_{\mA}^B\circ \tau = \id$, i.e., $\tau$ is a natural split.
\end{proof}

\subsubsection{The case of locally compact groups}

Let $G$ be a locally compact Hausdorff group. A \emph{proper $G$-Banach algebra $B$} is a $G$-Banach algebra which carries an action of $\Cont_0(Z)$ for some locally compact Hausdorff proper $G$-space $Z$ such that the two actions are compatible, i.e., such that $B$ is a $G$-$\Cont_0(Z)$-Banach algebra; to fit the definition of a proper Banach algebra in the groupoid setting, we also demand $B$ to be \emph{locally $\Cont_0(Z)$-convex}.\footnote{See Paragraph~\ref{Subsection:ComparisonKKbanRKKban}.} In this case, $B$ can also be regarded as a $G\ltimes Z$-Banach algebra which we then call $\hat{B}$ to match the above notation. In this situation, we have the following corollary of Theorem~\ref{Theorem:MainBostResult} which we state for those who do not like groupoids:

\begin{corollary}
If $B$ is a proper $G$-Banach algebra and $\mA(G)$ is a regular unconditional completion of $\Cont_c(G)$, then
\[
\mu_{\mA}^B \colon \KTh^{\top2,\ban}_*\left(G,B\right) \to \KTh_*(\mA(G,B))
\]
is split surjective. In particular, this is true for $\mA(G) = \Leb^1(G)$.
\end{corollary}

\nocite{BaumConHig:94} \nocite{Blackadar:98} \nocite{Bost:90} \nocite{Kasparov:88} \nocite{KasSkan:03} \nocite{CEO:03} \nocite{CEO:03} \nocite{ChaEch:01:Permanence}


\begin{thebibliography}{CEOO03}

\bibitem[BCH94]{BaumConHig:94}
Paul Baum, Alain Connes, and Nigel Higson.
\newblock Classifying space for proper actions and {$K$}-theory of group
  {$C^*$}-algebras.
\newblock {\em Contemporary Math.}, 167:241--291, 1994.

\bibitem[BEL07]{BarEchLueck:07}
Arthur Bartels, Siegfried Echterhoff, and Wolfgang L\"{u}ck.
\newblock Inheritance of isomorphism conjectures under colimits.
\newblock {\em Preprintreihe SFB 478 - Geometrische Strukturen in der
  Mathematik}, 252, 2007.

\bibitem[Bla96]{Blanchard:96}
{\'E}tienne Blanchard.
\newblock D\'eformations de {$C\sp *$}-alg\`ebres de {H}opf.
\newblock {\em Bull.\ Soc.\ Math.\ France}, 124(1):141--215, 1996.

\bibitem[Bla98]{Blackadar:98}
Bruce Blackadar.
\newblock {\em {$K$}-theory for operator algebras}, volume~5 of {\em
  Mathematical Sciences Research Institute Publications}.
\newblock Cambridge University Press, Cambridge, second edition, 1998.

\bibitem[Bos90]{Bost:90}
Jean-Benoit Bost.
\newblock Principe d'{O}ka, {$K$}-th\'eorie et syst\`emes dynamiques non
  commutatifs.
\newblock {\em Invent. Math.}, 101(2):261--333, 1990.

\bibitem[CE01]{ChaEch:01:Permanence}
J{\'e}r{\^o}me Chabert and Siegfried Echterhoff.
\newblock Permanence properties of the {B}aum-{C}onnes conjecture.
\newblock {\em Doc.\ Math.}, 6:127--183 (electronic), 2001.

\bibitem[CEOO03]{CEO:03}
J{\'e}r{\^o}me Chabert, Siegfried Echterhoff, and Herv{\'e} Oyono-Oyono.
\newblock Shapiro's lemma for topological {$K$}-theory of groups.
\newblock {\em Comment. Math. Helv.}, 78(1):203--225, 2003.

\bibitem[CMR07]{CMR:07}
Joachim Cuntz, Ralf Meyer, and Jonathan~M. Rosenberg.
\newblock {\em Topological and bivariant {$K$}-theory}, volume~36 of {\em
  Oberwolfach Seminars}.
\newblock Birkh\"auser Verlag, Basel, 2007.

\bibitem[DG83]{DuGil:83}
Maurice~J. Dupr\'{e} and Richard~M. Gillette.
\newblock {\em Banach bundles, Banach modules and automorphisms of
  {$C^*$}-algebras}.
\newblock Pitman Books Limited, 1983.

\bibitem[Gie82]{Gierz:82}
Gerhard Gierz.
\newblock {\em Bundles of topological vector spaces and their duality}, volume
  955 of {\em Lecture Notes in Mathematics}.
\newblock Springer-Verlag, Berlin, 1982.
\newblock With an appendix by the author and Klaus Keimel, Queen's Papers in
  Pure and Applied Mathematics, 57.

\bibitem[Hof72]{Hofmann:72}
Karl~Heinrich Hofmann.
\newblock Representations of algebras by continuous sections.
\newblock {\em Bull. Amer. Math. Soc.}, 78:291--373, 1972.

\bibitem[Kas88]{Kasparov:88}
Gennadi~G. Kasparov.
\newblock Equivariant {$KK$}-theory and the {N}ovikov conjecture.
\newblock {\em Invent.\ Math.}, 91:147--201, 1988.

\bibitem[Kas95]{Kasparov:81:Conspectus}
Gennadi~G. Kasparov.
\newblock {$K$}-theory, group {$C\sp *$}-algebras, and higher signatures
  (conspectus).
\newblock In {\em Novikov conjectures, index theorems and rigidity, Vol.\ 1
  (Oberwolfach, 1993)}, volume 226 of {\em London Math. Soc. Lecture Note
  Ser.}, pages 101--146. Cambridge Univ. Press, Cambridge, 1995.

\bibitem[KS03]{KasSkan:03}
Gennadi~G. Kasparov and Georges Skandalis.
\newblock Groups acting properly on ``bolic'' spaces and the {N}ovikov
  conjecture.
\newblock {\em Ann. of Math. (2)}, 158(1):165--206, 2003.

\bibitem[Laf02]{Lafforgue:02}
Vincent Lafforgue.
\newblock {$K$}-th\'{e}orie bivariante pour les alg\`{e}bres de {B}anach et
  conjecture de {B}aum-{C}onnes.
\newblock {\em Invent.\ Math.}, 149:1--95, 2002.

\bibitem[Laf06]{Lafforgue:06}
Vincent Lafforgue.
\newblock {$K$}-th\'{e}orie bivariante pour les alg\`{e}bres de {B}anach,
  groupo\"{i}des et conjecture de {B}aum-{C}onnes. {A}vec un appendice
  d'{H}erv\'{e} {O}yono-{O}yono.
\newblock {\em J.\ Inst.\ Math.\ Jussieu}, 2006.
\newblock Published online by Cambridge University Press 28 Nov 2006.

\bibitem[Par07]{Paravicini:07}
Walther Paravicini.
\newblock {\em $\KK$-Theory for {B}anach Algebras And Proper Groupoids}.
\newblock PhD thesis, Universit\"{a}t M\"{u}nster, 2007.
\newblock Persistent identifier: urn:nbn:de:hbz:6-39599660289.

\bibitem[Par08a]{Paravicini:07:Morita:erschienen}
Walther Paravicini.
\newblock Morita equivalences and {KK}-theory for banach algebras.
\newblock {\em J. of the Inst. of Math. of Jussieu}, Forthcoming(-1):1--29,
  2008.

\bibitem[Par08b]{Paravicini:07:CNullX:erschienen}
Walther Paravicini.
\newblock A note on {B}anach {$\Cont_0(X)$}-modules.
\newblock {\em M{\"u}nster J. of Math.}, 1:267--278, 2008.
\newblock Persistent identifier: urn:nbn:de:hbz:6-43529451393.

\bibitem[Par09a]{Paravicini:07:Trees:arxiv}
Walther Paravicini.
\newblock The {B}ost conjecture, opens subgroups and groups acting on trees.
\newblock {\em http://arxiv.org/abs/0902.4339}, 2009.
\newblock Preprint.

\bibitem[Par09b]{Paravicini:07:Induction:arxiv}
Walther Paravicini.
\newblock Induction for {B}anach algebras, groupoids and {KK}$^{\text{ban}}$.
\newblock {\em http://arxiv.org/abs/0902.4199}, 2009.
\newblock Preprint.

\bibitem[Ren80]{Renault:80}
Jean~N. Renault.
\newblock {\em A Groupoid Approach to {C$^*$}-Algebras}, volume 793.
\newblock Springer-Verlag, Berlin, 1980.
\newblock Lecture Notes in Mathematics.

\bibitem[Tu99]{Tu:99}
Jean-Louis Tu.
\newblock La conjecture de {N}ovikov pour les feuilletages hyperboliques.
\newblock {\em $K$-Theory}, 16(2):129--184, 1999.

\bibitem[Tu04]{Tu:04}
Jean-Louis Tu.
\newblock Non-{H}ausdorff groupoids, proper actions and {$K$}-theory.
\newblock {\em Doc. Math.}, 9:565--597 (electronic), 2004.

\bibitem[Var74]{Varela:74}
Januario Varela.
\newblock Sectional representation of {B}anach modules.
\newblock {\em Math. Z.}, 139:55--61, 1974.

\end{thebibliography}
%

\end{document}